\documentclass[12pt]{amsart}
\usepackage{amssymb}
\usepackage{amsfonts}
\usepackage{amsmath}
\usepackage{mathrsfs}
\usepackage{amsthm, amsmath, amssymb, enumerate}
\usepackage{amssymb, enumitem}
\usepackage{bbm}
\usepackage{amscd}
\usepackage[all]{xy}
\usepackage[hidelinks]{hyperref}
\usepackage{amsmath}
\usepackage{amssymb}
\usepackage{amsthm}
\usepackage[sort&compress,numbers]{natbib}
\usepackage{tikz}
\usepackage{tikz-cd}
\usepackage{wrapfig}

\newtheorem{theorem}{Theorem}[section]
\newtheorem{proposition}[theorem]{Proposition}
\newtheorem{lemma}[theorem]{Lemma}
\newtheorem{corollary}[theorem]{Corollary}
\theoremstyle{definition}
\newtheorem*{claim*}{Claim}
\newtheorem{definition}[theorem]{Definition}

\newtheorem{remark}[theorem]{Remark}

\AtBeginDocument{   \def\MR#1{}}

\begin{document}
\title[Borel-definable homological algebra]{Applications of Borel-definable
homological algebra to locally compact groups}
\author{Martino Lupini}
\address{Dipartimento di Matematica, Universit\`{a} di Bologna, Piazza di
Porta S. Donato, 5, 40126 Bologna,\ Italy}
\email{martino.lupini@unibo.it}
\urladdr{http://www.lupini.org/}
\urladdr{}
\thanks{The author\ was partially supported by the Starting Grant 101077154
\textquotedblleft Definable Algebraic Topology\textquotedblright\ from the
European Research Council, the Gruppo Nazionale per le Strutture Algebriche,
Geometriche e le loro Applicazioni (GNSAGA) of the Istituto Nazionale di
Alta Matematica (INDAM), and the University of Bologna. }
\date{\today }

\begin{abstract}
We apply the description of the left heart of the category $\mathbf{LCPAb}$
of locally compact Polish abelian groups in terms of groups with a Polish
cover and Borel-definable group homomorphisms to study injective and
projective objects in the left heart of $\mathbf{LCPAb}$, as well as in the
left hearts of subcategories of $\mathbf{LCPAb}$. We also prove that, for a
countable discrete ring $R$, the group $\mathrm{Ext}_{R}^{1}(C,A)$ admits a
natural presentation as a group with a Polish cover whenever $C$ is a
locally compact Polish $R$-module and $A$ is a pro-Lie Polish $R$-module.
The presentation is given by Polish groups of module cocycles and
coboundaries and specializes, for $R=\mathbb{Z}$, to the cocycle
presentation for locally compact abelian groups. In particular, we prove
that the left hearts of the following categories have no nonzero injective
objects: locally compact Polish abelian groups, compactly generated locally
compact Polish abelian groups, locally compact Polish abelian groups of
finite ranks, abelian Lie groups, totally disconnected locally compact
Polish abelian groups, topological torsion locally compact Polish abelian
groups, and locally compact Polish abelian topological $p$-groups for any
prime $p$.
\end{abstract}

\maketitle

\section{Introduction}

From the very inception of the field of homological algebra, researchers
have endeavored to extend its methods from the setting of discrete algebraic
structures to the subtler setting of algebraic structures endowed with a
nontrivial topology. Early such efforts include the works of Iwasawa \cite%
{iwasawa_some_1949}, Gleason \cite%
{gleason_compact_1951,gleason_structure_1949,gleason_note_1949,gleason_groups_1952,gleason_structure_1951}%
, Calabi \cite{calabi_sur_1949,calabi_sur_1951}, Montgomery--Zippin \cite%
{montgomery_small_1952, montgomery_small_1952-1}, Yamabe \cite%
{yamabe_generalization_1953}, van Est \cite{van_est_algebraic_1955}, Mackey 
\cite{mackey_les_1957}, \cite{hochschild_cohomology_1953}, Rieffel \cite%
{rieffel_extensions_1966}, Moskowitz \cite{moskowitz_homological_1967}, Fulp 
\cite{fulp_splitting_1972,fulp_homological_1970}, Fulp--Griffith \cite%
{fulp_extensions_1971}, Brown \cite{brown_extensions_1971}, Moore \cite%
{moore_group_1976,moore_group_1976-1,moore_group_1968,moore_extensions_1964}%
, Wigner \cite%
{wigner_algebraic_1970,wigner_algebraic_1973,casselman_continuous_1974}, and
Kehlet \cite{kehlet_extensions_1979}. More recent work is due to Waelbroeck 
\cite%
{waelbroeck_quotient_1982,waelbroeck_taylor_1982,waelbroeck_quotient_1989},
Vasilescu \cite%
{vasilescu_spectral_1987,vasilescu_spectral_1989,vasilescu_spectral_1988},
Prosmans \cite%
{prosmans_derived_2000,prosmans_topological_2000,prosmans_derived_1999},
Loth \cite{loth_topologically_2001,loth_pure_2003,loth_pure_2006},
Hoffmann--Spitzweck \cite{hoffmann_homological_2007}, B\"{u}hler \cite%
{buhler_algebraic_2011}, Wegner \cite{wegner_heart_2017}, Rump \cite%
{rump_almost_2001,rump_ample_2021,rump_morita_2022}, Braunling \cite%
{braunling_theory_2020,braunling_theory_2021}, Cabello S\'{a}nchez--Navarro
Garmendia \cite{cabello_sanchez_category_2023}, and Cabello S\'{a}%
nchez--Castillo \cite{cabello_sanchez_homological_2023}. One of the main
goals of this series of works was the development of the algebraic invariant 
$\mathrm{Ext}$ for locally compact abelian groups, which parameterizes the
isomorphism classes of topological extensions of pairs of such groups.

A number of difficulties arise in the process of extending the $\mathrm{Ext}$
functor from discrete to topological settings. A first one is that, unlike
discrete abelian groups, locally compact abelian groups, for example, do not
form an \emph{abelian category}. A number of relaxations of the notion of
abelian category, generous enough to include categories of topological
algebraic structures, but also well-behaved enough to be suitable for the
development of homological algebra, have been considered so far. These
include: \emph{exact categories} in the sense of Buchsbaum \cite%
{buchsbaum_exact_1955,buchsbaum_note_1959,heller_homological_1958,walker_relative_1966}%
; \emph{pre-abelian categories} \cite{richman_ext_1977,richman_valuated_1979}
called \emph{quasi-abelian} by Yoneda in \cite%
{yoneda_note_1958,yoneda_homology_1954,yoneda_ext_1960}; \emph{Bourn
protomodular} categories \cite%
{borceux_malcev_2004,bourn_protomodularity_1998}; \emph{ideally exact
categories} \cite{janelidze_ideally_2024,janelidze_semidirect_2024}; \emph{%
semi-abelian categories} \cite%
{rump_almost_2001,borceux_survey_2004,gran_second_2008,janelidze_semi-abelian_2002,everaert_baer_2004,everaert_baer_2004-1,palamodov_homological_1971,goedecke_satellites_2009,wengenroth_raikov_2012,rump_analysis_2011,wengenroth_palamodov_2006}%
; \emph{relative semi-abelian categories} \cite%
{janelidze_relative_2009,goedecke_homology_2013}; \emph{categories} with 
\emph{proper classes} of morphisms \cite%
{eilenberg_foundations_1965,generalov_relative_1996,enochs_relative_2011,enochs_relative_2011-1}%
; \emph{homological categories} \cite{borceux_malcev_2004}; \emph{relative
homological categories} \cite{janelidze_relative_2006}; \emph{additive
regular categories} \cite{henrard_left_2022}; \emph{quasi-abelian categories}
\cite%
{gruson_completion_1966,succi_cruciani_sulle_1973,succi_cruciani_funtori_1974,succi_cruciani_sul_1975,beilinson_faisceaux_1982,lupini_looking_2024,savage_koszul_2023,tattar_torsion_2021,fiorot_quasi-abelian_2021,prosmans_derived_2000,schneiders_quasi-abelian_1999,rump_morita_2022,bondal_generators_2003,rump_analysis_2011,andre_slope_2009}
called \emph{almost abelian} in \cite{rump_modules_2001,rump_almost_2001}
and semi-abelian in \cite{raikov_semiabelian_1976,kuzminov_semiabelian_1972}%
; \emph{left quasi-abelian categories} \cite{rump_almost_2001,
rump_abelian_2020}; \emph{Morita categories} \cite%
{rump_morita_2022,rump_ample_2021}; \emph{Barr-exact categories} \cite%
{neeman_derived_1990,buhler_exact_2010,bodzenta_abelian_2024}; \emph{%
suspended categories} or $S$-\emph{categories} \cite{keller_chain_1990}.

In this paper we consider \emph{quasi-abelian categories }as in \cite%
{schneiders_quasi-abelian_1999} as a framework to study categories of
algebraic structures endowed with a topology. Locally compact abelian groups
form a quasi-abelian category, as do many other important categories of
algebraic structures endowed with a topology, such as Polish abelian groups
or Banach spaces. The monograph \cite{schneiders_quasi-abelian_1999}
presents a canonical way of enlarging a quasi-abelian category to an abelian
one. Such an \textquotedblleft abelian envelope\textquotedblright\ $\mathrm{%
LH}(\mathcal{A})$ of a quasi-abelian category $\mathcal{A}$ is called the 
\emph{left heart of $\mathcal{A}$}. More precisely, $\mathrm{LH}(\mathcal{A}%
) $ is the heart of the derived category of $\mathcal{A}$ with respect to
its canonical left truncation structure. A different \textquotedblleft
completion\textquotedblright\ of the category of locally compact
groups---its \emph{ample closure }as a Morita category---is considered in 
\cite{rump_ample_2021}.

A second difficulty arising in the development of homological invariants
such as $\mathrm{Ext}$ for locally compact abelian groups is the fact that
this category, in stark contrast to the category of discrete abelian groups,
does not have enough injective or projective objects. This problem was
partially resolved by Hoffmann and Spitzweck \cite{hoffmann_homological_2007}%
, who showed that one may in this context replace injective resolutions with
suitable \textquotedblleft approximately injective\textquotedblright\
resolutions. They thereby succeeded in formulating the functor $\mathrm{Ext}$
for locally compact abelian groups within the framework of derived
homological algebra, by identifying it as the cohomological right derived
functor of $\mathrm{Hom}$ \cite[Section 2.2, Definition 1.4]%
{verdier_categories_1977}.

The derived functor of $\mathrm{Hom}$ for locally compact abelian Polish
groups will then form our main tool, together with the concrete description
of the left heart of the category of abelian Polish groups given in \cite%
{lupini_looking_2024}. This description identifies the objects of this
category as the groups with an abelian Polish cover (previously considered
in \cite{bergfalk_definable_2024-1,bergfalk_definable_2024}) and the
morphisms of this category as the Borel-definable group homomorphisms
between them. By construing $\mathrm{Ext}$ as a functor to the category of
groups with a Polish cover, we are able to apply Borel-definable methods
developed in \cite%
{lupini_looking_2024,bergfalk_definable_2024-1,bergfalk_definable_2024} to
compute this invariant in a number of examples. As an application, we study
the injective and projective objects in the left heart of the category $%
\mathbf{LCPAb}$ of locally compact Polish abelian groups, as well as in the
left hearts of a number of thick subcategories of $\mathbf{LCPAb}$,
including those spanned by totally disconnected groups, topological torsion
groups, topological $p$-groups, finite-rank topological $p$-groups,
topological torsion groups of finite ranks, compactly generated groups, and
Lie groups. In the case of the category $\mathbf{TorFLCPAb}$ of topological
torsion groups of finite ranks, we characterize the objects of $\mathbf{%
TorFLCPAb}$ that are injective in its left heart and all the projective
objects of its left heart.

We also extend the cocycle description to modules over a countable discrete
unital ring $R$. An extension of $R$-modules is encoded not only by its
additive cocycle but also by the defect of a Borel section from commuting
with scalar multiplication. A term-saturated family of control measures is
used to make all scalar substitutions continuous on the relevant spaces of
almost-everywhere classes. For a locally compact Polish $R$-module $C$ and a
pro-Lie Polish $R$-module $A$, this gives $\mathrm{Ext}_{R}^{1}(C,A)$ as a
quotient of a Polish group of module cocycles by a Polishable group of
coboundaries.

In this paper all the topological groups, unless otherwise indicated, will
be understood to be Polish. Particularly, all countable groups will be
considered as topological groups endowed with the discrete topology. The
rest of this paper is divided into five sections. In Section \ref%
{Section:categories} we recall the requisite preliminaries from category
theory and abstract homological algebra. In Section \ref%
{Section:groups-with-a-Polish-cover} we recall the description of the left
heart of abelian Polish groups in terms of groups with a Polish cover, and
develop further Borel-definable methods for derived categories, focusing on
the category of locally compact abelian Polish groups. In\ Section \ref%
{Section:derived} we present equivalent descriptions of $\mathrm{Hom}$ and $%
\mathrm{Ext}$ as groups with a Polish cover in terms of cocycles, showing
that the classically interchangeable definitions of $\mathrm{Ext}$ remain so
in the present context as well. In Section~\ref{Section:module-cocycles} we
give the corresponding cocycle and coboundary description for extensions of
Polish modules over a countable discrete ring. Finally,\ Section \ref%
{Section:injective} contains the main results of this paper, leading to a
complete characterization of injective and projective objects in the left
heart of the category of locally compact Polish abelian groups and of many
of its important subcategories.

\subsection*{Acknowledgments}

We are grateful to Matteo Casarosa, Ivan Di Liberti, Luca Marchiori, Marco
Moraschini, Filippo Sarti, and Alessio Savini for a number of valuable
discussions. We also thank Jeffrey Bergfalk for a large number of useful
comments and remarks on a preliminary version of this paper. ChatGPT by
OpenAI (version GPT-5.6 Sol) and Claude by Anthropic (version Opus 5) have
been used in the final stages of the preparation of this manuscript for
proof-reading and proof-checking.

\section{Exact and quasi-abelian categories and derived functors}

\label{Section:categories}

{In this section we recall some fundamental notions from homological algebra
and category theory as can be found in \cite%
{kashiwara_categories_2006,gelfand_methods_2003,mac_lane_homology_1995}; see
also \cite[Section 2]{casarosa_phantom_2025} for the case of categories
(enriched over categories of) topological modules.}

\subsection{Exact categories}

An \emph{additive category} is a category $\mathcal{A}$ that is enriched
over the category $\mathbf{Ab}$ of abelian groups, and that has all finite
products \cite[Section VIII.2]{mac_lane_categories_1998}. In an additive
category, $n$-ary products and $n$-ary coproducts coincide. For $n=0$ one
obtains the notion of \emph{null object}, for $n=2$ they are called \emph{%
biproducts }and denoted by $A\oplus B$ for objects $A,B$ of $\mathcal{A}$. A
functor between additive categories is \emph{additive} if it is $\mathbf{Ab}$%
-enriched or, equivalently, it carries biproducts to biproducts. An additive
category is called \emph{idempotent complete }if every idempotent arrow $%
p:A\rightarrow A$ in $\mathcal{A}$ has a kernel, in which case there is a
decomposition $A\cong A_{0}\oplus A_{1}$ with respect to which $p$ is
isomorphic to the canonical projection on $A_{0}$.

A \emph{kernel-cokernel pair} in an additive category $\mathcal{A}$ is a
pair $\left( i,p\right) $ of arrows $i:A^{\prime }\rightarrow A$ and $%
p:A\rightarrow A^{\prime \prime }$ in $\mathcal{A}$ such that $i$ is a
kernel of $p$ and $p$ is a cokernel for $i$. An \emph{exact structure} on a
category $\mathcal{A}$ \cite[Section 2]{buhler_exact_2010} is a class $%
\mathcal{E}$ of kernel-cokernel pairs $\left( i,p\right) $, where $i$ is
called an \emph{admissible monic }or \emph{inflation} for $\mathcal{E}$ and $%
p$ is called an \emph{admissible epic }or \emph{deflation} for $\mathcal{E}$%
, that is closed under isomorphism and satisfies the following axioms, as
well as their dual axioms obtained by reversing all the arrows:

\begin{enumerate}
\item For all objects $A$, the identity morphism $1_{A}$ is an admissible
monic;

\item The class of admissible monics is closed under composition;

\item The pushout of an admissible monic along an arbitrary morphism exists
and yields an admissible monic.
\end{enumerate}

They also imply that the pullback of an admissible monic along an admissible
epic yields an admissible monic, as well as its dual \cite[Proposition 2.15]%
{buhler_exact_2010}. An \emph{exact category }$\left( \mathcal{A},\mathcal{E}%
\right) $ is an additive category $\mathcal{A}$ endowed with an exact
structure $\mathcal{E}$. The elements of $\mathcal{E}$ are called the\emph{\
short exact sequences} for the given exact category. We now recall the
notion of exact functor between exact categories as in \cite[Definition 5.1]%
{buhler_exact_2010}; see also \cite[Section 4.1]{previdi_locally_2011}.

\begin{definition}
\label{Definition:exact}An additive functor between exact categories is 
\emph{exact }if it maps short exact sequences to short exact sequences (and
hence it preserves pushouts along admissible monics and pullbacks along
admissible epics) \cite[Section 5]{buhler_exact_2010}.
\end{definition}

\begin{remark}
\label{Remark:exact}There is another notion of exact functor, defined in the
context of categories with finite limits. According to \cite[Definition
3.3.1 and Proposition 3.3.3]{kashiwara_categories_2006}, if $\mathcal{C}$
and $\mathcal{D}$ are categories with finite limits and colimits, then a
functor $F:\mathcal{C}\rightarrow \mathcal{D}$ is exact if it commutes with
finite limits and colimits. In the case when $\mathcal{C}$ and $\mathcal{D}$
are exact categories with finite limits and colimits, an exact functor $F:%
\mathcal{C}\rightarrow \mathcal{D}$ does not generally preserve finite
limits and colimits \cite[Definition 5.1]{buhler_exact_2010}. For example,
the forgetful functor $\mathbf{PAb}\rightarrow \mathbf{Ab}$ from the
category of Polish abelian groups to the category $\mathbf{Ab}$ of discrete
abelian groups is exact but does not preserve finite colimits, since it does
not preserve the cokernel of the inclusion $\mathbb{Q}\rightarrow \mathbb{R}$%
---which is trivial in $\mathbf{PAb}$ and nontrivial in $\mathbf{Ab}$. By
considering the opposite categories, one sees that also finite limits are
not necessarily preserved by exact functors in the sense of \cite[Definition
5.1]{buhler_exact_2010}.
\end{remark}

A morphism $f:A\rightarrow B$ in an exact category is \emph{admissible} if
it admits an \textquotedblleft admissible factorization\textquotedblright\ $%
f=m\circ e$ where $e$ is an admissible epic and $m$ is an admissible monic,
in which case such a decomposition is unique up to a unique isomorphism \cite%
[Section 8]{buhler_exact_2010}. The class of admissible morphisms is
invariant with respect to pushout along admissible monics and pullback along
admissible epics. A sequence 
\begin{equation*}
A^{\prime }\overset{f}{\rightarrow }A\overset{g}{\rightarrow }A^{\prime
\prime }
\end{equation*}%
of admissible morphisms with admissible factorizations $f=i\circ p$ and $%
g=j\circ q$ is \emph{acyclic} if the sequence%
\begin{equation*}
\bullet \overset{i}{\rightarrow }A\overset{q}{\rightarrow }\bullet
\end{equation*}%
is short exact.

\begin{definition}
\label{Definition:fully-exact}A\emph{\ fully exact subcategory }$\mathcal{B}$
of an exact category $\mathcal{A}$ is a full subcategory whose collection of
objects is essentially closed under extensions, in the sense that if 
\begin{equation*}
0\rightarrow B^{\prime }\rightarrow A\rightarrow B^{\prime \prime
}\rightarrow 0
\end{equation*}%
is a short-exact sequence with $B^{\prime }$ and $B^{\prime \prime }$ in $%
\mathcal{B}$, then $A$ is isomorphic to an object of $\mathcal{B}$ \cite[%
Definition 10.21]{buhler_exact_2010}. We then have that $\mathcal{B}$ is
itself an exact category, where a short exact sequence in $\mathcal{B}$ is a
short exact sequence in $\mathcal{A}$ whose objects are in $\mathcal{B}$ 
\cite[Lemma 10.20]{buhler_exact_2010}.
\end{definition}

\subsection{Triangulated category}

A \emph{triangulated category graded by the translation functor $T$} is an
additive category $\mathcal{A}$ endowed with an additive automorphism $T$
together with a collection of \emph{distinguished triangles }$\left(
X,Y,Z,u,v,w\right) $ where $u\in \mathrm{Hom}^{0}\left( X,Y\right) $, $v\in 
\mathrm{Hom}^{0}\left( Y,Z\right) $, and $w\in \mathrm{Hom}^{1}\left(
Z,X\right) $, satisfying the axioms of triangulated categories \cite[Section
1.1]{verdier_categories_1977}; see also \cite[Definition 10.1.6]%
{kashiwara_categories_2006}. In this definition, we let $\mathrm{Hom}%
^{k}\left( X,Y\right) $ be the group $\mathrm{Hom}_{\mathcal{A}}\left(
X,T^{k}Y\right) $. The composition of $f\in \mathrm{Hom}^{n}\left(
X,Y\right) $ and $g\in \mathrm{Hom}^{m}\left( Y,Z\right) $ is defined to be $%
g\circ f:=T^{n}g\circ f\in \mathrm{Hom}^{n+m}\left( X,Z\right) $. This
defines a new category whose sets of morphisms are \emph{graded groups}.

A functor $F:\mathcal{A}\rightarrow \mathcal{B}$ between triangulated
categories is \emph{triangulated }if it is additive, it is graded (which
means that there is a natural isomorphism $\alpha _{F}:TF\Rightarrow FT$),
and it maps distinguished triangles to distinguished triangles. A \emph{%
morphism} $\mu :F\Rightarrow G$ between triangulated functors $\mathcal{A}%
\rightarrow \mathcal{B}$ is a natural transformation $\mu $ such that $\mu
T\circ \alpha _{F}\cong \alpha _{G}\circ T\mu $.

An additive functor $F:\mathcal{A}\rightarrow \mathcal{E}$ from a
triangulated category $\mathcal{A}$ to an abelian category $\mathcal{E}$ is 
\emph{cohomological} \cite[Section 1.1, Definition 3.1]%
{verdier_categories_1977} if for every distinguished triangle $\left(
X,Y,Z,u,v,w\right) $ of $\mathcal{A}$, the sequence%
\begin{equation*}
FX\overset{Fu}{\rightarrow }FY\overset{Fv}{\rightarrow }FZ
\end{equation*}%
is exact. This yields a long exact sequence%
\begin{equation*}
\cdots \rightarrow F^{d}X\overset{F^{d}u}{\rightarrow }F^{d}Y\overset{F^{d}v}%
{\rightarrow }F^{d}Z\overset{F^{d+1}u}{\rightarrow }F^{d+1}X\rightarrow
\cdots
\end{equation*}%
where we set $F^{d}:=F\circ T^{d}$ for $d\in \mathbb{Z}$.

\subsection{Derived categories}

Suppose that $\left( \mathcal{A},\mathcal{E}\right) $ is an idempotent
complete exact category. We let $\mathrm{Ch}\left( \mathcal{A}\right) $ be
the category of (cochain) complexes over $\mathcal{A}$, and $\mathrm{K}%
\left( \mathcal{A}\right) $ be the category of complexes of $\mathcal{A}$
whose morphisms are the homotopy classes of morphisms of complexes \cite[%
Section 9]{buhler_exact_2010}. Then $\mathrm{K}\left( \mathcal{A}\right) $
is a triangulated category with translation functor $T$ defined by $\left(
TA\right) ^{n}=A^{n+1}$ for $n\in \mathbb{Z}$. A \emph{strict triangle }in $%
\mathrm{K}\left( \mathcal{A}\right) $ over a morphism $f:A\rightarrow B$ in $%
\mathrm{Ch}\left( \mathcal{A}\right) $ is one of the form $\left( A,B,%
\mathrm{cone}\left( f\right) ,f,i,j\right) $ where $i:B\rightarrow \mathrm{%
cone}\left( f\right) $ and $j:\mathrm{cone}\left( f\right) \rightarrow TA$
are the canonical morphisms of complexes; see \cite[Exercise 1.2.8]%
{weibel_introduction_1995}. A distinguished triangle in $\mathrm{K}\left( 
\mathcal{A}\right) $ is one that is isomorphic to a strict triangle.

We say that a complex $A$ over $\mathcal{A}$ is \emph{left bounded }%
(respectively, \emph{right bounded}) if there exists $n\in \mathbb{Z}$ such
that $A^{k}=0$ for every $k<n$ (respectively, for every $k>n$). We say that $%
A$ is bounded if it is both left bounded and right bounded.

A complex $A$ over $\mathcal{A}$ is \emph{acyclic }if for every $n\in 
\mathbb{Z}$,%
\begin{equation*}
A^{n}\rightarrow A^{n+1}\rightarrow A^{n+2}
\end{equation*}%
is an exact sequence of admissible morphisms. The subcategory $\mathrm{N}%
\left( \mathcal{A}\right) $ of $\mathrm{K}\left( \mathcal{A}\right) $
spanned by acyclic complexes is a \emph{triangulated subcategory} which is
closed under taking direct summands \cite[Corollary 10.11]{buhler_exact_2010}%
; see also \cite{neeman_triangulated_2001}. By definition, a morphism $%
u=[f]:A\rightarrow B$ is a quasi-isomorphism if and only if its cone is
acyclic. Equivalently, in a distinguished triangle $A\overset{u}{\rightarrow}%
B\rightarrow Z\rightarrow A[1]$, the third object $Z$ is acyclic; see \cite[%
Definition 10.16 and Remark 10.17]{buhler_exact_2010}.

The \emph{derived category }of $\mathcal{A}$ is by definition the \emph{%
Verdier quotient} 
\begin{equation*}
\mathrm{D}\left( \mathcal{A}\right) :=\frac{\mathrm{K}\left( \mathcal{A}%
\right) }{\mathrm{N}\left( \mathcal{A}\right) }\text{.}
\end{equation*}%
This is defined as the \emph{category of fractions }$\mathrm{K}\left( 
\mathcal{A}\right) [\Sigma ^{-1}]$, where $\Sigma $ is the saturated \emph{%
multiplicative system} of \emph{quasi-isomorphisms} in $\mathrm{K}\left( 
\mathcal{A}\right) $; see \cite[Definition 7.1.5]{kashiwara_categories_2006}
and \cite[Definition 7.1.19]{kashiwara_categories_2006}. We let $Q_{\mathcal{%
A}}:\mathrm{K}\left( \mathcal{A}\right) \rightarrow \mathrm{D}\left( 
\mathcal{A}\right) $ be the canonical quotient functor. One similarly
defines the categories $\mathrm{D}^{+}\left( \mathcal{A}\right) $, $\mathrm{D%
}^{-}\left( \mathcal{A}\right) $, and $\mathrm{D}^{b}\left( \mathcal{A}%
\right) $ by considering left bounded, right bounded, and bounded complexes
respectively.

\subsection{Quasi-abelian categories}

A \emph{quasi-abelian category} is an additive category $\mathcal{A}$ such
that each morphism has a kernel and a cokernel, and such that the class of
kernels is invariant under push-out along arbitrary morphisms, and dually
the class of cokernels is invariant under pull-back along arbitrary
morphisms \cite[Section 4]{buhler_exact_2010}; see also \cite%
{schneiders_quasi-abelian_1999}. Every quasi-abelian category $\mathcal{A}$
can be regarded as an exact category, whose short exact sequences are all
the kernel-cokernel pairs \cite[Proposition 4.4]{buhler_exact_2010}. The 
\emph{admissible monics }in $\mathcal{A}$ are the kernels, the \emph{%
admissible epics }are the cokernels. The \emph{admissible morphisms }are the
morphisms $f$ whose canonical induced arrow $\widehat{f}:\mathrm{Coim}\left(
f\right) \rightarrow \mathrm{Im}\left( f\right) $ is an isomorphism. In the
following, we will regard a quasi-abelian category as an exact category with
respect to such an exact structure. An \emph{abelian category }is a
quasi-abelian category in which every morphism is admissible.

Let $\mathrm{LK}\left( \mathcal{A}\right) $ and $\mathrm{LH}\left( \mathcal{A%
}\right) $ be the full subcategories of $\mathrm{K}\left( \mathcal{A}\right) 
$ and $\mathrm{D}\left( \mathcal{A}\right) $, respectively, spanned by
complexes $A$ such that $A^{n}=0$ for $n\in \mathbb{Z}\setminus \left\{
-1,0\right\} $, and such that $\delta _{A}^{-1}:A^{-1}\rightarrow A^{0}$ is
monic. Then the collection $\mathrm{N}\left( \mathcal{A}\right) \cap \mathrm{%
LK}\left( \mathcal{A}\right) $ of quasi-isomorphisms between objects of $%
\mathrm{LK}\left( \mathcal{A}\right) $ is a saturated multiplicative system,
and the quotient $\mathrm{LK}\left( \mathcal{A}\right) /\left( \mathrm{N}%
\left( \mathcal{A}\right) \cap \mathrm{LK}\left( \mathcal{A}\right) \right) $
is equivalent to $\mathrm{LH}\left( \mathcal{A}\right) $ \cite[Corollary
1.2.21]{schneiders_quasi-abelian_1999}. We have that $\mathrm{LH}\left( 
\mathcal{A}\right) $ is an abelian category, called the \emph{left heart }of 
$\mathcal{A}$ (or, more precisely, the heart of the derived category of $%
\mathcal{A}$ with respect to its canonical left truncation-structure). The
inclusion $\mathcal{A}\rightarrow \mathrm{LH}\left( \mathcal{A}\right) $
obtained by identifying an object of $\mathcal{A}$ with a complex
concentrated in degree zero is a finitely continuous, exact (in the sense of
Definition \ref{Definition:exact}) fully faithful functor that satisfies the
following universal property: for every abelian category $\mathcal{M}$ and
finitely continuous exact functor $F:\mathcal{A}\rightarrow \mathcal{M}$
there exists an essentially unique finitely continuous exact functor $G:%
\mathrm{LH}\left( \mathcal{A}\right) \rightarrow \mathcal{M}$ whose
restriction to $\mathcal{A}$ is isomorphic to $F$; see \cite[Proposition
1.2.34]{schneiders_quasi-abelian_1999}. Furthermore, $\mathcal{A} $ is a
reflective subcategory of $\mathrm{LH}\left( \mathcal{A}\right) $ \cite[%
Proposition 1.2.27]{schneiders_quasi-abelian_1999}.

\begin{remark}
Notice that an exact functor in the sense of Definition \ref%
{Definition:exact} is not necessarily finitely continuous; see Remark \ref%
{Remark:exact}.
\end{remark}

\begin{definition}
Let $\mathcal{A}$ be a quasi-abelian category, and $\mathcal{B}$ be a (not
necessarily full) subcategory of $\mathcal{A}$. Then we say that $\mathcal{B}
$ is a quasi-abelian subcategory of $\mathcal{A}$ if $\mathcal{B}$ is a
quasi-abelian category and the inclusion $\mathcal{B}\rightarrow \mathcal{A}$
is finitely continuous and finitely cocontinuous. We say that $\mathcal{B}$
is a \emph{fully exact quasi-abelian} \emph{subcategory} of $\mathcal{A}$ if 
$\mathcal{B}$ is a quasi-abelian subcategory of $\mathcal{A}$ as well as a
fully exact subcategory of $\mathcal{A}$ in the sense of Definition \ref%
{Definition:fully-exact}.
\end{definition}

The notion of \emph{fully exact abelian} \emph{subcategory} of an abelian
category is obtained as a particular case. Such a category is called \emph{%
thick }in \cite[Definition 8.3.21(iv)]{kashiwara_categories_2006}. In
analogy with the abelian case, we call a fully exact quasi-abelian
subcategory of a quasi-abelian category a \emph{thick }subcategory.

\subsection{Derived functors}

Suppose that $F:\mathcal{A}\rightarrow \mathcal{B}$ is a functor between
exact categories. Then $F$ induces a triangulated functor $\mathrm{K}\left( 
\mathcal{A}\right) \rightarrow \mathrm{K}\left( \mathcal{B}\right) $, which
we also denote by $F$. If $F$ is exact, then its extension $\mathrm{K}\left( 
\mathcal{A}\right) \rightarrow \mathrm{K}\left( \mathcal{B}\right) $ maps
acyclic complexes to acyclic complexes, and hence it induces a canonical
triangulated functor $G:\mathrm{D}\left( \mathcal{A}\right) \rightarrow 
\mathrm{D}\left( \mathcal{B}\right) $ such that $GQ_{\mathcal{A}} $ is
isomorphic to $Q_{\mathcal{B}}F$.

In general, we say that $F$ has a \emph{total right derived functor }\cite[%
Section 2.2, Definition 1.2]{verdier_categories_1977} if there exists a
(essentially unique) triangulated functor $\mathrm{R}F:\mathrm{D}^{+}\left( 
\mathcal{A}\right) \rightarrow \mathrm{D}^{+}\left( \mathcal{B}\right) $
with a \emph{morphism} $\mu :Q_{\mathcal{B}}F\Rightarrow \left( \mathrm{R}%
F\right) Q_{\mathcal{A}}$ of triangulated functors (which is not necessarily
an isomorphism), such that for any other triangulated functor $G:\mathrm{D}%
^{+}\left( \mathcal{A}\right) \rightarrow \mathrm{D}^{+}\left( \mathcal{B}%
\right) $ and morphism $\nu :Q_{\mathcal{B}}F\Rightarrow GQ_{\mathcal{A}}$
of triangulated functors there exists a unique morphism $\sigma :\mathrm{R}%
F\Rightarrow G$ of triangulated functors such that $\sigma Q_{\mathcal{A}%
}\circ \mu =\nu $.

Let $\mathcal{E}$ be an abelian category and $\Phi :\mathrm{K}^{+}\left( 
\mathcal{A}\right) \rightarrow \mathcal{E}$ be a cohomological functor. Then
a \emph{cohomological right derived functor} of $F$ \cite[Section 2.2,
Definition 1.4]{verdier_categories_1977} is a cohomological functor $\mathrm{%
R}\Phi :\mathrm{D}^{+}\left( \mathcal{A}\right) \rightarrow \mathcal{E}$
together with a natural isomorphism $\mu :\Phi \Rightarrow \left( \mathrm{R}%
\Phi \right) Q_{\mathcal{A}}$ such that for any other cohomological functor $%
\Psi :\mathrm{D}^{+}\left( \mathcal{A}\right) \rightarrow \mathcal{E}$ and
natural transformation $\nu :\Phi \Rightarrow \Psi Q_{\mathcal{A}}$ there
exists a unique natural transformation $\sigma :\mathrm{R}\Phi \Rightarrow
\Psi $ such that $(\sigma Q_{\mathcal{A}})\circ \mu=\nu $. Such a
cohomological right derived functor, when it exists, is unique up to a
natural isomorphism.

\begin{definition}
Suppose that $\mathcal{A}$ is an exact category and $\mathcal{C}$ is a fully
exact subcategory of $\mathcal{A}$. We say that $\mathcal{C}$ is:

\begin{enumerate}
\item \emph{generating} \cite[Definition 8.3.21(vi)]%
{kashiwara_categories_2006} if for each object $A$ of $\mathcal{A}$ there
exists an admissible epimorphism $C\rightarrow A$ with $C$ in $\mathcal{C}$;

\item \emph{cogenerating} if $\mathcal{C}^{\mathrm{op}}$ is generating in $%
\mathcal{A}^{\mathrm{op}}$;

\item \emph{closed under quotients }if for every short exact sequence 
\begin{equation*}
0\rightarrow C^{\prime }\rightarrow C\rightarrow A^{\prime \prime
}\rightarrow 0
\end{equation*}%
in $\mathcal{A}$ with $C$ and $C^{\prime }$ in $\mathcal{C}$, $A^{\prime
\prime }$ is isomorphic to an object of $\mathcal{C}$;

\item \emph{closed under subobjects }if $\mathcal{C}^{\mathrm{op}}$ is
closed under quotients in $\mathcal{A}^{\mathrm{op}}$.
\end{enumerate}
\end{definition}

We have the following result \cite[Theorem 10.22 and Remark 10.23]%
{buhler_exact_2010}; see also \cite[Proposition 13.2.2]%
{kashiwara_categories_2006}, \cite[Lemma 1.3.3 and Proposition 1.3.4]%
{schneiders_quasi-abelian_1999}, and \cite[Corollary 3.10]%
{hoffmann_homological_2007}.

\begin{lemma}
\label{Lemma:same-derived}Suppose that $\mathcal{A}$ is an exact category
and $\mathcal{C}$ is a cogenerating fully exact subcategory of $\mathcal{A}$
closed under quotients. Then for any bounded complex $A$ in $\mathcal{A}$
there exists a bounded complex $C$ in $\mathcal{C}$ and a quasi-isomorphism $%
\eta :A\rightarrow C$ with $\eta ^{k}:A^{k}\rightarrow C^{k}$ an admissible
monic for every $k\in \mathbb{Z}$.

Furthermore, the inclusion $\mathcal{C}\rightarrow \mathcal{A}$ induces an
equivalence of categories%
\begin{equation*}
\mathrm{D}^{b}\left( \mathcal{C}\right) \rightarrow \mathrm{D}^{b}\left( 
\mathcal{A}\right) \text{.}
\end{equation*}
\end{lemma}

As a consequence of Lemma \ref{Lemma:same-derived} we have the following;
see \cite[Proposition 1.3.5]{schneiders_quasi-abelian_1999}.

\begin{proposition}
\label{Proposition:explicitly-right-derivable}Suppose that $F:\mathcal{A}%
\rightarrow \mathcal{R}$ is a functor between exact categories, and $%
\mathcal{D}$ is a cogenerating fully exact subcategory of $\mathcal{A}$
closed under quotients such that $F|_{\mathcal{D}}$ is exact. For a bounded
complex $A$ over $\mathcal{A}$, pick a bounded complex $D_{A}$ over $%
\mathcal{D}$ together with a quasi-isomorphism $\eta _{A}:A\rightarrow D_{A}$%
.\ Define $\left( \mathrm{R}F\right) \left( A\right) :=F\left( D_{A}\right) $
for a bounded complex $A$ over $\mathcal{A}$. This yields a triangulated
functor $\mathrm{R}F:\mathrm{D}^{b}\left( \mathcal{A}\right) \rightarrow 
\mathrm{D}^{b}\left( \mathcal{R}\right) $. Defining $\mu _{A}:=F\left( \eta
_{A}\right) :F\left( A\right) \rightarrow \left( \mathrm{R}F\right) \left(
A\right) $ for each bounded complex $A$ over $\mathcal{A}$ yields a morphism 
$\mu :Q_{\mathcal{R}}F\Rightarrow \left( \mathrm{R}F\right) Q_{\mathcal{A}}$
of triangulated functors. We have that $\left( \mathrm{R}F,\mu \right) $ is
the total right derived functor of $F$. Furthermore, if $\mathcal{R}$ is a
quasi-abelian category, then $\mathrm{H}^{0}\circ \mathrm{R}F$ is a
cohomological right derived functor of $\mathrm{H}^{0}\circ F$.
\end{proposition}

For future use, we isolate the following:

\begin{lemma}
\label{Lemma:pushouts-monic}Let $\mathcal{A}$ be an exact category.\ Then
the pushout of a monic along an admissible monic is a monic.
\end{lemma}

\begin{proof}
Consider a pushout%
\begin{equation*}
\begin{array}{ccc}
A & \overset{i}{\rightarrow } & A^{\prime } \\ 
m\downarrow &  & \downarrow m^{\prime } \\ 
B & \underset{j}{\rightarrow } & P%
\end{array}%
\end{equation*}%
where $i$ is an admissible monic, and $m$ is a monic. By \cite[Proposition
2.12]{buhler_exact_2010}, this gives an exact sequence%
\begin{equation*}
A\overset{(i,-m)}{\longrightarrow }A^{\prime }\oplus B\overset{(m^{\prime
},j)}{\longrightarrow }P.
\end{equation*}%
If $x:X\rightarrow A^{\prime }$ satisfies $m^{\prime }x=0$, then $(m^{\prime
},j)(x,0)=0$. Since $(i,-m)$ is the kernel of $(m^{\prime },j)$, there is $%
h:X\rightarrow A$ such that $ih=x$ and $mh=0$. Since $m$ is monic, $h=0$ and
hence $x=0$.
\end{proof}

\begin{lemma}
\label{Lemma:pushout-LH}Let $\mathcal{A}$ be a quasi-abelian category.
Consider a pushout diagram in $\mathcal{A}$%
\begin{equation*}
\begin{array}{ccc}
X & \overset{f}{\rightarrow } & Y \\ 
i\downarrow &  & \downarrow j \\ 
Z & \underset{f^{\prime }}{\rightarrow } & P%
\end{array}%
\end{equation*}%
where $i$ is a strict monic in $\mathcal{A}$. Then this is also a pushout
diagram in $\mathrm{LH}\left( \mathcal{A}\right) $.
\end{lemma}

\begin{proof}
Notice that $P$ is the cokernel of the arrow%
\begin{equation*}
\alpha =%
\begin{bmatrix}
f \\ 
-i%
\end{bmatrix}%
:X\rightarrow Y\oplus Z\text{;}
\end{equation*}%
see \cite[Proposition 2.12]{buhler_exact_2010}. Then $\alpha $ can be
written as the composition of the split (hence, strict) monic%
\begin{equation*}
\begin{bmatrix}
f \\ 
-\mathrm{id}_{X}%
\end{bmatrix}%
:X\rightarrow Y\oplus X
\end{equation*}%
and of the strict monic%
\begin{equation*}
\mathrm{id}_{Y}\oplus i:Y\oplus X\rightarrow Y\oplus Z\text{.}
\end{equation*}%
Since composition of strict monics is a strict monic \cite[Proposition 1.1.7]%
{schneiders_quasi-abelian_1999}, it follows that $\alpha $ is a strict
monic. Thus,%
\begin{equation*}
0\rightarrow X\overset{\alpha }{\rightarrow }Y\oplus Z\rightarrow
P\rightarrow 0
\end{equation*}%
is an exact sequence in $\mathcal{A}$. Since the inclusion $\mathcal{A}%
\rightarrow \mathrm{LH}\left( \mathcal{A}\right) $ is exact, it maps exact
sequences to exact sequences. Thus,%
\begin{equation*}
0\rightarrow X\overset{\alpha }{\rightarrow }Y\oplus Z\rightarrow
P\rightarrow 0
\end{equation*}%
is also an exact sequence in $\mathrm{LH}\left( \mathcal{A}\right) $,
yielding the desired conclusion.
\end{proof}

\begin{lemma}
\label{Lemma:presentation}Let $\mathcal{A}$ be a quasi-abelian category.
Suppose that $\mathcal{D}$ is a cogenerating subcategory of $\mathcal{A}$.
Let also $\mathcal{C}$ be a subcategory of $\mathcal{A}$ such that for every
object $A$ of $\mathcal{D}$ there exists a \emph{strict epimorphism} $%
q:A\rightarrow C$ such that for every monic arrow $m$ with source $A$, the
pushout of $m$ along $q$ in $\mathcal{A}$ is also the pushout of $m$ along $%
q $ in $\mathrm{LH}\left( \mathcal{A}\right) $. Then every object of $%
\mathrm{LH}\left( \mathcal{A}\right) $ is isomorphic to one of the form $G/N$
where $G$ is in $\mathcal{A}$ and $N$ is in $\mathcal{C}$.
\end{lemma}

\begin{proof}
Consider an object $G/N$ of $\mathrm{LH}\left( \mathcal{A}\right) $
represented by a monic arrow $N\rightarrow G$ in $\mathcal{A}$. Since $%
\mathcal{D}$ is cogenerating, there exists a strict monic%
\begin{equation*}
N\rightarrow D
\end{equation*}%
for some object $D$ of $\mathcal{D}$. Consider the pushout%
\begin{equation*}
\begin{array}{ccc}
N & \rightarrow & G \\ 
\downarrow &  & \downarrow \\ 
D & \rightarrow & G_{D}%
\end{array}%
\end{equation*}%
Then by Lemma \ref{Lemma:pushout-LH} this is also a pushout in $\mathrm{LH}%
\left( \mathcal{A}\right) $.\ Hence, $D\rightarrow G_{D}$ represents the
same object as $N\rightarrow G$ in $\mathrm{LH}\left( \mathcal{A}\right) $,
i.e.,%
\begin{equation*}
G_{D}/D\cong G/N\text{.}
\end{equation*}%
Thus, after replacing $N$ with $D$ and $G$ with $G_{D}$, we can assume that $%
N\in \mathcal{D}$. Then by hypothesis, there exists a strict epimorphism $%
q:N\rightarrow C$ for some object $C$ of $\mathcal{C}$, such that the pushout%
\begin{equation*}
\begin{array}{ccc}
N & \rightarrow & G \\ 
\downarrow &  & \downarrow \\ 
C & \rightarrow & G_{C}%
\end{array}%
\end{equation*}%
in $\mathcal{A}$ is also a pushout in $\mathrm{LH}\left( \mathcal{A}\right) $%
. Then again $N\rightarrow G$ and $C\rightarrow G_{C}$ represent the same
element in $\mathrm{LH}\left( \mathcal{A}\right) $, i.e., 
\begin{equation*}
G/N\cong G_{C}/C\text{.}
\end{equation*}%
This concludes the proof.
\end{proof}

\subsection{Bifunctors}

Suppose that $F:\mathcal{A}^{\mathrm{op}}\times \mathcal{B}\rightarrow 
\mathcal{R}$ is an additive functor, where $\mathcal{A}$, $\mathcal{B}$, and 
$\mathcal{R}$ are exact categories. Let $\mathrm{Ch}\left( \mathcal{A}%
\right) $ be the category of complexes over $\mathcal{A}$ and $\mathrm{Ch}%
^{2}\left( \mathcal{A}\right) $ be the category of \emph{double} complexes
over $\mathcal{A}$. Then we have that $F$ induces an additive functor $F:%
\mathrm{Ch}\left( \mathcal{A}\right) ^{\mathrm{op}}\times \mathrm{Ch}\left( 
\mathcal{B}\right) \rightarrow \mathrm{Ch}^{2}\left( \mathcal{R}\right) $
defined by%
\begin{equation*}
F\left( A,B\right) ^{i,j}:=F\left( A_{i},B^{j}\right)
\end{equation*}%
see \cite[Section 11.6]{kashiwara_categories_2006}. The vertical and
horizontal maps in $F\left( A,B\right) $ are given by%
\begin{equation*}
\delta _{\mathrm{v}}^{ij}:=F\left( \partial _{i}^{A},B^{j}\right) :F\left(
A_{i},B^{j}\right) \rightarrow F\left( A_{i+1},B^{j}\right)
\end{equation*}%
\begin{equation*}
\delta _{\mathrm{h}}^{ij}:=F(A_{i},\delta _{B}^{j}):F\left(
A_{i},B^{j}\right) \rightarrow F\left( A_{i},B^{j+1}\right) \text{.}
\end{equation*}%
The functor $F^{\bullet }:\mathrm{Ch}^{b}\left( \mathcal{A}\right) ^{\mathrm{%
op}}\times \mathrm{Ch}^{b}\left( \mathcal{B}\right) \rightarrow \mathrm{Ch}%
\left( \mathcal{R}\right) $ is defined by%
\begin{equation*}
F^{\bullet }\left( A,B\right) =\mathrm{Tot}\left( F\left( A,B\right) \right)
\end{equation*}%
where $\mathrm{Tot}\left( F\left( A,B\right) \right) $ denotes the \emph{%
total complex} of $F\left( A,B\right) $. This induces a \emph{triangulated}
bifunctor $F^{\bullet }:\mathrm{K}^{b}\left( \mathcal{A}\right) ^{\mathrm{op}%
}\times \mathrm{K}^{b}\left( \mathcal{B}\right) \rightarrow \mathrm{K}%
^{b}\left( \mathcal{R}\right) $; see \cite{neeman_triangulated_2001}.

If $\Phi :\mathrm{K}^{b}\left( \mathcal{A}\right) ^{\mathrm{op}}\times 
\mathrm{K}^{b}\left( \mathcal{B}\right) \rightarrow \mathrm{K}^{b}\left( 
\mathcal{R}\right) $ is a triangulated bifunctor, then a \emph{total right
derived functor }$\mathrm{R}\Phi $ of $\Phi $ is a triangulated bifunctor $%
\mathrm{D}^{b}\left( \mathcal{A}\right) ^{\mathrm{op}}\times \mathrm{D}%
^{b}\left( \mathcal{B}\right) \rightarrow \mathrm{D}^{b}\left( \mathcal{R}%
\right) $ with a morphism of triangulated functors $g:Q_{\mathcal{R}}\Phi
\Rightarrow \mathrm{R}\Phi \left( Q_{\mathcal{A}}\times Q_{\mathcal{B}%
}\right) $ such that for any other triangulated bifunctor $\Psi :\mathrm{D}%
^{b}\left( \mathcal{A}\right) ^{\mathrm{op}}\times \mathrm{D}^{b}\left( 
\mathcal{B}\right) \rightarrow \mathrm{D}^{b}\left( \mathcal{R}\right) $
with morphism of triangulated functors $h:Q_{\mathcal{R}}\Phi \Rightarrow
\Psi \left( Q_{\mathcal{A}}\times Q_{\mathcal{B}}\right) $ there exists a
unique morphism $\gamma :\mathrm{R}\Phi \Rightarrow \Psi $ of triangulated
functors such that $\gamma \left( Q_{\mathcal{A}}\times Q_{\mathcal{B}%
}\right) \circ g=h$. The same proof as \cite[Corollary 3.7 and Corollary 3.10%
]{hoffmann_homological_2007} gives the following.

\begin{proposition}
\label{Proposition:explicitly-right-derivable2}Suppose that $F:\mathcal{A}^{%
\mathrm{op}}\times \mathcal{B}\rightarrow \mathcal{R}$ is an additive
functor, where $\mathcal{A}$, $\mathcal{B}$, and $\mathcal{R}$ are exact
categories, $\mathcal{C}$ is a fully exact generating subcategory of $%
\mathcal{A}$ closed under subobjects, and $\mathcal{D}$ is a fully exact
cogenerating subcategory of $\mathcal{B}$ closed under quotients. Suppose
that $F\left( C,-\right) $ is exact on $\mathcal{D}$ for every object $C$ of 
$\mathcal{C}$, and $F\left( -,D\right) $ is exact on $\mathcal{C}$ for every
object $D$ of $\mathcal{D}$. For a bounded complex $A$ over $\mathcal{A}$
and a bounded complex $B$ over $\mathcal{B}$, pick a bounded complex $C_{A}$
over $\mathcal{C}$ and a bounded complex $D_{B}$ over $\mathcal{D}$ together
with quasi-isomorphisms $\eta _{A}:C_{A}\rightarrow A$ and $\eta
_{B}:B\rightarrow D_{B}$.\ Define $\left( \mathrm{R}F^{\bullet }\right)
\left( A,B\right) :=F^{\bullet }\left( C_{A},D_{B}\right) $. This yields a
triangulated functor $\mathrm{R}F^{\bullet }:\mathrm{D}^{b}\left( \mathcal{A}%
\right) ^{\mathrm{op}}\times \mathrm{D}^{b}\left( \mathcal{B}\right)
\rightarrow \mathrm{D}^{b}\left( \mathcal{R}\right) $. Defining $\mu
_{A}:=F^{\bullet }\left( \eta _{A},\eta _{B}\right) :F^{\bullet }\left(
A,B\right) \rightarrow \left( \mathrm{R}F^{\bullet }\right) \left(
A,B\right) $ for each bounded complex $A$ over $\mathcal{A}$ and bounded
complex $B$ over $\mathcal{B}$ yields a morphism $\mu :Q_{\mathcal{R}%
}F^{\bullet }\Rightarrow \left( \mathrm{R}F^{\bullet }\right) (Q_{\mathcal{A}%
}\times Q_{\mathcal{B}})$ of triangulated functors. We have that $\left( 
\mathrm{R}F^{\bullet },\mu \right) $ is the total right derived functor of $%
F^{\bullet }$. Furthermore, if $\mathcal{R}$ is a quasi-abelian category,
then $\mathrm{H}^{0}\circ \mathrm{R}F^{\bullet }$ is a cohomological right
derived functor of $\mathrm{H}^{0}\circ F^{\bullet }$.
\end{proposition}

\subsection{Injectives and projectives}

Let $\left( \mathcal{A},\mathcal{E}\right) $ be an exact category. An object 
$J$ of $\mathcal{A}$ is \emph{injective} in $\mathcal{A}$ if the functor $%
\mathrm{Hom}_{\mathcal{A}}\left( -,J\right) :\mathcal{A}^{\mathrm{op}%
}\rightarrow \mathbf{Ab}$ is exact. Equivalently, for every short exact
sequence $0\rightarrow A_{0}\rightarrow A_{1}\rightarrow A_{2}\rightarrow 0$
we have that the induced homomorphism%
\begin{equation*}
\mathrm{Hom}\left( A_{1},J\right) \rightarrow \mathrm{Hom}\left(
A_{0},J\right)
\end{equation*}%
is surjective. If $\mathcal{A}$ is quasi-abelian, then we have that $J$ is
injective in $\mathrm{LH}\left( \mathcal{A}\right) $ if and only if for
every monomorphism $u:A_{0}\rightarrow A_{1}$ in $\mathcal{A}$ (which we may
not in general assume to be a kernel) we have that the induced homomorphism%
\begin{equation*}
\mathrm{Hom}\left( A_{1},J\right) \rightarrow \mathrm{Hom}\left(
A_{0},J\right)
\end{equation*}%
is surjective; see \cite[Proposition 1.3.26(a)]%
{schneiders_quasi-abelian_1999}.

Dually, an object $P$ of $\left( \mathcal{A},\mathcal{E}\right) $ is
projective in $\mathcal{A}$ if the functor $\mathrm{Hom}_{\mathcal{A}}\left(
P,-\right) :\mathcal{A}\rightarrow \mathbf{Ab}$ is exact. Equivalently, for
every short exact sequence $0\rightarrow A_{0}\rightarrow A_{1}\rightarrow
A_{2}\rightarrow 0$ the induced homomorphism%
\begin{equation*}
\mathrm{Hom}\left( P,A_{1}\right) \rightarrow \mathrm{Hom}\left(
P,A_{2}\right)
\end{equation*}%
is surjective. By \cite[Proposition 1.3.24(a)]{schneiders_quasi-abelian_1999}%
, if $\mathcal{A}$ is quasi-abelian, then an object of $\mathcal{A}$ is
projective in $\mathcal{A}$ if and only if it is projective in $\mathrm{LH}%
\left( \mathcal{A}\right) $.

The category $\mathcal{A}$ \emph{has enough injectives }if for every object $%
X$ of $\mathcal{A}$ there exists an admissible monic $X\rightarrow I$ where $%
I$ is injective. Dually, $\mathcal{A}$ \emph{has enough projectives }if for
every object $X$ of $\mathcal{A}$ there exists an admissible epic $%
P\rightarrow X$ where $P$ is projective.

Suppose that $\mathrm{Hom}_{\mathcal{A}}^{\bullet }:\mathrm{K}^{b}\left( 
\mathcal{A}\right) ^{\mathrm{op}}\times \mathrm{K}^{b}\left( \mathcal{A}%
\right) \rightarrow \mathbf{Ab}$ has a total right derived functor $R\mathrm{%
Hom}_{\mathcal{A}}^{\bullet }$. Set $\mathrm{Ext}_{\mathcal{A}}^{n}:=\mathrm{%
H}^{n}\circ R\mathrm{Hom}_{\mathcal{A}}^{\bullet }$ for $n\in \mathbb{Z}$.
Then $J$ is injective in $\mathcal{A}$ if and only if $\mathrm{Ext}_{%
\mathcal{A}}^{1}\left( X,J\right) =0$ for every object $X$ of $\mathcal{A}$,
and injective in $\mathrm{LH}\left( \mathcal{A}\right) $ if and only if $%
\mathrm{Ext}_{\mathcal{A}}^{1}\left( X,J\right) =0$ for every object $X$ of $%
\mathrm{LH}\left( \mathcal{A}\right) $. Likewise, we have that $P$ is
projective in $\mathcal{A}$ (or, equivalently, in $\mathrm{LH}\left( 
\mathcal{A}\right) $) if and only if $\mathrm{Ext}_{\mathcal{A}}^{1}\left(
P,Y\right) =0$ for every object $Y$ of $\mathcal{A}$.

\subsection{Enriched (quasi-)abelian categories}

We would like to consider categories whose hom-sets are topological modules.
To formalize this notion within the context of enriched categories, the
abstract notion of \emph{category of modules} $\mathcal{M}$ has been
introduced in \cite[Section 2.15]{casarosa_phantom_2025}. If $\mathcal{M}$
is a category of modules, then its hom-sets have a canonical abelian group
structure. It is thus meaningful to define a (quasi-)abelian category of
modules, in reference to this canonical $\mathbf{Ab}$-enrichment.

Since $\mathcal{M}$ is in particular required to be a monoidal tensor
category, one can consider the notion of $\mathcal{M}$-\emph{category}. When 
$\mathcal{M}$ is a quasi-abelian category of modules, the notion of
(quasi-)abelian $\mathcal{M}$-category $\mathcal{A}$ admits a natural
definition. This is obtained by requiring that the (quasi-)abelian structure
on $\mathcal{A}$ be compatible with the $\mathcal{M}$-enrichment. Similar
considerations apply to triangulated categories.

It is observed in \cite[Section 2.15]{casarosa_phantom_2025} that starting
from a category of modules $\mathcal{M}$, one can obtain new categories of
modules through some canonical constructions. Indeed, the category of
epimorphic towers $\boldsymbol{\Pi }\left( \mathcal{M}\right) $ is also a
category of modules \cite[Section 2.15]{casarosa_phantom_2025}, as is its
left heart $\mathrm{LH}\left( \boldsymbol{\Pi }\left( \mathcal{M}\right)
\right) $.

Suppose that $\mathcal{M}$ is a quasi-abelian category of modules, and $%
\mathcal{A}$ is a quasi-abelian $\mathcal{M}$-category. Then $\mathrm{Hom}_{%
\mathcal{A}}$ is a functor $\mathcal{A}^{\mathrm{op}}\times \mathcal{A}%
\rightarrow \mathcal{M}$. The corresponding functor $\mathrm{Hom}^{\bullet }:%
\mathrm{K}\left( \mathcal{A}\right) ^{\mathrm{op}}\times \mathrm{K}\left( 
\mathcal{A}\right) \rightarrow \mathrm{K}\left( \mathcal{M}\right) $ defined
as in Proposition \ref{Proposition:explicitly-right-derivable2} is such that 
$\left( \mathrm{H}^{0}\circ \mathrm{Hom}^{\bullet }\right) \left( A,B\right) 
$ is naturally isomorphic to $\mathrm{Hom}_{\mathrm{K}^{b}\left( \mathcal{A}%
\right) }\left( A,B\right) $; see \cite[Proposition 11.7.3]%
{kashiwara_categories_2006}.

\section{Groups with a Polish cover and definable functors}

\label{Section:groups-with-a-Polish-cover}In this section we recall the
definition of several categories of locally compact Polish abelian groups,
and discuss their parameterizations by points of a standard Borel space.
Standard references for descriptive set theory and standard Borel spaces
include \cite{kechris_classical_1995,gao_invariant_2009}.

\subsection{The left heart of abelian Polish groups}

Recall that a \emph{Polish space }is a topological space whose topology is
second-countable and induced by a complete metric. A Polish group is a group
object in the category of Polish spaces, namely a Polish space endowed with
a continuous group operation (whence the map assigning to each object its
inverse is also continuous). It follows from the Open\ Mapping Theorem for
Polish groups \cite[Theorem 2.3.3]{gao_invariant_2009} that abelian Polish
groups form a quasi-abelian category $\mathbf{PAb}$ with continuous group
homomorphisms as morphisms. The corresponding left heart $\mathrm{LH}\left( 
\mathbf{PAb}\right) $ is shown in \cite{lupini_looking_2024} to be
equivalent to the category of abelian \emph{groups with a Polish cover}
(called \emph{pseudo-Polish groups} by Moore in \cite{moore_group_1976})
with \emph{Borel-definable group homomorphisms} as morphisms (see \cite%
{bergfalk_definable_2024-1,bergfalk_definable_2024} for further discussion
of this category). An abelian group with a Polish cover is a group $G$
explicitly presented as a quotient $\hat{G}/N$ where $\hat{G}$ is an abelian
Polish group and $N$ is a \emph{Polish subgroup} of $\hat{G}$. This means
that $N$ is endowed with a (necessarily unique) Polish group topology that
makes the inclusion $N\rightarrow \hat{G}$ continuous or, equivalently, a
Borel isomorphism onto the image (endowed with the Borel structure induced
from $\hat{G}$). A group homomorphism $\varphi :G\rightarrow H$ between
groups with a Polish cover $G=\hat{G}/N$ and $H=\hat{H}/M$ is
Borel-definable if there exists a Borel function $\hat{\varphi}:\hat{G}%
\rightarrow \hat{H}$ that is a \emph{lift} of $\varphi $, in the sense that $%
\hat{\varphi}(x)+M=\varphi \left( x+N\right) $ for every $x\in \hat{G}$.
Several equivalent characterizations of Borel-definable group homomorphisms
are provided in \cite[Theorem 4.6 and Proposition 4.7]{lupini_looking_2024}.

More generally, we have from \cite[Theorem 6.14]{lupini_looking_2024} the
following:

\begin{theorem}
If $\mathcal{B}$ is a thick subcategory of $\mathbf{PAb}$, then $\mathrm{LH}%
\left( \mathcal{B}\right) $ is equivalent to the full subcategory of the
category of abelian groups with a Polish cover spanned by the groups with a $%
\mathcal{B}$-cover, which are the groups with a Polish cover of the form $G=%
\hat{G}/N$ where both $\hat{G}$ and $N$ are objects of $\mathcal{B}$.
\end{theorem}

For future reference, we notice that if $\psi :H\rightarrow G$ is a
continuous homomorphism between Polish abelian groups, then its image $N$ is
a Polish subgroup of $G$. Indeed, $M:=\mathrm{Ker}\left( \psi \right) $ is a
closed subgroup of $H$, and the quotient topology on $N$ renders it a Polish
group isomorphic to $H/M$. Furthermore, the inclusion $N\rightarrow G$ is
continuous by continuity of $\psi :H\rightarrow G$.

\subsection{Borel-definable sets and groups}

A generalization of the notion of group with a Polish cover was introduced
in \cite{bergfalk_definable_2024} in the context of Borel-definable sets. A
slightly more restrictive notion has been considered in \cite[Section 13.2]%
{casarosa_phantom_2025}, which includes all the known examples.

Recall that an equivalence relation $E$ on a standard Borel space $\hat{X}$
is \emph{Borel} if it is a Borel subset of $\hat{X}\times \hat{X}$, and 
\emph{idealistic }if there exists a Borel function $s:\hat{X}\rightarrow 
\hat{X}$ and a nontrivial Borel-on-Borel assignment $C\mapsto \mathcal{F}%
_{C} $ from $E$-classes to $\sigma $-filters. By definition, this means
that, for every $C$, $\mathcal{F}_{C}$ is a $\sigma $-filter of subsets of $%
C $ that does not contain the empty set, and for every Borel subset $A$ of $%
\hat{X}\times \hat{X}$ the set%
\begin{equation*}
\left\{x\in\hat{X}:\left\{y\in[x]_{E}:\left(s(x),y\right)\in A\right\} \in%
\mathcal{F}_{[x]_{E}}\right\}
\end{equation*}%
is a Borel subset of $\hat{X}$. From the notion of \emph{idealistic}
equivalence relation one obtains the notion of \emph{stably idealistic }%
equivalence relation, where one requires not only $E$ to be idealistic, but
also its infinite product $E^{\mathbb{N}}$, as an equivalence relation on $%
\hat{X}^{\mathbb{N}}$.

A Borel-definable set as defined in \cite[Section 13.2]%
{casarosa_phantom_2025} is a set $X$ explicitly presented as a quotient $%
\hat{X}/E$ where $\hat{X}$ is a standard Borel space and $E$ is an
equivalence relation on $\hat{X}$ that is Borel and idealistic. A subset $A$
of a Borel-definable set $\hat{X}/E$ is \emph{Borel} if $\hat{A}:=\{x\in 
\hat{X}:[x]_{E}\in A\}$ is a Borel subset of $\hat{X}$, and $\boldsymbol{%
\Sigma }_{1}^{1}$ or analytic if $\hat{A}$ is an analytic subset of $\hat{X}$%
. The product $X\times Y$ of two Borel-definable sets $X=\hat{X}/E$ and $Y=%
\hat{Y}/F$ is the Borel-definable set $(\hat{X}\times \hat{Y})/(E\times F)$.

A function $f:X\rightarrow Y$ between Borel-definable sets $X=\hat{X}/E$ and 
$Y=\hat{Y}/F$ is \emph{Borel-definable }if there exists a Borel function $%
\hat{f}:\hat{X}\rightarrow \hat{Y}$ that is a \emph{lift }of $f$, in the
sense that $f\left( [x]_{E}\right) =[\hat{f}(x)]_{F}$ for every $x\in \hat{X}
$. By the version of the \textquotedblleft large section\textquotedblright\
uniformization theorem \cite[Theorem 18.6]{kechris_classical_1995} presented
in the proof of \cite[Lemma 3.7]{kechris_borel_2016}, this is equivalent to
the assertion that the graph $\Gamma \left( f\right) $ of $f$ is a Borel
subset of $X\times Y$. This gives the notion of morphism in the category $%
\mathbf{DSet}$ of Borel-definable sets. A Borel-definable group is simply a
group object in the category of Borel-definable sets.

More generally, one can consider the notion of $\boldsymbol{\Sigma }_{1}^{1}$%
-\emph{definable set}, which is a set presented as a quotient $X=\hat{X}/E$
where $\hat{X}$ is a Polish space and $E$ is an equivalence relation on $%
\hat{X}$ that is merely required to be analytic. The notions of Borel and $%
\boldsymbol{\Sigma }_{1}^{1}$ subset of a $\boldsymbol{\Sigma }_{1}^{1}$%
-definable set, and of Borel-definable function between $\boldsymbol{\Sigma }%
_{1}^{1}$-definable sets, are defined as in the case of Borel-definable
sets. We say that a function between $\boldsymbol{\Sigma }_{1}^{1}$%
-definable sets is $\boldsymbol{\Sigma }_{1}^{1}$-definable if its graph is $%
\boldsymbol{\Sigma }_{1}^{1}$. The category of $\boldsymbol{\Sigma }_{1}^{1}$%
-definable sets and $\boldsymbol{\Sigma }_{1}^{1}$-definable functions will
be denoted by $\mathbf{\Sigma }_{1}^{1}$-$\mathbf{DSet}$.

We let $\mathbf{DAb}$ be the full subcategory of $\mathbf{DSet}$ spanned by
Borel-definable abelian groups, and let $\boldsymbol{\Sigma }_{1}^{1}$-$%
\mathbf{DAb}$ be the full subcategory of $\boldsymbol{\Sigma }_{1}^{1}$-$%
\mathbf{DSet}$ spanned by $\boldsymbol{\Sigma }_{1}^{1}$-definable abelian
groups. An abelian group with a Polish cover is, in particular, a
Borel-definable group, whence $\mathrm{LH}\left( \mathbf{PAb}\right) $ is
(equivalent to) a full subcategory of $\mathbf{DAb}$.

\begin{lemma}
\label{Lemma:Borel-definable}Suppose that $X=\hat{X}/E$ and $Y=\hat{Y}/F$
are $\boldsymbol{\Sigma }_{1}^{1}$-definable sets, and $f:X\rightarrow Y$ is
a $\boldsymbol{\Sigma }_{1}^{1}$-definable function.

\begin{enumerate}
\item If $F$ is Borel, then $\Gamma \left( f\right) $ is Borel.

\item If $E$ is Borel and $f$ is injective, then $\Gamma \left( f\right) $
is Borel.

\item If $E$ is Borel and $f$ is bijective, then $F$ is Borel;

\item If $\Gamma \left( f\right) $ is Borel and $F$ is idealistic, then $f$
is Borel-definable.
\end{enumerate}
\end{lemma}

\begin{proof}
(1) In order to prove that $\Gamma \left( f\right) $ is Borel, we need to
show that 
\begin{equation*}
\hat{\Gamma}\left( f\right) :=\{(x,y)\in \hat{X}\times \hat{Y}:f\left(
[x]_{E}\right) =[y]_{F}\}
\end{equation*}%
is co-analytic. We have that, for $(x,y)\in \hat{X}\times \hat{Y}$, $%
(x,y)\in \hat{\Gamma}\left( f\right) $ if and only if 
\begin{equation*}
\forall y_{1}\in \hat{Y}\text{, }\left( x,y_{1}\right) \in \hat{\Gamma}%
\left( f\right) \Rightarrow y_{1}Fy\text{.}
\end{equation*}

\noindent (2) As in (1), we need to prove that $\hat{\Gamma}\left( f\right) $
is co-analytic. Since $f$ is injective, we have that $(x,y)\in \hat{\Gamma}%
\left( f\right) $ if and only if%
\begin{equation*}
\forall x_{1}\in \hat{X}\text{, }\left( x_{1},y\right) \in \hat{\Gamma}%
\left( f\right) \Rightarrow x_{1}Ex\text{.}
\end{equation*}%
(3) If $f$ is bijective, then for $y,y^{\prime }\in \hat{Y}$ we have that $%
yFy^{\prime }$ if and only if 
\begin{equation*}
\forall x,x^{\prime }\in \hat{X}\text{, }(x,y)\in \hat{\Gamma}\left(
f\right) \wedge \left( x^{\prime },y^{\prime }\right) \in \hat{\Gamma}\left(
f\right) \Rightarrow xEx^{\prime }\text{.}
\end{equation*}%
This shows that $F$ is co-analytic, and hence Borel.

(4) This follows from the \textquotedblleft large section\textquotedblright\
uniformization theorem \cite[Theorem 18.6]{kechris_classical_1995} presented
in the proof of \cite[Lemma 3.7]{kechris_borel_2016}.
\end{proof}

It follows in particular from Lemma \ref{Lemma:Borel-definable} that $%
\mathbf{DSet}$ is a full subcategory of $\boldsymbol{\Sigma }_{1}^{1}$-$%
\mathbf{DSet}$, and $\mathbf{DAb}$ is a full subcategory of $\boldsymbol{%
\Sigma }_{1}^{1}$-$\mathbf{DAb}$.

\subsection{Quasi-abelian standard Borel categories\label%
{Subsection:Borel-quasi-abelian}}

Following \cite{chen_representing_2019}---see also \cite{lupini_polish_2017}%
---we say that a category $\mathcal{C}$ is a \emph{(standard) Borel category}
if its morphisms form a standard Borel set such that the set $\mathrm{Ob}%
\left( \mathcal{C}\right) $ of objects (identified with the corresponding
identity arrows) is a Borel subset of $\mathrm{Mor}\left( \mathcal{C}\right) 
$, the set $\mathrm{Inv}(\mathcal{C})$ of invertible arrows is a Borel
subset of $\mathrm{Mor}\left( \mathcal{C}\right) $, inversion of invertible
arrows is a Borel function $\mathrm{Inv}\left( \mathcal{C}\right)
\rightarrow \mathrm{Inv}\left( \mathcal{C}\right) $, source and target maps
are Borel functions $\mathrm{Mor}\left( \mathcal{C}\right) \rightarrow 
\mathrm{Ob}(\mathcal{C})$, and composition of arrows is a Borel function 
\begin{equation*}
\mathrm{Mor}\left( \mathcal{C}\right) ^{\left( 2\right) }:=\left\{ \left(
g,f\right) \in \mathrm{Mor}\left( \mathcal{C}\right) \times \mathrm{\mathrm{%
Mor}}\left( \mathcal{C}\right) :t\left( f\right) =s\left( g\right) \right\}
\rightarrow \mathrm{Mor}\left( \mathcal{C}\right)
\end{equation*}%
A \emph{Borel additive category} $\mathcal{C}$ is a Borel category which is
also additive, and such that:

\begin{itemize}
\item addition of arrows is a Borel function 
\begin{equation*}
\left\{ \left( f_{0},f_{1}\right) \in \mathrm{Mor}\left( \mathcal{C}\right)
\times \mathrm{Mor}\left( \mathcal{C}\right) :s\left( f_{0}\right) =s\left(
f_{1}\right) ,t\left( f_{0}\right) =t\left( f_{1}\right) \right\}
\rightarrow \mathrm{Mor}\left( \mathcal{C}\right) \text{,}
\end{equation*}

\item there exist Borel functions $\mathrm{Mor}\left( \mathcal{C}\right)
\rightarrow \mathrm{Mor}\left( \mathcal{C}\right) $, $f\mapsto \mathrm{ker}%
\left( f\right) $ and $f\mapsto \mathrm{coker}\left( f\right) $ such that $%
\mathrm{ker}\left( f\right) $ and $\mathrm{coker}\left( f\right) $ are a
kernel and a cokernel of $f$, respectively;

\item there exists a Borel function $\mathrm{Ob}\left( \mathcal{C}\right)
\times \mathrm{Ob}\left( \mathcal{C}\right) \rightarrow \mathrm{Ob}\left( 
\mathcal{C}\right) \times \mathrm{Mor}\left( \mathcal{C}\right) \times 
\mathrm{Mor}\left( \mathcal{C}\right) $, 
\begin{equation*}
\left( A,B\right) \mapsto \left( A\oplus B,\pi _{A},\pi _{B},\iota
_{A},\iota _{B}\right)
\end{equation*}
such that $A\oplus B$ is a biproduct of $A$ and $B$ with canonical
projections $\pi _{A}$ and $\pi _{B}$ and canonical coprojections $\iota
_{A} $ and $\iota _{B}$.
\end{itemize}

Informally, the last two hypotheses assert that kernels, cokernels, and
biproducts in a Borel additive category are given by Borel functions. The
expression of limits and colimits in terms of kernels, cokernels, products,
and coproducts shows that in fact all finite limits and colimits are given
by Borel functions.

Suppose now that $\mathcal{A}$ is a Borel quasi-abelian category. Then the
category of complexes $\mathrm{K}^{b}\left( \mathcal{A}\right) $ over $%
\mathcal{A}$ is a $\mathbf{DAb}$-enriched triangulated category. Indeed, for
bounded complexes $A^{\bullet },B^{\bullet }$ over $\mathcal{A}$, we have
that $\mathrm{Hom}_{\mathrm{K}^{b}\left( \mathcal{A}\right) }\left(
A^{\bullet },B^{\bullet }\right) $ is the quotient of the abelian Polish
group of chain maps $A^{\bullet }\rightarrow B^{\bullet }$ modulo the Polish
subgroup of nullhomotopic chain maps. Likewise, the derived category $%
\mathrm{D}^{b}\left( \mathcal{A}\right) $ is a $\boldsymbol{\Sigma }_{1}^{1}$%
-$\mathbf{DAb}$-enriched triangulated category, and the quotient functor $Q_{%
\mathcal{A}}:\mathrm{K}^{b}\left( \mathcal{A}\right) \rightarrow \mathrm{D}%
^{b}\left( \mathcal{A}\right) $ is $\boldsymbol{\Sigma }_{1}^{1}$-$\mathbf{%
DAb}$-enriched. Indeed, for bounded complexes $A^{\bullet },B^{\bullet }$
over $\mathcal{A}$, we have that $\mathrm{Hom}_{\mathrm{D}^{b}\left( 
\mathcal{A}\right) }\left( A^{\bullet },B^{\bullet }\right) $ is the
quotient of the standard Borel space of pairs $\left( f,\sigma \right) $
where $\sigma :A^{\prime \bullet }\rightarrow A^{\bullet }$ is a
quasi-isomorphism and $f:A^{\prime \bullet }\rightarrow B^{\bullet }$ is a
chain map, modulo the analytic equivalence relation $E$ defined by setting $%
\left( f_{0},\sigma _{0}\right) E_{\mathrm{D}^{b}\left( \mathcal{A}\right)
}\left( f_{1},\sigma _{1}\right) $ if and only if there exist
quasi-isomorphisms $\tau _{0},\tau _{1}$ such that $f_{0}\tau _{0}\sim
f_{1}\tau _{1}$ and $\sigma _{0}\tau _{0}\sim \sigma _{1}\tau _{1}$, where
we let $\sim $ denote chain homotopy. The composition of arrows $A^{\bullet
}\rightarrow B^{\bullet }$ and $B^{\bullet }\rightarrow C^{\bullet }$ is
given by the Borel-definable function $\mathrm{Hom}_{\mathrm{D}^{b}\left( 
\mathcal{A}\right) }\left( A^{\bullet },B^{\bullet }\right) \times \mathrm{%
Hom}_{\mathrm{D}^{b}\left( \mathcal{A}\right) }\left( B^{\bullet
},C^{\bullet }\right) \rightarrow \mathrm{Hom}_{\mathrm{D}^{b}\left( 
\mathcal{A}\right) }\left( A^{\bullet },C^{\bullet }\right) $ defined by
setting $[\left( f,\sigma \right) ]\circ \lbrack \left( g,\tau \right)
]=[\left( t,\lambda \right) ]$ if and only if for some or, equivalently,
every quasi-isomorphism $\rho $ and chain map $h$ such that $f\rho \sim \tau
h$ one has that $\left( gh,\sigma \rho \right) E_{\mathrm{D}^{b}\left( 
\mathcal{A}\right) }\left( t,\lambda \right) $. The translation functor on $%
\mathrm{D}^{b}\left( \mathcal{A}\right) $ is clearly $\boldsymbol{\Sigma }%
_{1}^{1}$-$\mathbf{DAb}$-enriched as well. By the foregoing discussion, the
left heart $\mathrm{LH}\left( \mathcal{A}\right) $ of $\mathcal{A}$, as a
full subcategory of $\mathrm{D}^{b}\left( \mathcal{A}\right) $, becomes a $%
\boldsymbol{\Sigma }_{1}^{1}$-$\mathbf{DAb}$-enriched abelian category.

More precisely, $\mathrm{LH}\left( \mathcal{A}\right) $ is the full
subcategory of $\mathrm{D}^{b}\left( \mathcal{A}\right) $ spanned by $H^{0}$%
-complexes, which are by definition the complexes $A$ such that $H^{n}\left(
A\right) =0$ for $n\in \mathbb{Z}\setminus \left\{ 0\right\} $ or,
equivalently, are isomorphic in $\mathrm{D}^{b}\left( \mathcal{A}\right) $
to a monic $\delta ^{-1}:B^{-1}\rightarrow B^{0}$, regarded as a complex $B$
over $\mathcal{A}$ such that $B^{n}=0$ for $n\notin \left\{ -1,0\right\} $;
see \cite[Section 1.2]{schneiders_quasi-abelian_1999} and also \cite[III.5.1]%
{gelfand_methods_2003}.

\subsection{Classes of Polish abelian groups\label{Subsection:classes}}

In this subsection we introduce a number of thick subcategories of the
quasi-abelian category $\mathbf{PAb}$ of Polish abelian groups. For more
information on these classes, we refer the reader to \cite%
{armacost_structure_1981,moskowitz_homological_1967,hoffmann_homological_2007,fulp_extensions_1971}
as well as \cite[Section 6]{lupini_looking_2024}. We let $\mathbf{LCPAb}$ be
the thick subcategory of $\mathbf{PAb}$ consisting of \emph{locally compact }%
Polish abelian groups.

A Polish abelian group is \emph{non-Archimedean }if it has a basis of zero
neighborhoods consisting of subgroups.

If $G$ is a locally compact Polish abelian group, we let $c\left( G\right) $
be the connected component of zero. We have that $c\left( G\right) $ is the
intersection of all the open subgroups of $G$ \cite[Theorem 7.8]%
{hewitt_abstract_1979}. Furthermore, $c\left( G\right) $ is a closed
connected subgroup such that $G/c\left( G\right) $ is totally disconnected 
\cite[Theorem 7.3]{hewitt_abstract_1979}. Conversely, if $H$ is a closed
connected subgroup of $G$ such that $G/H$ is totally disconnected, then $%
H=c\left( G\right) $.

Non-Archimedean Polish abelian groups form a thick subcategory of $\mathbf{%
PAb}$ {\cite[Theorem 6.18]{lupini_looking_2024}}. A locally compact Polish
abelian group is non-Archimedean if and only if it is \emph{totally
disconnected}. We let $\mathbf{TDLCPAb}$ be the thick subcategory of $%
\mathbf{LCPAb}$ consisting of totally disconnected locally compact Polish
abelian groups.

A \emph{vector group }is a locally compact Polish abelian group isomorphic
to $\mathbb{R}^{n}$ for some $n\in \omega $. A \emph{torus }is a locally
compact Polish abelian group isomorphic to $\mathbb{T}^{\sigma }$ for some $%
\sigma \in \omega +1$.

If $\mathcal{A}$ is a thick subcategory of $\mathbf{LCPAb}$, we let $%
\mathcal{A}_{\mathrm{c}}$ be the thick subcategory consisting of compact
groups in $\mathcal{A}$.

A locally compact Polish abelian group $G$ is \emph{compactly generated}%
\footnote{%
(The terminological overlap with the notion \emph{compactly generated} for
more general topological spaces is unfortunate, but entrenched.)} if there
exists a compact (symmetric) zero neighborhood $U$ such that $U$ is an
algebraic set of generators for $G$. Let $\mathbf{LCPAb}_{\mathrm{cg}}$ be
the thick subcategory of $\mathbf{LCPAb}$ consisting of compactly generated
locally compact Polish abelian groups; see \cite[Theorem 2.6]%
{moskowitz_homological_1967}. A Polish abelian group is a compactly
generated locally compact Polish abelian group if and only if it is
isomorphic to $V\oplus C\oplus F$ where $V$ is a finite-dimensional vector
group, $C$ is a compact group, and $F$ is a finite-rank free abelian group 
\cite[Theorem 2.5]{moskowitz_homological_1967}.

\enlargethispage{2pt} A locally compact Polish abelian group $G$ is a \emph{%
Lie group} if and only if its connected component of zero $c\left( G\right) $
is a Lie group, if and only if $G$ is isomorphic to $V\oplus C\oplus D$
where $V$ is a vector group, $C$ is a finite-dimensional torus, and $D$ is a
countable discrete abelian group. We let $\mathbf{LieAb}$ be the thick
subcategory of $\mathbf{LCPAb}$ consisting of abelian Lie groups; see \cite[%
Theorem 6.18]{lupini_looking_2024}. As observed in \cite[Section 2.9]%
{casarosa_homological_2026}, the category $\mathbf{LieAb}$ is a category of
modules in the sense of \cite[Section 2.15]{casarosa_phantom_2025}. The
point is that $\mathbf{LieAb}$ is endowed with a canonical monoidal
structure obtained by defining $G\otimes H$ to be the Lie abelian Polish
group whose Pontryagin dual is $\mathrm{Hom}(G,H^{\vee })$, where $H^{\vee }$
is the Pontryagin dual of $H$; see \cite[Section\ IV]%
{moskowitz_homological_1967} and \cite{garling_tensor_1966}. By duality, we
have 
\begin{equation*}
\mathrm{Hom}(G\otimes H,\mathbb{T})\cong \mathrm{Hom}(G,\mathrm{Hom}\left( H,%
\mathbb{T}\right) )
\end{equation*}%
and more generally for any other Lie abelian Polish group $L$,%
\begin{equation*}
\mathrm{Hom}(G\otimes H,L)\cong \mathrm{Hom}(G,\mathrm{Hom}\left( H,L\right)
)
\end{equation*}%
The \emph{category of modules }structure on $\mathbf{LieAb}$ induces a
category of modules structure on the category $\mathbf{\Pi }\left( \mathbf{%
LieAb}\right) $ of \emph{epimorphic towers }of $\mathbf{LieAb}$; see \cite[%
Section 2.15]{casarosa_phantom_2025}. This is equivalent to the category $%
\mathbf{proLiePAb}$ of \emph{pro-Lie }Polish abelian groups from \cite%
{casarosa_phantom_2025}, which is the category of Polish abelian group $G$
with the property that every zero neighborhood $U$ contains a closed
subgroup $H$ such that $G/H$ is a Lie group. In turn, this renders its left
heart $\mathrm{LH}\left( \mathbf{proLiePAb}\right) $ a category of modules
as well.

The category $\mathbf{proLiePAb}$ contains the category $\mathbf{LCPAb}$ of
locally compact Polish abelian groups as a full subcategory \cite[Examples
5.1(i)]{hofmann_lie_2007}.\ Furthermore, if $G$ and $H$ are locally compact
Polish abelian groups, then their hom-set $\mathrm{Hom}\left( G,H\right) $
endowed with the compact-open topology is pro-Lie \cite[Proposition 4.22]%
{casarosa_homological_2026}. This renders the category $\mathbf{LCPAb}$ a
quasi-abelian $\mathbf{proLiePAb}$-category in the sense of \cite[Section
2.15]{casarosa_phantom_2025}.

A locally compact Polish abelian group $G$ is a \emph{topological $p$-group}
if $\lim_{n\rightarrow \infty }{}p^{n}x=0$ for every $x\in G $. This is
equivalent to the assertion that $G$ has a basis of zero neighborhoods
consisting of open subgroups $U$ such that $G/U$ is a countable $p$-group.
We let $\mathbf{LCPAb}\left( p\right) $ be a thick subcategory of $\mathbf{%
TDLCPAb}$ consisting of topological $p$-groups.

A locally compact Polish abelian group $G$ is a \emph{topological torsion
group} if $\lim_{n\rightarrow \infty }{}n!x=0$ for every $x\in G$, or,
equivalently, if $G$ has a basis of zero neighborhoods consisting of open
subgroups $U$ such that $G/U$ is a countable torsion group. We let $\mathbf{%
TorLCPAb}$ be the thick subcategory of $\mathbf{TDLCPAb}$ consisting of
topological torsion groups.

For every topological torsion group $G$, the union $G_{p}$ of the images of
continuous homomorphisms $\mathbb{Z}_{p}\rightarrow G$ is a subgroup, called
the \emph{$p$-component} of $G$. This group $G_{p}$ is a topological $p$%
-group and $G\cong G_{p}^{\prime }\oplus G_{p}$ where $G_{p}^{\prime }$ has
no nonzero closed subgroups that are $p$-groups. If $U$ is a compact open
subgroup of $G$, then $G$ is isomorphic to the restricted product $%
\prod_{p}\left( G_{p}:G_{p}\cap U\right) $ ranging over all primes. By
definition, $\prod_{p}\left( G_{p}:G_{p}\cap U\right) $ is the set of
sequences $\left( x_{p}\right) $ indexed by primes, such that the set $%
\left\{ p:x_{p}\notin U\right\} $ is finite. The topology on $%
\prod_{p}\left( G_{p}:G_{p}\cap U\right) $ is the coarsest topology that is
finer than the product topology and that has $\prod_{p}(G_{p}\cap U)$ as an
open set.

For a prime $p$ and $k\in \omega $ we let $\mathbb{Z}(p^{k})$ be the group $%
\mathbb{Z}/p^{k}\mathbb{Z}$. We also denote by $\mathbb{Z}\left( p^{\infty
}\right) $ the subgroup of $\mathbb{T}=\mathbb{R}/\mathbb{Z}$ consisting of
roots of unity of $p$-power order; hence $\mathbb{Z}\left( p^{\infty
}\right) =\mathrm{co\mathrm{lim}}_{k}{}\mathbb{Z}(p^{k})$.

A locally compact Polish abelian group $G$ has:

\begin{itemize}
\item type $\mathbb{Z}$ if it is countable and torsion-free;

\item type $\mathbb{S}^{1}$ if it is compact and connected;

\item type $\mathbb{A}$ if its connected component $c\left( G\right) $ is a
vector group (in which case $c\left( G\right) $ is a direct summand of $G$),
and $G/c\left( G\right) $ is a topological torsion group.
\end{itemize}

We let $\mathbf{LCPAb}_{\mathbb{A}}$ be the full subcategory of $\mathbf{%
LCPAb}$ consisting of type $\mathbb{A}$ groups. Notice that $\mathbf{LCPAb}_{%
\mathbb{A}}$ is not thick, since $\mathbb{T}$ is a quotient by a closed
subgroup of the object $\mathbb{R}$ of $\mathbf{LCPAb}_{\mathbb{A}}$.

Every locally compact abelian Polish group $G$ admits canonical subgroups $%
G_{\mathbb{S}^{1}}\subseteq F_{\mathbb{Z}}G\subseteq G$ such that $G_{%
\mathbb{S}^{1}}$ has type $\mathbb{S}^{1}$, $G_{\mathbb{A}}:=F_{\mathbb{Z}%
}G/G_{\mathbb{S}^{1}}$ has type $\mathbb{A}$, and $G/F_{\mathbb{Z}}G:=G_{%
\mathbb{Z}}$ has type $\mathbb{Z}$ {\cite[Proposition 2.2]%
{hoffmann_homological_2007}}. One lets $G_{\mathbb{R}}$ be the connected
component of $G_{\mathbb{A}}$ and $G_{\mathrm{t}}$ be $G/c\left( G\right) $.
For a prime number $p$, we let the $p$-component $G_{p}$ of $G$ be the $p$%
-component of $G_{\mathrm{t}}$. One says that $G$:

\begin{itemize}
\item \emph{has finite $\mathbb{Z}$-rank} if $G_{\mathbb{Z}}$ has finite
torsion-free rank;

\item \emph{has finite $\mathbb{S}^{1}$-rank} if $G_{\mathbb{S}^{1}}$ has
finite dimension;

\item \emph{has finite $p$-rank} for some prime $p$ if the map $x\mapsto px$
on $G_{p}$ has finite kernel and cokernel or, equivalently, $G_{p}$ is a
finite direct sum of copies of $\mathbb{Z}_{p}$, $\mathbb{Q}_{p}$, $\mathbb{Z%
}\left( p^{\infty }\right) $, and $\mathbb{Z}(p^{k})$ for $k\in \omega +1$;

\item \emph{has finite ranks} if it has finite $\mathbb{Z}$-rank, $\mathbb{S}%
^{1}$-rank, and $p$-rank for every prime $p$.
\end{itemize}

We let $\mathbf{FLCPAb}$ be the thick subcategory of $\mathbf{LCPAb}$
consisting of locally compact abelian Polish groups with finite ranks. We
also let 
\begin{equation*}
\mathbf{TorFLCPAb}:=\mathbf{TorLCPAb}\cap \mathbf{FLCPAb}
\end{equation*}%
\begin{equation*}
\mathbf{FLCPAb}\left( p\right) :=\mathbf{LCPAb}\left( p\right) \cap \mathbf{%
FLCPAb}
\end{equation*}

Notice that $\mathbf{FLCPAb}\left( p\right) $ is an abelian category.
Indeed, if $G/H$ belongs to its left heart, then we can assume without loss
of generality that $H$ is countable, since it contains a compact open
subgroup; see also Lemma \ref{Lemma:presentation}. Furthermore, we have that 
$G\cong G_{0}\oplus G_{1}$ where $G_{0}$ is torsion-free and $G_{1}$ is
countable. Since $H$ is a countable $p$-group, we have that $H\subseteq
G_{1} $ and $G/H\cong G_{0}\oplus \left( G_{1}/H\right) \in \mathbf{FLCPAb}%
\left( p\right) $.

\subsection{The standard Borel category of locally compact Polish abelian
groups\label{Subsection:standard-Borel-category-lc}}

We outline how the category $\mathbf{LCPAb}$ of locally compact Polish
abelian groups can be seen as a Borel quasi-abelian category. For a Polish
space $X$ with a compatible complete metric $d$ bounded by $1$, let $\mathrm{%
F}\left( X\right) $ be the space of closed subsets of $X$ endowed with the
Effros Borel structure \cite[Section 12.C]{kechris_classical_1995}, which is
the $\sigma $-algebra generated by the sets of the form $\left\{ F\in 
\mathrm{F}\left( X\right) :F\cap U\neq \varnothing \right\} $ for some open
subset $U$ of $X$. This is also the $\sigma $-algebra associated with the
topology obtained by letting a net $\left( C_{i}\right) $ converge to $C$ if
and only if for every $x\in X$ we have that $d\left( C_{i},x\right)
\rightarrow d\left( C,x\right) $, where 
\begin{equation*}
d\left( F,x\right) =\mathrm{inf}\left\{ d\left( y,x\right) :y\in F\right\}
\end{equation*}%
for $x\in X$ and $F\in \mathrm{F}\left( X\right) $ \cite{beer_polish_1991}.
We also let $\mathrm{K}\left( X\right) $ be the space of compact subsets of $%
X$, which is a Borel subset of $\mathrm{F}\left( X\right) $. Finally, we let 
$\mathrm{F}_{\mathrm{lc}}\left( X\right) $ be the space of locally compact
closed subspaces of $X$. Recall that a \emph{regular closed subset} $F$ of a
Polish space $X$ is a set that is the closure of its interior. The
Kuratowski--Ryll-Nardzewski Theorem \cite[Theorem 12.13]%
{kechris_classical_1995} guarantees that there exists a Borel function $%
\mathrm{F}\left( X\right) \setminus \left\{ \varnothing \right\} \rightarrow
X^{\omega }$, $F\mapsto \left( x_{F,n}\right) _{n\in \omega }$ such that,
for every closed nonempty subset $F$ of $X$, $\left( x_{F,n}\right) _{n\in
\omega }$ is an enumeration of a dense subset of $F$.

\begin{lemma}
\label{Lemma:selection-locally-compact}We have that:

\begin{enumerate}
\item $\mathrm{F}_{\mathrm{lc}}\left( X\right) $ is a Borel subset of $%
\mathrm{F}\left( X\right) $;

\item there exists a Borel function $\mathrm{F}_{\mathrm{lc}}\left( X\right)
\rightarrow \mathrm{K}\left( X\right) ^{\omega }$, $F\mapsto \left(
F_{n}\right) $ such that the sequence $\left( F_{n}\right) $ is increasing, $%
F_{n}$ is a compact regular closed subset of $F$ for every $n\in \omega $,
and the relative interiors of the $F_{n}$ cover $F$; in particular, $F$ is
the union of the sequence $\left( F_{n}\right) $;

\item if $F\in \mathrm{F}_{\mathrm{lc}}\left( X\right) $, then the function $%
\mathrm{F}\left( X\right) \rightarrow \mathrm{F}_{\mathrm{lc}}\left(
X\right) $, $C\mapsto F\cap C$ is Borel.
\end{enumerate}
\end{lemma}

\begin{proof}
Fix a countable basis $\mathcal{B}$ of open sets for $X$, and a compatible
metric $d$ on $X$ bounded by $1$. Notice that for every open subset $U$ of $%
X $ the function $\mathrm{F}\left( X\right) \rightarrow \mathrm{F}\left(
X\right) $, $F\mapsto \overline{F\cap U}$ is Borel. Indeed, for any open set 
$W\subseteq X$, 
\begin{equation*}
\overline{F\cap U}\cap W\neq \varnothing \Leftrightarrow F\cap \left( U\cap
W\right) \neq \varnothing \text{.}
\end{equation*}

\begin{enumerate}
\item Notice that an element $F$ of $\mathrm{F}\left( X\right) $ is locally
compact if and only if there exist sequences $\left( K_{n}\right) $ in $%
\mathrm{K}\left( X\right) $ and $\left( D_{n}\right) $ in $\mathrm{F}\left(
X\right) $ such that, setting $U_{n}:=F\setminus D_{n}$, one has that%
\begin{equation*}
K_{n}\subseteq U_{n}\subseteq K_{n+1}
\end{equation*}%
and%
\begin{equation*}
\bigcup_{n}U_{n}=\bigcup_{n}K_{n}=F\text{.}
\end{equation*}%
This easily implies using the results of \cite[Section 12.C]%
{kechris_classical_1995} that $\mathrm{F}_{\mathrm{lc}}\left( X\right) $ is
analytic.

We also have that, for $F\in \mathrm{F}\left( X\right) $, $F$ is locally
compact if and only if for every $x\in X$, if $x\in F$ then there exists $%
U\in \mathcal{B}$ such that $x\in U$ and the closure $\overline{F\cap U}$ of 
$F\cap U$ is compact. This shows that $\mathrm{F}_{\mathrm{lc}}\left(
X\right) $ is also coanalytic.

\item If $F\subseteq X$ is locally compact, then we can find a sequence $%
\left( U_{n}\right) $ in $\mathcal{B}$ such that $F$ is the union of the
sequence $\left( \overline{F\cap U_{n}}\right) _{n\in \omega }$ and $%
\overline{F\cap U_{n}}$ is compact for every $n\in \omega $. For a given
open set $U$, the set $\left\{ F\in \mathrm{F}\left( X\right) :\overline{%
F\cap U}\in \mathrm{K}\left( X\right) \right\} $ is Borel. Thus, by the
\textquotedblleft small section\textquotedblright\ Uniformization Theorem 
\cite[Theorem 18.10]{kechris_classical_1995} there exist Borel functions $%
F\mapsto k_{F,n}$ for $n\in \omega $ such that $\left( U_{k_{F,n}}\right) $
is an enumeration of the elements $U$ of $\mathcal{B}$ such that $\overline{%
F\cap U}$ is compact. Each $\overline{F\cap U_{k_{F,n}}}$ is regular closed
in $F$, being the closure of a relatively open subset of $F$. Furthermore,
since $F$ is locally compact, every $x\in F$ belongs to $F\cap U$ for some $%
U\in \mathcal{B}$ with $\overline{F\cap U}$ compact, and $F\cap U$ is
relatively open in $F$; hence the relative interiors of the sets $\overline{%
F\cap U_{k_{F,n}}}$ cover $F$. As a finite union of compact regular closed
subsets of $F$ is again compact and regular closed, replacing $(\overline{%
F\cap U_{k_{F,n}}})_{n\in \omega }$ with the sequence of its successive
finite unions yields a function $F\mapsto \left( F_{n}\right) $ as required.

\item Suppose initially that $F\in \mathrm{K}\left( X\right) $. Let $\left\{
z_{n}:n\in \omega \right\} $ be an enumeration of a dense subset of $F$. We
need to prove that if $U\subseteq X$ is open, then $\left\{ C\in \mathrm{F}%
\left( X\right) :C\cap F\cap U\neq \varnothing \right\} $ is Borel. Without
loss of generality, we can assume that $U=\left\{ x\in X:d\left(
y_{0},x\right) <r\right\} $ for some $r>0$. Since $U$ is the union of $L_{s}$
for $s\in \mathbb{Q}\cap \left( 0,r\right) $ where $L_{s}=\left\{ x\in
X:d\left( y_{0},x\right) \leq s\right\} $, it suffices to prove that, for
every $s\in \mathbb{Q}\cap \left( 0,r\right) $, 
\begin{equation*}
\mathcal{W}_{F,s}:=\left\{ C\in \mathrm{F}\left( X\right) :C\cap F\cap
L_{s}\neq \varnothing \right\}
\end{equation*}%
is Borel. Since $F$ is compact, we have that, for $C\in \mathrm{F}\left(
X\right) $, $C\cap F\cap L_{s}\neq \varnothing $ if and only if for every $%
k\in \omega $ there exists $\ell \in \omega $ such that $d\left( z_{\ell
},C\right) <2^{-k}$ and $d\left( z_{\ell },y_{0}\right) \leq s$. Thus, $%
\mathcal{W}_{F,s}$ is Borel.

For arbitrary $F\in \mathrm{F}_{\mathrm{lc}}\left( X\right) $, we can write%
\begin{equation*}
F\cap C=\bigcup_{n}\left( F_{n}\cap C\right)
\end{equation*}%
where the Borel function $F\mapsto \left( F_{n}\right) $ is as in (2). This
concludes the proof.
\end{enumerate}
\end{proof}

By \cite[Theorem 2]{shkarin_universal_1999}, there exists a universal
abelian Polish group $\mathbb{G}$, i.e., an abelian Polish group that
contains an isomorphic copy of any other abelian Polish group; see \cite[%
Section 2.5]{gao_invariant_2009}. As in \cite[Section 2.6]%
{lupini_polish_2017}, the category $\mathbf{LCPAb}$ is equivalent to the
category $\mathrm{SG}_{\mathrm{lc}}\left( \mathbb{G}\right) $ of \emph{%
closed locally compact} subgroups of $\mathbb{G}$. This can be seen as the
Borel category whose arrows are the subsets $f$ of $\mathrm{F}_{\mathrm{lc}%
}\left( \mathbb{G}^{2}\right) $ that are the graphs of continuous group
homomorphisms between closed locally compact subgroups of $\mathbb{G}$. The
argument in \cite[Section 2.6]{lupini_polish_2017} shows that $\mathrm{SG}_{%
\mathrm{lc}}\left( \mathbb{G}\right) $ is indeed a Borel category, which is
easily seen to be, in fact, a Borel quasi-abelian category.

One can equivalently parametrize abelian Polish groups via \emph{countable
pseudometric abelian groups}, which are pairs $\left( \ell ,s\right) $ where 
$s:\mathbb{N}^{2}\rightarrow \mathbb{N}$ defines an abelian group operation
on $\mathbb{N}$ and $\ell $ is an invariant pseudonorm (or pseudolength
function) on the corresponding group $\left( \mathbb{N},s\right) $; see \cite%
[Section 1]{shkarin_universal_1999}. Allowing pseudonorms is harmless: one
first divides by the zero-distance subgroup and then takes the metric
completion. The triple $G_{\left( s,\ell \right) }:=\left( \mathbb{N},s,\ell
\right) $ thus defines a Polish abelian group $\hat{G}_{\left( s,\ell
\right) }$ with a distinguished dense sequence $(x_{n}^{(s,\ell )})$
corresponding to the elements of $\mathbb{N}$. Every abelian Polish group
arises in this fashion, and the convention also permits finite groups to be
represented using the fixed underlying set $\mathbb{N}$. A continuous group
homomorphism $\varphi :\hat{G}_{\left( s,\ell \right) }\rightarrow \hat{G}%
_{\left( s^{\prime },\ell ^{\prime }\right) }$ is coded by the relation $%
\rho _{\varphi }:\mathbb{N}\times \mathbb{N}\rightarrow \mathbb{R}$ defined
by setting 
\begin{equation*}
\rho _{\varphi }\left( n,m\right) :=d(\varphi (x_{n}^{(s,\ell
)}),x_{m}^{(s^{\prime },\ell ^{\prime })})
\end{equation*}%
for $n,m\in \omega $. These object codes and distance-matrix morphism codes
can be arranged in a Borel category $\mathbf{MAG}$ presenting abelian Polish
groups. Indeed, their defining conditions are countable combinations of
inequalities between coordinates, while Borel separation--completion is
provided by \cite[Proposition 6.2]{vaes_wouters_borel_2025}.

The codes with locally compact completion form a Borel full subcategory. To
see this directly, for a code $c=\left( s,\ell \right) $ put $%
d_{c}(n,m):=\ell \left( n^{-1}m\right) $. Its completion is locally compact
if and only if 
\begin{equation*}
\begin{aligned} &\exists r\in \mathbb{Q}_{>0}\quad \forall k\in \omega\quad
\exists N\in \omega\quad \forall n\in \omega\quad \\ &\qquad \left( \ell
\left( n\right) <r\ \Longrightarrow \exists i\leq N\ d_{c}\left( n,i\right)
<2^{-k}\right) . \end{aligned}
\end{equation*}%
This is a Borel condition, and it says precisely that some neighborhood of
zero in the dense coded subgroup is totally bounded, or equivalently that
the completion has a compact neighborhood of zero. We denote the resulting
full Borel subcategory by $\mathbf{MAG}_{\mathrm{lc}}$. Thus $\mathrm{SG}_{%
\mathrm{lc}}\left( \mathbb{G}\right) $ and $\mathbf{MAG}_{\mathrm{lc}}$ are
two Borel presentations of $\mathbf{LCPAb}$.

\begin{lemma}
\label{Lemma:Borel-classes}The isomorphism classes of locally compact Polish
abelian groups are Borel.
\end{lemma}

\begin{proof}
By \cite[Theorem 3.1]{rosendal_compact_2018}, the relation of isomorphism of
locally compact Polish abelian groups is Borel reducible to the orbit
equivalence relation $E$ of a continuous action of a Polish group on a
Polish space. Since $E$ has Borel equivalence classes \cite[Proposition
3.1.10]{gao_invariant_2009}, the conclusion follows.
\end{proof}

\begin{lemma}
\label{Lemma:Borel-connected-component}The functor $G\mapsto c\left(
G\right) $ assigning to a locally compact abelian Polish group its connected
component of zero is a Borel functor on $\mathrm{SG}\left( \mathbb{G}\right) 
$.
\end{lemma}

\begin{proof}
It suffices to notice that the relation on $\mathrm{Ob}\left( \mathrm{SG}%
\left( \mathbb{G}\right) \right) $ consisting of pairs $\left( G,H\right) $
such that $H\subseteq G$ is closed, connected, and $G/H$ is totally
disconnected is Borel by items (3) and (4) of Lemma \ref%
{Lemma:Borel-categories} below.
\end{proof}

Recall that if $\mathcal{C}$ is a full subcategory of the category of Polish
abelian groups, then a subfunctor of the identity on $\mathcal{C}$ is a
functor $F:\mathcal{C}\rightarrow \mathcal{C}$ such that $F\left( G\right)
\subseteq G$ and $F\left( f\right) =f|_{F\left( G\right) }$ for every pair
of objects $G,H$ of $\mathcal{C}$ and continuous group homomorphism $%
f:G\rightarrow H$.

\begin{lemma}
\label{Lemma:Borel-types}The functors $G\mapsto G_{\mathbb{S}^{1}}$ and $%
G\mapsto F_{\mathbb{Z}}G$ assigning to a locally compact Polish abelian
group $G$ closed subgroups such that $G_{\mathbb{S}^{1}}$, $G_{\mathbb{A}%
}:=F_{\mathbb{Z}}G/G_{\mathbb{S}^{1}}$ and $G_{\mathbb{Z}}:=G/F_{\mathbb{Z}%
}G $ are of type $\mathbb{S}^{1}$, $\mathbb{A}$, and $\mathbb{Z}$
respectively are Borel.
\end{lemma}

\begin{proof}
This follows from items (7), (8), and (9) of Lemma \ref%
{Lemma:Borel-categories} together with \cite[Proposition 2.2]%
{hoffmann_homological_2007}, considering the uniqueness of the chain of
subgroups $G_{\mathbb{S}^{1}}\subseteq F_{\mathbb{Z}}G\subseteq G$ subject
to the requirement that $G_{\mathbb{S}^{1}}$, $G_{\mathbb{A}}:=F_{\mathbb{Z}%
}G/G_{\mathbb{S}^{1}}$ and $G_{\mathbb{Z}}:=G/F_{\mathbb{Z}}G$ are of type $%
\mathbb{S}^{1}$, $\mathbb{A}$, and $\mathbb{Z}$ respectively.
\end{proof}

\begin{lemma}
\label{Lemma:Borel-categories}The following classes of locally compact
Polish abelian groups are Borel:

\begin{enumerate}
\item countable groups;

\item compact groups;

\item connected groups;

\item totally disconnected groups;

\item topological $p$-groups;

\item topological torsion groups;

\item groups of type $\mathbb{S}^{1}$;

\item groups of type $\mathbb{A}$;

\item groups of type $\mathbb{Z}$;

\item divisible groups;

\item profinite groups;

\item pro-$p$ groups;

\item abelian Lie groups;

\item compactly generated groups.
\end{enumerate}
\end{lemma}

\begin{proof}
Fix a countable basis $\mathcal{B}$ of open sets for $\mathbb{G}$, and a
compatible length function $\ell $ on $\mathbb{G}$ bounded by $1$.

(1) Let $\left( x_{n}^{G}\right) $ be an enumeration of a dense subset of $G$
given by the Kuratowski--Ryll-Nardzewski Theorem cited above; then $G$ is
countable if and only if there exists $U\in \mathcal{B}$ such that $0\in U$
and $\forall n\in \omega $, $x_{n}^{G}\in U\Rightarrow x_{n}^{G}=0$.

(2) This follows from the fact that the collection $K\left( \mathbb{G}%
\right) $ of compact subsets of $\mathbb{G}$ is Borel.

(3) By definition, a locally compact abelian Polish group $G$ is connected
if and only if for every pair of closed subsets $F,F^{\prime }$ of $G$ with $%
F\cap F^{\prime }=\varnothing $ and $F\cup F^{\prime }=G$ we have that $%
F=\varnothing $ or $F^{\prime }=\varnothing $. This shows that the class of
locally compact connected abelian Polish groups is coanalytic.

By \cite[Theorem 4.8]{armacost_structure_1981}, a locally compact abelian
Polish group $G$ is connected if and only if its $\mathbb{R}$-subgroups
generate $G$. Thus, a locally compact abelian Polish group $G$ with dense
subset $\left( g_{n}\right) $ and compatible metric $d_{G}$ is connected if
and only if for every $n,k\in \omega $ there exists a homomorphism $\phi :%
\mathbb{R}\rightarrow G$ and $t\in \mathbb{R}$ such that $d_{G}\left( \phi
\left( t\right) ,g_{n}\right) <2^{-k}$. This shows that the class of locally
compact connected abelian Polish groups is analytic.

(4) We have that a locally compact Polish abelian group $G$ with basis of
zero neighborhoods $\left( U_{n}\right) $ is totally disconnected if and
only if it is zero-dimensional, if and only if for every $n$ there exists an
open subgroup $H$ of $G$ contained in $U_{n}$ \cite[Theorem 7.7]%
{hewitt_abstract_1979}.\ This shows that the class of totally disconnected
locally compact Polish abelian groups is analytic. We also have that $G$ is
totally disconnected if and only if $\mathrm{Hom}\left( \mathbb{R},G\right)
=0$, if and only if $\forall \varphi \in \mathrm{Hom}\left( \mathbb{R}%
,G\right) $, $\forall t\in \mathbb{R}$, $\varphi \left( t\right) =0$.\ This
shows that the class of totally disconnected groups is also coanalytic.

(5) Notice that the class of \emph{countable} $p$-groups is Borel. By \cite[%
Theorem 2.12]{armacost_structure_1981} we have that a locally compact
abelian group $G$ with basis of zero neighborhoods $\left( U_{n}\right) $ is
a topological $p$-group if and only if for every $n\in \omega $ there exists
a clopen subgroup $H$ of $G$ contained in $U_{n}$ such that $G/H$ is a
countable $p$-group. This shows that the class of locally compact abelian
Polish topological $p$-groups is analytic. We also have that $G$ is a
topological $p$-group if and only if for every $x\in G$ and $\varepsilon >0$
there exists $n_{0}\in \mathbb{N}$ such that for every $n\geq n_{0}$, $\ell
\left( p^{n}x\right) <\varepsilon $. This shows that the class of
topological $p$-groups is also coanalytic.

(6) Notice that the class of \emph{countable }torsion groups is Borel. By 
\cite[Theorem 3.5]{armacost_structure_1981}, we have that a locally compact
abelian group $G$ with basis of zero neighborhoods $\left( U_{n}\right) $ is
a topological torsion group if and only if for every $n\in \omega $ there
exists a clopen subgroup $H$ of $G$ contained in $U_{n}$ such that $G/H$ is
a countable torsion group. This shows that the class of locally compact
abelian Polish topological torsion groups is analytic. We also have that $G$
is a topological torsion group if and only if for every $x\in G$ and $%
\varepsilon >0$ there exists $n_{0}\in \mathbb{N}$ such that for every $%
n\geq n_{0}$, $\ell \left( n!x\right) <\varepsilon $. This shows that the
class of topological torsion groups is also coanalytic.

(7) By definition, a locally compact abelian Polish group is of type $%
\mathbb{S}^{1}$ if and only if it is compact and connected. The conclusion
thus follows from (2) and (3).

(8) We have that a locally compact abelian Polish group $G$ is of type $%
\mathbb{A}$ if and only if $c\left( G\right) \cong \mathbb{R}^{n}$ for some $%
n\in \omega $ and $G/c\left( G\right) $ is a topological torsion group. This
item thus follows from Lemma \ref{Lemma:Borel-classes}, Lemma \ref%
{Lemma:Borel-connected-component}, and (6).

(9) A locally compact abelian Polish group $G$ is of type $\mathbb{Z}$ if
and only if it is countable and torsion-free. This item thus follows from
(1).

(10) We have that $G$ is divisible if and only if for every $k\in \mathbb{Z}$%
, $kG=G$, if and only if $kG$ is dense in $G$, which is easily seen to be a
Borel condition.

(11) A group is profinite if and only if it is compact and topological
torsion. This item thus follows from (2) and (6).

(12) A group is pro-$p$ if and only if it is a compact topological $p$%
-group. This item thus follows from (2) and (5).

(13) An abelian locally compact Polish group is a Lie group if and only if $%
c\left( G\right) $ is a Lie group, which happens if and only if $c\left(
G\right) \cong \mathbb{R}^{n}\oplus \mathbb{T}^{k}$ for some $n,k\in \omega $%
. This item thus follows from Lemma \ref{Lemma:Borel-classes}.

(14) We have that $G$ is compactly generated if and only if $G\cong C\oplus 
\mathbb{R}^{n}\oplus \mathbb{Z}^{k}$ for $n,k\in \omega $, if and only if $%
G_{\mathbb{Z}}$ has finite rank and $G_{\mathrm{t}}:=G_{\mathbb{A}}/c\left(
G_{\mathbb{A}}\right) $ is compact. This item thus follows from Lemma \ref%
{Lemma:Borel-types} together with (2).
\end{proof}

\begin{remark}
\label{Remark:tdlc}A number of equivalent Borel presentations of the
category of \emph{totally disconnected }locally compact Polish abelian
groups can be considered as in \cite%
{nies_coarse_2022,lupini_computable_2023,melnikov_computably_2023,kechris_complexity_2018}%
. The methods in \cite{lupini_computable_2023} give an alternative, more
concrete coding of the restriction of Pontryagin duality to this
subcategory. Proposition \ref{Proposition:Borel-dual} below gives a uniform
coding which applies without the totally disconnected hypothesis.
\end{remark}

\begin{proposition}
\label{Proposition:Borel-dual}Pontryagin duality admits a Borel
presentation. More precisely, there is a Borel functor 
\begin{equation*}
D:\mathrm{SG}_{\mathrm{lc}}( \mathbb{G}) ^{\mathrm{op}%
}\longrightarrow \mathbf{MAG}_{\mathrm{lc}}
\end{equation*}%
whose realization on completions is naturally isomorphic to the functor $%
( -) ^{\vee }: \mathbf{LCPAb}^{\mathrm{op}}\rightarrow \mathbf{%
LCPAb}$. Thus both $G\mapsto G^{\vee }$ on objects and $f\mapsto f^{\vee }$
on arrows have Borel codes. In particular, Pontryagin duality is represented
by a Borel functor that is an anti-equivalence on the underlying categories.
\end{proposition}

\begin{proof}
Put $X:=\mathrm{Ob}( \mathrm{SG}_{\mathrm{lc}}( \mathbb{G})
) $ and consider the tautological field 
\begin{equation*}
\mathscr{G}:=\{ ( G,x) :G\in X,\ x\in G\}
\longrightarrow X\text{.}
\end{equation*}%
With the Borel structure inherited from $X\times \mathbb{G}$, the metric
induced by a fixed compatible complete metric on $\mathbb{G}$, and Borel
dense sections supplied by the Kuratowski--Ryll-Nardzewski theorem \cite[
Theorem 12.13]{kechris_classical_1995}, $\mathscr{G}$ is a Borel field of locally compact Polish abelian groups. By 
\cite[Proposition 10.1]{vaes_wouters_borel_2025}, the field 
\begin{equation*}
\widehat{\mathscr{G}}:=\bigsqcup_{G\in X}G^{\vee }
\end{equation*}%
of Pontryagin duals has a unique standard Borel structure for which the
projection to $X$ and the evaluation map 
\begin{equation}  \label{Equation:Borel-dual-evaluation}
\widehat{\mathscr{G}}\times _{X}\mathscr{G}\longrightarrow \mathbb{T},
\qquad ( \chi ,x) \longmapsto \chi ( x)
\end{equation}%
are Borel. With this structure, $\widehat{\mathscr{G}}$ is a Borel field of
locally compact Polish abelian groups. In particular, we can choose Borel
sections $G\mapsto \chi _{i}^{G}$, for $i\in \omega $, whose values are
dense in $G^{\vee }$ for the compact-open topology.

We record a convenient Borel metric on this field. Fix a bounded compatible
invariant metric $d_{\mathbb{T}}$ on $\mathbb{T}$ and Borel dense sections $%
G\mapsto x_{n}^{G}$ of $\mathscr{G}$. By Lemma \ref%
{Lemma:selection-locally-compact}, we obtain, in a Borel fashion in $G$, an
increasing sequence $( K_{j}^{G}) _{j\in \omega }$ of compact
regular closed subsets of $G$ whose interiors cover $G$. Define 
\begin{equation}
\partial _{G}( \chi ,\psi ) :=\sum_{j\in \omega }2^{-j-1}\sup
\{ d_{\mathbb{T}}(\chi (x_{n}^{G}),\psi (x_{n}^{G})):x_{n}^{G}\in
K_{j}^{G}\} .  \label{Equation:Borel-dual-metric}
\end{equation}%
where the supremum of the empty set is understood to be $0$. Since $K_{j}^{G}
$ is regular closed, the points $x_{n}^{G}$ that belong to $K_{j}^{G}$ are
dense in $K_{j}^{G}$. Moreover, since the $K_{j}^{G}$ are increasing and
their interiors cover $G$, every compact subset of $G$ is contained in some $%
K_{j}^{G}$. Hence \eqref{Equation:Borel-dual-metric} induces the
compact-open topology. A $\partial _{G}$-Cauchy sequence is uniformly Cauchy
on every $K_{j}^{G}$; its compatible uniform limits define a continuous
character, and convergence to this character in $\partial _{G}$ follows by
controlling the tail of the series. Thus $\partial _{G}$ is a complete
invariant metric. It is Borel on $\widehat{\mathscr{G}}\times _{X}\widehat{%
\mathscr{G}}$, since evaluation in \eqref{Equation:Borel-dual-evaluation}
and the relation $x_{n}^{G}\in K_{j}^{G}$ are Borel.

Fix the countable group $A:=\mathbb{Z}^{(\omega )}$ and an enumeration of $A$
which transports its group operation to $\mathbb{N}$. For $a\in A$, set 
\begin{equation*}
\chi _{a}^{G}:=\sum_{i\in \omega }a(i)\chi _{i}^{G},\qquad \ell _{G}^{\vee
}(a):=\partial _{G}( \chi _{a}^{G},0) ,
\end{equation*}%
where the sum is finite. The map $( G,a) \mapsto \ell _{G}^{\vee
}(a)$ is Borel and $\ell _{G}^{\vee }$ is an invariant pseudonorm on $A$.
The homomorphism $A\rightarrow G^{\vee }$, $a\mapsto \chi _{a}^{G}$, has
dense range and kernel $\{ a:\ell _{G}^{\vee }(a)=0\} $.
Consequently, the completion of its metric quotient is canonically $G^{\vee
} $. Thus $G\mapsto D(G):=( A,\ell _{G}^{\vee }) $ is a Borel
object map from $\mathrm{SG}_{\mathrm{lc}}( \mathbb{G}) $ to $%
\mathbf{MAG}_{\mathrm{lc}}$.

We next consider arrows. If $f:G\rightarrow H$ and $\chi \in H^{\vee }$, put 
\begin{equation*}
P(f,\chi ):=\chi \circ f\in G^{\vee }.
\end{equation*}%
The map $P$ is Borel. Indeed, evaluation of a continuous map coded by its
graph at a Borel point is Borel, since the corresponding Borel relation has
singleton sections. Hence, for every $n\in \omega $, 
\begin{equation*}
P(f,\chi )( x_{n}^{G}) =\chi ( f( x_{n}^{G})
)
\end{equation*}%
is Borel in $( f,\chi ) $. The map 
\begin{equation*}
\widehat{\mathscr{G}}\longrightarrow X\times \mathbb{T}^{\omega }, \qquad
\chi \longmapsto ( G,( \chi ( x_{n}^{G}) ) _{n\in
\omega }) \quad ( \chi \in G^{\vee })
\end{equation*}%
is Borel and injective. By the Lusin--Souslin theorem it is a Borel
isomorphism onto its image, so the displayed coordinate identities imply
that $P$ is Borel. Notice also that $P(f,-):H^{\vee }\rightarrow G^{\vee }$
is continuous for the compact-open topologies, since $f$ maps compact sets
to compact sets.

Relative to the dense subgroups coded above, the morphism $%
D(f):D(H)\rightarrow D(G)$ is represented by the matrix 
\begin{equation*}
\rho _{D(f)}(a,b):=\partial _{G}( P( f,\chi _{a}^{H}) ,\chi
_{b}^{G}) ,\qquad a,b\in A.
\end{equation*}%
This matrix depends in a Borel fashion on $f$ and is unchanged upon passage
to the zero-distance quotients. Hence $f\mapsto D(f)$ is Borel. Since
precomposition satisfies $D(1_{G})=1_{D(G)}$ and $D(g\circ f)=D(f)\circ D(g)$%
, $D$ is a Borel functor on $\mathrm{SG}_{\mathrm{lc}}( \mathbb{G}%
) ^{\mathrm{op}}$.

The completion represented by $D(G)$ is canonically $G^{\vee }$, and under
this identification $D(f)$ is precomposition by $f$. Hence the realization
of $D$ is naturally isomorphic to the ordinary Pontryagin-dual functor. By
the Pontryagin duality theorem, this functor is an anti-equivalence of the
underlying categories. Thus $D$ is a Borel presentation of the Pontryagin
anti-equivalence, as claimed.
\end{proof}

If $G$ is a topological torsion locally compact abelian Polish group, then
its $p$-component is the closed subgroup $G_{p}$ consisting of $x\in G$ such
that $\lim_{n\rightarrow \infty }p^{n}x=0$; see \cite[Lemma 3.8]%
{armacost_structure_1981}.\ We also have that $G_{p}$ is the set of elements
of $G$ that belong to the range of a continuous group homomorphism $\mathbb{Z%
}_{p}\rightarrow G$; see \cite[Example 4.13(a)]{armacost_structure_1981}.

\begin{lemma}
\label{Lemma:Borel-p-component}The functor $G\mapsto G_{p}$ assigning to a
topological torsion locally compact abelian Polish group its $p$-component
is Borel.
\end{lemma}

\begin{proof}
Let $\mathbf{TorLCPAb}$ be the category of topological torsion locally
compact abelian Polish groups. Just as for the Pontryagin dual functor, the
functor $G\mapsto \mathrm{Hom}\left( G,\mathbb{Z}\left( p^{\infty }\right)
\right) $ on $\mathbf{TorLCPAb}$, where $\mathbb{Z}\left( p^{\infty }\right)
\subseteq \mathbb{T}$ is the group of roots of unity of $p$-power order, is
Borel. It follows that $\left\{ G:G_{p}=0\right\} =\left\{ G:\mathrm{Hom}%
\left( G,\mathbb{Z}\left( p^{\infty }\right) \right) =0\right\} $ is a Borel
subset of $\mathbf{TorLCPAb}$. (Notice that for a topological $p$-group $A$,
every continuous homomorphism $A\rightarrow \mathbb{T}$ has values in $%
\mathbb{Z}\left( p^{\infty }\right) $, and $A^{\vee }=\mathrm{Hom}\left( A,%
\mathbb{Z}\left( p^{\infty }\right) \right) $ is also a topological $p$%
-group; see \cite[Theorem 2.12 and Corollary 2.13]{armacost_structure_1981}%
). Since $G_{p}$ is the unique closed subgroup of $G $ that is a topological 
$p$-group and such that $G/G_{p}$ has trivial $p$-component, the conclusion
follows.
\end{proof}

Recall the notion of a locally compact Polish abelian group \emph{with
finite ranks} \cite[Definitions 2.5 and 2.6]{hoffmann_homological_2007}.

\begin{lemma}
\label{Lemma:Borel-finite-ranks}The class of locally compact Polish abelian
groups with finite ranks is Borel.
\end{lemma}

\begin{proof}
By Lemma \ref{Lemma:Borel-types} and \cite[Proposition 2.9]%
{hoffmann_homological_2007} it suffices to consider the cases of groups of
type $\mathbb{Z}$, $\mathbb{S}^{1}$, and $\mathbb{A}$. For groups of type $%
\mathbb{Z}$, we have that a countable torsion-free group has finite ranks if
and only if it has finite $\mathbb{Z}$-rank, if and only if it has finite
torsion-free rank, which is easily seen to be a Borel condition. The case of
groups of type $\mathbb{S}^{1}$ follows by duality in view of Proposition %
\ref{Proposition:Borel-dual}.

We now consider the case of groups of type $\mathbb{A}$, which we can take
to be topological torsion groups. For such a group $G$, we have that $G$ has
finite ranks if and only if $G_{p}$ has finite $p$-rank for every prime $p$.
It therefore suffices to consider the case of topological $p$-groups for a
given prime $p$. In this case the conclusion follows from Lemma \ref%
{Lemma:Borel-classes}, noticing that there exist only countably many classes
of topological $p$-groups with finite ranks by \cite[Lemma 2.8]%
{hoffmann_homological_2007}.
\end{proof}

\begin{lemma}
\label{Lemma:connected-is-divisible}Let $C$ be a locally compact Polish
abelian group. If $C$ is connected, then $C$ is divisible.
\end{lemma}

\begin{proof}
If $C$ is connected, then $C\cong K\oplus \mathbb{R}^{n}$ where $K$ is
compact and connected. Thus, without loss of generality we can assume that $%
C $ is compact. Then for every prime number $p$, we have that $C/pC$ is a
compact bounded topological $p$-group. (Recall that an abelian group $G$ is
bounded if there exists $n\in \mathbb{Z}$ such that $nx=0$ for every $x\in G$%
.) Thus, $\left( C/pC\right) ^{\vee }$ is a countable bounded topological $p$%
-group. Thus we have that $\left( C/pC\right) ^{\vee }\cong \left( \mathbb{Z}%
/p\mathbb{Z}\right) ^{\left( n\right) }$ for some $n\in \omega +1$. Thus we
have that $C/pC\cong \left( \mathbb{Z}/p\mathbb{Z}\right) ^{n}$.

For $i<n$ consider the projection $\left( \mathbb{Z}/p\mathbb{Z}\right)
^{n}\rightarrow \left( \mathbb{Z}/p\mathbb{Z}\right) ^{i}$ on the first $i$
coordinates. Since $\left( \mathbb{Z}/p\mathbb{Z}\right) ^{i}$ is finite, we
have that the composition 
\begin{equation*}
C\rightarrow C/pC\cong \left( \mathbb{Z}/p\mathbb{Z}\right) ^{n}\rightarrow
\left( \mathbb{Z}/p\mathbb{Z}\right) ^{i}
\end{equation*}%
is trivial. Thus, we have that%
\begin{equation*}
pC=\mathrm{Ker}\left( C\rightarrow C/pC\right) =\bigcap_{i<n}\mathrm{Ker}%
(C\rightarrow \left( \mathbb{Z}/p\mathbb{Z}\right) ^{i})=C
\end{equation*}%
This shows that $pC=C$, and hence that $C$ is divisible.
\end{proof}

\begin{lemma}
\label{Lemma:divisible-types}A locally compact Polish abelian group $C$ is
divisible if and only if $C_{\mathrm{t}}$ and $C_{\mathbb{Z}}$ are divisible.
\end{lemma}

\begin{proof}
If $C$ is divisible, then $C_{\mathbb{Z}}$ is divisible, being a quotient of 
$C$. Since $C_{\mathbb{Z}}$ is torsion-free, we have that $F_{\mathbb{Z}}C=%
\mathrm{ker}\left( C\rightarrow C_{\mathbb{Z}}\right) $ is divisible.
Therefore, we have that $C_{\mathbb{A}}=C_{\mathbb{R}}\oplus C_{\mathrm{t}}$
is also divisible, and $C_{\mathrm{t}}$ is divisible.

Conversely, suppose that $C_{\mathrm{t}}$ and $C_{\mathbb{Z}}$ are
divisible. Then we have that $C_{\mathbb{S}^{1}}$ is divisible, being
connected, and $C_{\mathbb{A}}=C_{\mathrm{t}}\oplus C_{\mathbb{R}}$ is
divisible. Thus, $C$ is divisible, being obtained from divisible groups by
taking extensions.
\end{proof}

A locally compact Polish abelian group is called \emph{codivisible }if its
Pontryagin dual is divisible.

\begin{corollary}
\label{Corollary:codivisible-types}A Polish abelian group $C$ is codivisible
if and only if $C_{\mathrm{t}}$ and $C_{\mathbb{S}^{1}}$ are codivisible.
\end{corollary}

\begin{lemma}
The class of codivisible locally compact Polish abelian groups is Borel.
\end{lemma}

\begin{proof}
By Corollary \ref{Corollary:codivisible-types} and Lemma \ref%
{Lemma:Borel-types} it suffices to observe that the classes of topological
torsion codivisible groups and compact connected codivisible groups are
Borel. By Proposition \ref{Proposition:Borel-dual}, this follows from Lemma %
\ref{Lemma:Borel-categories}(10).
\end{proof}

\section{Extensions and cocycles}

\label{Section:derived}

{The study of locally compact abelian groups by the modern methods of
homological algebra was pioneered by Hoffmann and Spitzweck \cite%
{hoffmann_homological_2007}. There, the authors manage to produce for the
functor }$\mathrm{Hom}$ on this category a cohomological right derived
functor \cite%
{mac_lane_homology_1995,gelfand_methods_2003,iversen_cohomology_1986},
through a small tour-de-force involving their type decomposition of locally
compact abelian groups. An analogous result for the category of pro-Lie
Polish abelian groups was obtained in \cite{casarosa_homological_2026} using
projective resolutions. In this section, we restrict the analysis to the
\textquotedblleft well-behaved\textquotedblright\ category of those locally
compact abelian groups that are \emph{Polish }(i.e., second countable, as
most compact abelian groups arising in practice are). For such groups, we
obtain descriptions of the cohomological derived functor $\mathrm{Ext}$ in
terms of cocycles, in the spirit of Moore's measurable group cohomology \cite%
{moore_group_1976,moore_group_1976-1,moore_group_1968}.

\subsection{Resolutions for locally compact groups}

\label{Subsection:derived-lc}

We let $\mathbf{LCPAb}$ be the quasi-abelian category of locally compact
abelian Polish groups. Injective and projective objects in this
quasi-abelian category were characterized by Moskowitz in \cite%
{moskowitz_homological_1967}.

\begin{proposition}
\label{Proposition:injectives}Suppose that $I\in \mathbf{LCPAb}$. The
following assertions are equivalent:

\begin{enumerate}
\item $I$ is injective for $\mathbf{LCPAb}$;

\item $\mathrm{Ext}\left( \mathbb{T},I\right) =0$;

\item $I\cong V\oplus T$ where $V$ is a vector group and $T$ is a torus;

\item $I$ is path-connected.
\end{enumerate}

Dually, for $P\in \mathbf{LCPAb}$, the following assertions are equivalent:

\begin{enumerate}
\item $P$ is projective for $\mathbf{LCPAb}$;

\item $\mathrm{Ext}\left( P,\mathbb{Z}\right) =0$;

\item $P\cong V\oplus F$ where $V$ is a vector group and $F$ is a countable
free abelian group.
\end{enumerate}
\end{proposition}

\begin{proof}
We prove the first assertion, the second one following by duality.

The implications (1)$\Rightarrow $(2) and (3)$\Rightarrow $(4) are obvious.
The implication (3)$\Rightarrow $(1) is established in \cite[Theorem 3.2]%
{moskowitz_homological_1967}. The implication (4)$\Rightarrow $(3) is a
consequence of \cite[Theorem 8.27 and Remark 8.25(c)]%
{armacost_structure_1981}. We postpone the proof of the implication (2)$%
\Rightarrow $(3) to the end of this subsection.
\end{proof}

\begin{corollary}
\label{Corollary:quotient-injective}The collection of injective locally
compact Polish abelian groups is closed under taking quotients by closed
subgroups.
\end{corollary}

The total derived functor of $\mathrm{Hom}$ on the category $\mathbf{%
ProLiePAb}$ of \emph{pro-Lie }Polish abelian groups was produced in \cite[%
Section 3]{casarosa_homological_2026} by way of projective resolutions.
Since $\mathbf{ProLiePAb}$ is a \emph{hereditary }quasi-abelian category, $%
\mathrm{Ext}^{n}=0$ for $n\geq 2$ on $\mathbf{ProLiePAb}$; see also \cite[%
Theorem 2.9]{fulp_extensions_1971} and {\cite[Proposition 4.15(iv)]%
{hoffmann_homological_2007}}. For \emph{locally compact }Polish abelian
groups $C$ and $A$, $\mathrm{Hom}\left( C,A\right) $ is a (pro-Lie) Polish
abelian group. In the next subsection, we will show that one can analogously
regard $\mathrm{Ext}\left( C,A\right) $ as a group with a Polish cover.

\begin{remark}
\label{Remark:ext-dual}By \cite[Theorem 2.1]{moskowitz_homological_1967}, if 
$A\rightarrow B\rightarrow C$ is a short exact sequence of locally compact
Polish abelian groups, then $C^{\vee }\rightarrow B^{\vee }\rightarrow
A^{\vee }$ is a short exact sequence of Polish abelian groups. This
correspondence establishes a definable isomorphism between $\mathrm{Ext}%
\left( C,A\right) $ and $\mathrm{Ext}\left( A^{\vee },C^{\vee }\right) $.
\end{remark}

We now conclude the proof of Proposition \ref{Proposition:injectives}:

Suppose that $I$ satisfies (2) as in the first circuit of equivalences in
Proposition \ref{Proposition:injectives}. We have an exact sequence%
\begin{equation*}
\mathrm{Hom}\left( \mathbb{R},I\right) \rightarrow \mathrm{Hom}\left( 
\mathbb{Z},I\right) \cong I\rightarrow \mathrm{Ext}\left( \mathbb{T}%
,I\right) =0\text{.}
\end{equation*}%
Thus, every element of $I$ belongs to a one-parameter subgroup of $I$ (which
is a continuous group homomorphism $\mathbb{R}\rightarrow I$), and $I$ is
connected. Hence, $I\cong V\oplus C$ where $V$ is a vector group and $C$ is
compact and connected. Since $C^{\vee }$ is a countable torsion-free group
and 
\begin{equation*}
0=\mathrm{Ext}\left( \mathbb{T},C\right) \cong \mathrm{Ext}\left( C^{\vee },%
\mathbb{Z}\right)
\end{equation*}%
we have that $C^{\vee }$ is a free abelian group by \cite[Theorem 3.2]%
{friedenberg_extensions_2013}. Hence, $C$ is a torus.

\subsection{\textrm{Ext }and cocycles\label{Subsection:cocycles}}

In this subsection we present an equivalent description in terms of cocycles
of the groups $\mathrm{Hom}(G,A)$ and $\mathrm{Ext}\left( B,C\right) $ when $%
G,B,C$ are locally compact Polish abelian groups and $A$ is an abelian group
with a Polish cover. For this purpose, we recall the definition of the space 
$L^{0}$ of almost-everywhere equivalence classes of Borel functions.

By a \emph{$\sigma $-finite measure space} we mean a standard Borel space
endowed with a $\sigma $-finite Borel measure. Let $\left( X,\mu \right) $
be a $\sigma $-finite measure space and $A$ be a Polish space. One sets $%
L^{0}\left( X;A\right) $ (denoted by $U\left( X,A\right) $ in \cite%
{moore_group_1976}) to be the space of Borel functions $X\rightarrow A$ up
to the relation of equality almost-everywhere. $L^{0}\left( X;A\right) $ is
then a Polish space when endowed with the topology of convergence in $\nu $%
-measure, with respect to any Borel probability measure $\nu $ on $X$ that
is equivalent to $\mu $ and any compatible metric $d$ on $A$ that is bounded
by $1$. Thus, a net $\left( f_{i}\right) $ converges to $f$ if and only if 
\begin{equation*}
\nu \left( \left\{ x\in X:d\left( f_{i}\left( x\right) ,f\left( x\right)
\right) >\varepsilon \right\} \right) \rightarrow 0
\end{equation*}%
for every $\varepsilon >0$. A sequence $\left( f_{n}\right) $ converges to $%
f $ if and only if every subsequence $\left( f_{n_{k}}\right) $ has a
subsequence converging to $f$ almost-everywhere. (This shows that the
topology is independent of the choices of $\nu $ and $d$.) A compatible
complete metric on $L^{0}\left( X;A\right) $ is given by%
\begin{equation*}
d\left( f,g\right) :=\int_{X}d\left( f\left( x\right) ,g\left( x\right)
\right) d\nu \left( x\right) \text{;}
\end{equation*}%
see \cite[Proposition 6]{moore_group_1976}. The set of (equivalence classes
of) Borel functions $X\rightarrow A$ that attain finitely many values is
dense in $L^{0}\left( X;A\right) $. (When $X$ is purely atomic, this is easy
to see. One can thus assume without loss of generality that $X=[0,1]$ with
the Lebesgue measure, and then apply Luzin's Theorem on measurable functions 
\cite[Theorem 17.12]{kechris_classical_1995}.)

The Borel structure on $L^{0}\left( X;A\right) $ is generated by the
real-valued functions%
\begin{equation*}
f\mapsto \int_{x\in Y}\phi \left( f\left( x\right) \right) d\nu \left(
x\right)
\end{equation*}%
where $Y\subseteq X$ is a Borel set and $\phi :A\rightarrow \mathbb{R}$ is a 
\emph{bounded} \emph{Borel function} \cite[Proposition 8]{moore_group_1976}.
A Borel function $\phi :A\rightarrow B$ between Polish spaces induces a
Borel function $L^{0}\left( X;A\right) \rightarrow L^{0}\left( X;B\right) $,
which is continuous if $\phi $ is continuous.

In the statement and proof of the following results, we will use the measure
quantifiers $\forall ^{\ast }x$ and $\exists ^{\ast }x$. When the variable $%
x $ is ranging within a $\sigma $-finite measure space, $\forall ^{\ast }x$, 
$P\left( x\right) $ asserts that the property $P$ is satisfied by a conull
set of elements, while $\exists ^{\ast }x$, $P\left( x\right) $ asserts that
the property $P$ is satisfied by a non-null set of elements. We will also
regard a locally compact Polish abelian group as a $\sigma $-finite measure
space with respect to a Haar measure.

Recall that if $A$ and $G$ are abelian groups, then a (symmetric) $2$%
-cocycle on $G$ with values in $A$ is a function $c:G^{2}\rightarrow A$
satisfying, for every $x,y,z\in G$:

\begin{enumerate}
\item $c\left( x,y\right) =c\left( y,x\right) $;

\item $c\left( x,y\right) +c\left( x+y,z\right) =c\left( x,y+z\right)
+c\left( y,z\right) $.
\end{enumerate}

It is furthermore \emph{normalized }if $c\left( x,0\right) =0$ for every $%
x\in G$. We say that $c$ is Borel or continuous if it is so as a function $%
G^{2}\rightarrow A$. For a function $f:G\rightarrow H$ between abelian
groups, we define $\delta f:G^{2}\rightarrow H$ by setting $\delta f\left(
x,y\right) :=f(x)+f(y)-f(x+y)$. For a measurable function $f$, we let $[f]$
be its a.e.-class.

\begin{lemma}
\label{Lemma:Weil-discrete}Suppose that $G$ is a locally compact Polish
abelian group, $D$ is a (not necessarily countable) abelian group, and $%
f:G\rightarrow D$ is a function satisfying $\forall ^{\ast }x$, $\forall
^{\ast }y$, $f\left( x+y\right) =f\left( x\right) +f\left( y\right) $. Then
there exists a unique homomorphism $\varphi :G\rightarrow D$ such that $%
\forall ^{\ast }x$, $\varphi \left( x\right) =f\left( x\right) $.
\end{lemma}

\begin{proof}
Define 
\begin{equation*}
R:=\left\{ \left( x,y\right) \in G^{2}:f\left( x+y\right) =f\left( x\right)
+f\left( y\right) \right\} \subseteq G^{2}\text{.}
\end{equation*}%
Let 
\begin{equation*}
U:=\left\{ x\in G:\forall ^{\ast }y,\left( x,y\right) \in R\right\} \text{,}
\end{equation*}%
which is a conull subset of $G$.

Then for $x,y\in U$ with $x+y\in U$, $\left( x,y\right) \in R$. Indeed, $%
\forall ^{\ast }z\in G$ we have that 
\begin{equation*}
\left( x+y,z\right) ,\left( y,z\right) ,\left( x,z+y\right) \in R
\end{equation*}%
and hence $\forall ^{\ast }z\in G$,%
\begin{equation*}
f\left( x\right) +f\left( y\right) +f\left( z\right) =f\left( x\right)
+f\left( y+z\right) =f\left( x+y+z\right) =f\left( x+y\right) +f\left(
z\right) \text{.}
\end{equation*}

We form a discrete group $\dot{G}$ by specifying its generators and
relations. The generators are denoted by $\left( x\right) $ for $x\in U$ and
are subject to the relations $\left( x\right) +\left( y\right) =\left(
x+y\right) $ for $x,y,x+y$ in $U$. By definition, the function $\eta :\left(
x\right) \mapsto x$ extends to a homomorphism $\dot{G}\rightarrow G$: by 
\cite[Lemma 6]{weil_sur_1964}, this is in fact an isomorphism.

The assignment%
\begin{equation*}
\left( x\right) \mapsto f\left( x\right)
\end{equation*}%
for $x\in U$ respects the defining relations of $\dot{G}$. Therefore, it
extends to a \emph{homomorphism}%
\begin{equation*}
\dot{f}:\dot{G}\rightarrow D\text{.}
\end{equation*}%
Define now%
\begin{equation*}
\varphi :=\dot{f}\circ \eta ^{-1}:G\rightarrow D\text{.}
\end{equation*}%
Since $\eta ^{-1}$ and $\dot{f}$ are homomorphisms, so is $\varphi $.
Furthermore, for $x\in U$ we have that $\varphi \left( x\right) =f\left(
x\right) $. Since $U\subseteq G\ $is conull, this shows that $\forall ^{\ast
}x$, $\varphi \left( x\right) =f\left( x\right) $.

For uniqueness, it suffices to prove that if $f:G\rightarrow D$ is a group
homomorphism such that $\forall ^{\ast }x$, $f\left( x\right) =0$, then $f=0$%
. Suppose that $W:=\left\{ x\in G:f\left( x\right) =0\right\} $ is conull.
Fix $z\in G$. Then $W\cap \left( z-W\right) $ is also conull. This implies
that $z=x+y$ for some $x,y\in W$, and hence $f\left( z\right) =f\left(
x\right) +f\left( y\right) =0$.
\end{proof}

\begin{corollary}
\label{Corollary:Weil-discrete}Suppose that $G$ is a locally compact Polish
abelian group and $D$ is a (not necessarily countable) abelian group. If $%
c_{0},c_{1}:G^{2}\rightarrow D$ are \emph{normalized} $2$-cocycles such that 
$\forall ^{\ast }x$, $\forall ^{\ast }y$, $c_{0}\left( x,y\right)
=c_{1}\left( x,y\right) $, then there exists a function $u:G\rightarrow D$
such that:

\begin{enumerate}
\item $\delta u=c_{1}-c_{0}$, and

\item $\forall ^{\ast }x\in G$, $u\left( x\right) =0$.
\end{enumerate}
\end{corollary}

\begin{proof}
Set $c:=c_{0}-c_{1}$. By hypothesis, we have that $\forall ^{\ast }x$, $%
\forall ^{\ast }y$, $c\left( x,y\right) =0$. Consider the extension%
\begin{equation*}
0\rightarrow D\rightarrow E\overset{p}{\rightarrow }G\rightarrow 0
\end{equation*}%
determined by $c$. Thus, by definition $E$ is the abelian group that has $%
D\times G$ as set of elements, with group operation defined by%
\begin{equation*}
\left( a,x\right) +\left( b,y\right) :=\left( a+b+c\left( x,y\right)
,x+y\right) \text{.}
\end{equation*}%
Consider the function%
\begin{equation*}
s_{0}:G\rightarrow E\text{, }x\mapsto \left( 0,x\right) \text{.}
\end{equation*}%
Then $\forall ^{\ast }x$, $\forall ^{\ast }y$, 
\begin{equation*}
s_{0}\left( x+y\right) =s_{0}\left( x\right) +s_{0}\left( y\right) \text{.}
\end{equation*}%
Therefore, the existence claim in Lemma \ref{Lemma:Weil-discrete} produces a
homomorphism $s:G\rightarrow E$ such that $\forall ^{\ast }x$, $s\left(
x\right) =s_{0}\left( x\right) $. Consider now the homomorphisms $p\circ
s:G\rightarrow G$ and $\mathrm{id}_{G}:G\rightarrow G$. Then we have that $%
\forall ^{\ast }x$, $\left( p\circ s\right) \left( x\right) =x$. Then the
uniqueness claim in Lemma \ref{Lemma:Weil-discrete} shows that $p\circ s=%
\mathrm{id}_{G}$.

Let $v:G\rightarrow D$ be the function such that, for $x\in G$, $s\left(
x\right) =\left( v\left( x\right) ,x\right) $. Since $\forall ^{\ast }x$, $%
s\left( x\right) =s_{0}\left( x\right) $, we have $\forall ^{\ast }x$, $%
v\left( x\right) =0$. Furthermore, since $s$ is a homomorphism, we have $%
\forall x,y\in G$,%
\begin{equation*}
c\left( x,y\right) =v\left( x+y\right) -v\left( x\right) -v\left( y\right) 
\text{.}
\end{equation*}%
This shows that $c_{1}-c_{0}=-c=\delta v$.
\end{proof}

\begin{lemma}
\label{Lemma:Weil-topological}Suppose that $G$ is a locally compact Polish
abelian group, $H$ is a Polish abelian group, and $f:G\rightarrow H$ is
measurable and satisfies $f\left( x+y\right) =f\left( x\right) +f\left(
y\right) $ for almost every $x,y\in G$. Then there exists a unique group
homomorphism $\varphi :G\rightarrow H$ such that $f\left( x\right) =\varphi
\left( x\right) $ for almost every $x\in G$, and such a $\varphi $ is
necessarily continuous.
\end{lemma}

\begin{proof}
By Lemma \ref{Lemma:Weil-discrete}, there exists a unique group homomorphism 
$\varphi :G\rightarrow H$ such that $f\left( x\right) =\varphi \left(
x\right) $ for almost every $x\in G$. Such a $\varphi $ is continuous by the
analogue for measure of Pettis' Theorem \cite[Theorem 22.18]%
{hewitt_abstract_1979}.
\end{proof}

Suppose that $G$ is a locally compact abelian\ Polish group, and $A=\hat{A}%
/M $ is a group with a locally compact Polish cover. Define the continuous
group homomorphism%
\begin{equation*}
\delta :L^{0}(G,\hat{A})\rightarrow L^{0}(G\oplus G,\hat{A})\text{, }%
[f]\mapsto \lbrack \delta f]\text{.}
\end{equation*}%
Notice also that considering $M$ as a Polish subgroup of $\hat{A}$ allows
one to identify $L^{0}\left( G,M\right) $ with a Polish subgroup of $L^{0}(G,%
\hat{A})$, as we do in the following proposition.

\begin{lemma}
\label{Lemma:Borel-cocycle}Suppose that $G$ is a locally compact abelian
Polish group and $A$ is an abelian Polish group. Suppose that a measurable
function $c:G\oplus G\rightarrow A$ satisfies the cocycle and symmetry
identities almost everywhere. Then there exists a strict Borel symmetric
cocycle $c_{0}:G\oplus G\rightarrow A$ such that $c=c_{0}$ almost everywhere.
\end{lemma}

\begin{proof}
By \cite[Theorem 5]{moore_group_1976}, which establishes the isomorphism
between measurable group cohomology defined in terms of Borel cocycles and
measurable cocycles, applied in the case of degree $2$, we obtain a Borel
function $d:G\oplus G\rightarrow A$ satisfying $d\left( x+y,z\right)
+d\left( x,y\right) =d\left( x,y+z\right) +d\left( y,z\right) $ for every $%
x,y,z\in G$, and that represents the same (measurable) cohomology class as $%
c $. Thus, there exists a Borel function $f:G\rightarrow A$ such that $c$ is
equal almost everywhere to $c_{0}:=d+\delta f$. Notice that $c_{0}$ is still
a strict Borel cocycle. Since $\forall ^{\ast }x$, $\forall ^{\ast }y$, $%
c\left( x,y\right) =c\left( y,x\right) $, the same holds for $c_{0}$. Define
now as in \cite[Corollary 1.9]{ludeking_cocycles_1994} the \emph{bicharacter}%
\begin{equation*}
\omega \left( x,y\right) :=c_{0}\left( x,y\right) -c_{0}\left( y,x\right) 
\text{;}
\end{equation*}%
see also \cite{kleppner_multipliers_1965,kleppner_multipliers_1993}. Notice
that $\forall ^{\ast }x$, $\forall ^{\ast }y$, $\omega \left( x,y\right) =0$%
. Since $c_{0}$ is a strict cocycle and $G$ is abelian, a straightforward
computation shows that for every $z\in G$, the functions $\omega \left(
-,z\right) $ and $\omega \left( z,-\right) $ are group homomorphisms. By
continuity of Borel homomorphisms between Polish groups \cite[Theorem 9.10]%
{kechris_classical_1995}, $\omega $ is separately continuous. Since a
nonempty open set in $G$ has positive measure, this implies $\omega =0$,
concluding the proof.
\end{proof}

\begin{lemma}
\label{Lemma:strict-lift}Suppose that $G$ and $\hat A$ are locally compact
abelian Polish groups and that $M$ is a locally compact Polish subgroup of $%
\hat A$. If $f:G\rightarrow\hat A$ is Borel and 
\begin{equation*}
\delta f(x,y)\in M
\end{equation*}%
for almost every $(x,y)\in G^{2}$, then there is a Borel function $%
f_{0}:G\rightarrow\hat A$ such that $f_{0}=f$ almost everywhere and 
\begin{equation*}
\delta f_{0}(x,y)\in M
\end{equation*}%
for every $(x,y)\in G^{2}$.
\end{lemma}

\begin{proof}
If $G$ is discrete, ``almost everywhere'' means everywhere, and there is
nothing to prove. We may therefore assume that $G$ is nondiscrete.

Choose a Borel $M$-valued function $c$ agreeing almost everywhere with $%
\delta f$. Since $\delta f$ is a strict $\hat A$-valued cocycle, $c$
satisfies the cocycle and symmetry identities almost everywhere. Lemma \ref%
{Lemma:Borel-cocycle}, applied with coefficients in $M$, gives a strict
Borel symmetric $M$-valued cocycle $c_{0}$ such that 
\begin{equation*}
c_{0}=\delta f\quad\text{almost everywhere.}
\end{equation*}%
Changing $f$ at $0$, and changing $c_{0}$ by the coboundary of an $M$-valued
function supported at $0$, we may suppose that both $\delta f$ and $c_{0}$
are normalized. These changes do not alter their almost-everywhere classes.

Apply Corollary \ref{Corollary:Weil-discrete}, with coefficients in $\hat A$%
, to $c_{0}$ and $\delta f$. It gives a function $u:G\rightarrow\hat A$ such
that 
\begin{equation*}
\delta u=\delta f-c_{0} \quad\text{and}\quad u=0\quad\text{almost everywhere.%
}
\end{equation*}%
In the proof of that corollary, $u$ is obtained from an almost-everywhere
homomorphic section of the extension determined by the Borel cocycle $%
c_{0}-\delta f$. Since $G$ and $\hat A$ are locally compact Polish, Mackey's
theorem \cite[Theorem 2]{mackey_les_1957} equips that extension with a
locally compact Polish group topology having the product Borel structure.
Lemma \ref{Lemma:Weil-topological} then shows that the homomorphic section,
and hence $u$, may be taken Borel. Therefore $f_{0}:=f-u$ is Borel, $f_{0}=f$
almost everywhere, and $\delta f_{0}=c_{0}$ is $M$-valued everywhere.
\end{proof}

\begin{proposition}
\label{Proposition:Yoneda-Hom}Suppose that $G$ is a locally compact abelian\
Polish group and $A=\hat{A}/M$ is an abelian group with a locally compact
Polish cover. Let $\mathrm{Z}\left( \mathrm{Hom}_{\mathrm{Yon}}\left(
G,A\right) \right) $ be the group of pairs%
\begin{equation*}
\left( \lbrack f],[c]\right) \in L^{0}(G,\hat{A})\oplus L^{0}\left( G\oplus
G,M\right)
\end{equation*}%
such that%
\begin{equation*}
\delta \lbrack f]=[c]\text{.}
\end{equation*}%
Define also%
\begin{equation*}
\mathrm{B}\left( \mathrm{Hom}_{\mathrm{Yon}}\left( G,A\right) \right)
\subseteq \mathrm{Z}\left( \mathrm{Hom}_{\mathrm{Yon}}\left( G,A\right)
\right)
\end{equation*}%
to be the subgroup comprising the pairs $\left( [f],[c]\right) \in \mathrm{Z}%
\left( \mathrm{Hom}_{\mathrm{Yon}}\left( G,A\right) \right) $ such that $%
[f]\in L^{0}\left( G,M\right) $.

Then:

\begin{enumerate}
\item $\mathrm{Z}\left( \mathrm{Hom}_{\mathrm{Yon}}\left( G,A\right) \right) 
$ is a closed subgroup of $L^{0}(G,\hat{A})\times L^{0}\left( G^{2},M\right) 
$;

\item $\mathrm{B}\left( \mathrm{Hom}_{\mathrm{Yon}}\left( G,A\right) \right) 
$ is a Polish subgroup of $\mathrm{Z}\left( \mathrm{Hom}_{\mathrm{Yon}%
}\left( G,A\right) \right) $;

\item $\mathrm{Hom}_{\mathrm{Yon}}\left( G,A\right) :=\mathrm{Z}\left( 
\mathrm{Hom}_{\mathrm{Yon}}\left( G,A\right) \right) /\mathrm{B}\left( 
\mathrm{Hom}_{\mathrm{Yon}}\left( G,A\right) \right) $ is an abelian group
with a Polish cover;

\item $\mathrm{Hom}_{\mathrm{Yon}}\left( G,A\right) $ is naturally
isomorphic in $\mathbf{DAb}$ to $\mathrm{Hom}_{\mathrm{LH}\left( \mathbf{%
LCPAb}\right) }\left( G,A\right) $.
\end{enumerate}
\end{proposition}

\begin{proof}
Notice that $\mathrm{Z}\left( \mathrm{Hom}_{\mathrm{Yon}}\left( G,A\right)
\right) $ is the pullback in $\mathbf{PAb}$ of%
\begin{equation*}
\delta :L^{0}(G,\hat{A})\rightarrow L^{0}\left( G\oplus G,\hat{A}\right)
\end{equation*}%
and the inclusion%
\begin{equation*}
L^{0}\left( G\oplus G,M\right) \rightarrow L^{0}(G\oplus G,\hat{A})\text{.}
\end{equation*}%
Likewise, $\mathrm{B}\left( \mathrm{Hom}_{\mathrm{Yon}}\left( G,A\right)
\right) $ is the image of the continuous group homomorphism%
\begin{equation*}
L^{0}\left( G,M\right) \rightarrow L^{0}(G,\hat{A})\oplus L^{0}\left(
G\oplus G,M\right) \text{, }[g]\mapsto \left( \lbrack g],\delta \lbrack
g]\right) \text{.}
\end{equation*}

We now describe a natural Borel-definable isomorphism 
\begin{equation*}
\mathrm{Hom}_{\mathrm{LH}\left( \mathbf{LCPAb}\right) }\left( G,A\right)
\rightarrow \mathrm{Hom}_{\mathrm{Yon}}\left( G,A\right) .
\end{equation*}%
An element of $\mathrm{Hom}_{\mathrm{LH}\left( \mathbf{LCPAb}\right) }\left(
G,A\right) $ is by definition of the left heart of a quasi-abelian category 
\cite{schneiders_quasi-abelian_1999} represented by a triple $\left(
Z,g,\sigma \right) $ where $Z=\hat{Z}/N$ is an abelian group with a Polish
cover, $\sigma :Z\rightarrow G$ is a liftable isomorphism with a continuous
homomorphism $\hat{\sigma}:\hat{Z}\rightarrow G $ as a lift, and $%
g:Z\rightarrow A$ is a liftable homomorphism with a continuous homomorphism $%
\hat{g}:\hat{Z}\rightarrow \hat{A}$ as a lift. Let $t:G\rightarrow \hat{Z}$
be a Borel right inverse for $\hat{\sigma}:\hat{Z}\rightarrow G$. Then $f:=%
\hat{g}\circ t:G\rightarrow \hat{A}$ defines an element $\left( [f],[\delta
f]\right) $ of $\mathrm{Z}\left( \mathrm{Hom}_{\mathrm{Yon}}\left(
G,A\right) \right) $. Conversely, consider an element $\left( [f],[\delta
f]\right) $ of $\mathrm{Z}\left( \mathrm{Hom}_{\mathrm{Yon}}\left(
G,A\right) \right) $. Choose a Borel representative $f:G\rightarrow \hat{A}$%
. Initially, $\delta f$ is known to be $M$-valued only almost everywhere. By
Lemma \ref{Lemma:strict-lift}, after changing $f$ on a null set we may
suppose that 
\begin{equation*}
\delta f(x,y)\in M
\end{equation*}%
for every $x,y\in G$. Define%
\begin{equation*}
\hat{Z}:=\left\{ \left( x,a\right) \in G\oplus \hat{A}:f\left( x\right)
\equiv a\ \mathrm{mod}M\right\}
\end{equation*}%
and%
\begin{equation*}
N=\left\{ 0\right\} \oplus M\subseteq \hat{Z}\text{.}
\end{equation*}%
The pointwise condition on $\delta f$ makes $\hat{Z}$ a group. Equivalently,
using $f$ it is the Borel extension of $G$ by $M$ associated with the strict
Borel cocycle $\delta f$. Mackey's theorem \cite[Theorem 2]{mackey_les_1957}
supplies the corresponding locally compact Polish group topology. We obtain
an element of $\mathrm{Hom}_{\mathrm{LH}\left( \mathbf{LCPAb}\right) }\left(
G,A\right) $ by letting $Z:=\hat{Z}/N$ be the graph of the induced
homomorphism and taking $\sigma $ and $g$ to be the first and second
coordinate projections, respectively. These constructions describe an
isomorphism 
\begin{equation*}
\Psi :\mathrm{Hom}_{\mathrm{LH}\left( \mathbf{LCPAb}\right) }\left(
G,A\right) \rightarrow \mathrm{Hom}_{\mathrm{Yon}}\left( G,A\right)
\end{equation*}%
together with its inverse. In order to argue that $\Psi $ is
Borel-definable, it suffices to show that it is $\boldsymbol{\Sigma }%
_{1}^{1} $-definable {by \cite[Proposition 4.7]{lupini_looking_2024}}. For
this, observe that the lift of the graph of $\Psi $ is the set of pairs $%
\left( \left( Z,g,\sigma \right) ,f\right) $ such that, adopting the
notations above, there exists $t\in \mathrm{Z}\left( \mathrm{Hom}_{\mathrm{%
Yon}}\left( G,Z\right) \right) $ that represents an isomorphism in $\mathrm{%
Hom}_{\mathrm{Yon}}\left( G,Z\right) $ and such that $f=\hat{g}\circ t$.
\end{proof}

\begin{proposition}
\label{Proposition:Yoneda-Ext}Suppose that $G$ is a locally compact abelian
Polish group and $A$ is a locally compact abelian Polish group. Define 
\begin{equation*}
\mathrm{Z}\left( \mathrm{Ext}_{\mathrm{Yon}}^{1}\left( G,A\right) \right)
\subseteq L^{0}\left( G\oplus G,A\right)
\end{equation*}%
to be the group of $[c]\in L^{0}\left( G\oplus G,A\right) $ for which the
cocycle and symmetry identities hold almost everywhere. Define also $\mathrm{%
B}\left( \mathrm{Ext}_{\mathrm{Yon}}^{1}\left( G,A\right) \right) $ to be
the image of the continuous homomorphism%
\begin{equation*}
L^{0}\left( G,A\right) \rightarrow L^{0}\left( G\oplus G,A\right) \text{, }%
[f]\mapsto \delta \lbrack f]\text{.}
\end{equation*}%
Then:

\begin{enumerate}
\item $\mathrm{Z}\left( \mathrm{Ext}_{\mathrm{Yon}}^{1}\left( G,A\right)
\right) $ is a closed subgroup of $L^{0}\left( G\oplus G,A\right) $;

\item $\mathrm{B}\left( \mathrm{Ext}_{\mathrm{Yon}}^{1}\left( G,A\right)
\right) $ is a Polish subgroup of $\mathrm{Z}\left( \mathrm{Ext}_{\mathrm{Yon%
}}^{1}\left( G,A\right) \right) $;

\item The quotient 
\begin{equation*}
\mathrm{Ext}_{\mathrm{Yon}}^{1}\left( G,A\right) :=\mathrm{Z}\left( \mathrm{%
Ext}_{\mathrm{Yon}}^{1}\left( G,A\right) \right) /\mathrm{B}\left( \mathrm{%
Ext}_{\mathrm{Yon}}^{1}\left( G,A\right) \right)
\end{equation*}%
is a group with a Polish cover naturally isomorphic in $\mathbf{DAb}$ to 
\begin{equation*}
\mathrm{Ext}^{1}\left( G,A\right) \cong \mathrm{Hom}_{\mathrm{D}^{b}\left( 
\mathbf{LCPAb}\right) }\left( G,A[1]\right) .
\end{equation*}
\end{enumerate}
\end{proposition}

\begin{proof}
Notice that $\mathrm{Z}\left( \mathrm{Ext}_{\mathrm{Yon}}^{1}\left(
G,A\right) \right) $ is a closed subgroup of $L^{0}\left( G\oplus G,A\right) 
$, being the kernel of the continuous homomorphism%
\begin{equation*}
L^{0}\left( G\oplus G,A\right) \rightarrow L^{0}\left( G\oplus G\oplus
G,A\right) \oplus L^{0}\left( G\oplus G,A\right)
\end{equation*}%
\begin{equation*}
\lbrack c]\mapsto \left( \lbrack \delta c],[\sigma c]\right)
\end{equation*}%
where%
\begin{equation*}
\delta c\left( x,y,z\right) =c\left( x,y\right) +c\left( x+y,z\right)
-c\left( x,y+z\right) -c\left( y,z\right)
\end{equation*}%
and%
\begin{equation*}
\sigma c\left( x,y\right) =c\left( x,y\right) -c\left( y,x\right)
\end{equation*}%
for $x,y,z\in G$. Likewise, we have that $\mathrm{B}\left( \mathrm{Ext}_{%
\mathrm{Yon}}^{1}\left( G,A\right) \right) $ is the image of a continuous
group homomorphism%
\begin{equation*}
L^{0}\left( G,A\right) \rightarrow L^{0}\left( G\oplus G,A\right)
\end{equation*}%
and hence a Polish subgroup of $L^{0}\left( G\oplus G,A\right) $.

By Lemma \ref{Lemma:Borel-cocycle}, every element of $\mathrm{Z}\left( 
\mathrm{Ext}_{\mathrm{Yon}}^{1}\left(G,A\right)\right)$ has a strict Borel
symmetric cocycle as a representative of the same almost-everywhere class.

We now show that 
\begin{equation*}
\mathrm{Ext}_{\mathrm{Yon}}^{1}\left( G,A\right)
\end{equation*}%
is naturally isomorphic in $\mathbf{DAb}$ to 
\begin{equation*}
\mathrm{Ext}^{1}\left( G,A\right) \cong \mathrm{Hom}_{\mathrm{D}^{b}\left( 
\mathbf{LCPAb}\right) }\left( G,A[1]\right) .
\end{equation*}%
By the usual Yoneda description of $\mathrm{Ext}^{1}$ in an exact category,
equivalently by \cite[Theorem III.5.5]{gelfand_methods_2003}, an element of $%
\mathrm{Ext}^{1}\left( G,A\right) $ is represented by a strict short exact
sequence in $\mathbf{LCPAb}$%
\begin{equation*}
0\rightarrow A\overset{\alpha }{\rightarrow }R\overset{\beta }{\rightarrow }%
G\rightarrow 0
\end{equation*}%
where $R$ is a locally compact Polish abelian group. Choose a Borel section $%
t:G\rightarrow R$ of $\beta $. There is then a unique Borel function $%
c:G^{2}\rightarrow A$ satisfying 
\begin{equation*}
\alpha \left( c(x,y)\right) =t(x)+t(y)-t(x+y).
\end{equation*}%
Associativity and commutativity in $R$ show pointwise that $c$ is a strict
Borel symmetric cocycle. Replacing $t$ changes $c$ by a Borel coboundary,
and equivalent extensions give the same class in $\mathrm{Ext}_{\mathrm{Yon}%
}^{1}\left( G,A\right) $.

Conversely, given an element $[c]$ of $\mathrm{Z}\left( \mathrm{Ext}_{%
\mathrm{Yon}}^{1}\left( G,A\right) \right) $, Lemma \ref{Lemma:Borel-cocycle}
allows us to choose a strict Borel symmetric cocycle in the same
almost-everywhere class. By a standard argument, $c\left( x,0\right) $ is
independent of the choice of $x\in G$. Thus, after replacing $c$ with $c_{%
\mathrm{norm}}$ defined by%
\begin{equation*}
c_{\mathrm{norm}}\left( x,y\right) :=c\left( x,y\right) -c\left( 0,0\right)
\end{equation*}
we may replace $c$ by a cohomologous strict normalized cocycle. Define the
corresponding extension%
\begin{equation*}
0\rightarrow A\overset{\alpha }{\rightarrow }R\overset{\beta }{\rightarrow }%
G\rightarrow 0
\end{equation*}%
where $R=A\times G$ is the standard Borel group endowed with the group
operation 
\begin{equation*}
\left( a,g\right) +\left( a^{\prime },g^{\prime }\right) =\left( a+a^{\prime
}+c\left( g,g^{\prime }\right) ,g+g^{\prime }\right) \text{,}
\end{equation*}%
the product Borel structure, and the product of Haar measures on $A$ and $G$%
. By \cite[Theorem 1]{mackey_les_1957} there exists a unique locally compact
Polish group topology on $R$ whose Borel sets are the given Borel sets and
for which the displayed measure is a Haar measure. This topology renders the
extension%
\begin{equation*}
0\rightarrow A\overset{\alpha }{\rightarrow }R\overset{\beta }{\rightarrow }%
G\rightarrow 0
\end{equation*}%
an extension of locally compact abelian Polish groups. In turn, this defines
an element of $\mathrm{Ext}^{1}\left( G,A\right) $ as in \cite[Theorem
III.5.5]{gelfand_methods_2003}.

These constructions describe an isomorphism $\mathrm{Ext}_{\mathrm{Yon}%
}^{1}\left( G,A\right) \rightarrow \mathrm{Ext}^{1}\left( G,A\right) $
together with its inverse. It remains to prove that this isomorphism is
Borel-definable. By \cite[Proposition 4.7]{lupini_looking_2024}, it suffices
to show that it is $\boldsymbol{\Sigma }_{1}^{1}$-definable. We have that
the lift of its graph is given by pairs $\left( c,\left( f,\sigma \right)
\right) $ where $\left( f,\sigma \right) $ represents the element of $%
\mathrm{Ext}^{1}\left( G,A\right) $ associated with the extension%
\begin{equation*}
0\rightarrow A\overset{\alpha }{\rightarrow }R\overset{\beta }{\rightarrow }%
G\rightarrow 0
\end{equation*}%
as in \cite[Theorem III.5.5]{gelfand_methods_2003} such that there exists $%
t\in L^{0}\left( G,R\right) $ satisfying $\beta (t(x))=x$ for almost every $%
x\in G$ and, for almost every $(x,y)\in G^{2}$, 
\begin{equation*}
\alpha \left( c\left( x,y\right) \right) =t\left( y\right) -t\left(
x+y\right) +t\left( x\right) \text{.}
\end{equation*}%
By \cite[Theorem 17.25]{kechris_classical_1995}, this shows that the lift of
the graph is analytic, concluding the proof.
\end{proof}

\subsection{Morphisms}

We isolate a characterization of morphisms in the left heart of abelian
locally compact Polish groups, building on \cite%
{lupini_looking_2024,bergfalk_definable_2024,casarosa_homological_2026}.
Notice that if $G/N$ is a group with a locally compact Polish cover, then $%
N=V\oplus M$ where $V$ is a finite-dimensional vector group and $M$ has a
compact open subgroup. Thus, $G/N$ is isomorphic to $H/\left( V\oplus
D\right) $ where $H=C/U$ and $D=M/U$ is countable. We call this a group with
a locally compact Polish cover in \emph{standard form}.\ Following \cite%
{casarosa_homological_2026}, we let $\mathbf{ProLiePAb}$ be the category of
pro-Lie Polish abelian groups. It is shown in \cite%
{casarosa_homological_2026} that $\mathbf{ProLiePAb}$ is a quasi-abelian
category with enough projectives and homological dimension $1$, containing $%
\mathbf{LCPAb}$ as a thick subcategory.

We recall the following lemma from \cite[Theorem 2]{karube_local_1958}; see
also \cite%
{gleason_spaces_1950,mostert_local_1953,mostert_sections_1956,kehlet_cross_1984}%
.

\begin{lemma}
\label{Lemma:local-section}Suppose that $\pi :Z\rightarrow G$ is a
surjective continuous homomorphism between locally compact Polish groups. If 
$\mathrm{Ker}\left( \pi \right) $ is a Lie group, then there exist an
identity neighborhood $U$ of $G$ and a Borel right inverse $\varphi
:U\rightarrow Z$ for $\pi $ such that $\varphi |_{U}$ is continuous.
\end{lemma}

\begin{lemma}
\label{Lemma:factor-out}Let $G/N$ and $H/M$ be groups with a pro-Lie Polish
cover. Suppose that $M=M_{0}\oplus W$ where $W$ is a vector group. Let $%
f:G/N\rightarrow H/M$ be a Borel-definable group homomorphism. Then $f$
factors in $\mathrm{LH}\left( \mathbf{ProLiePAb}\right) $ through the
quotient map $H/M_{0}\rightarrow H/M$.
\end{lemma}

\begin{proof}
Let $\varphi :G\rightarrow H$ be a Borel lift for $f$. Consider the
corresponding $M$-valued Borel cocycle $c$ on $G$, defined by%
\begin{equation*}
c(x,y):=\varphi (x)+\varphi (y)-\varphi \left( x+y\right) \text{.}
\end{equation*}%
Then we can write $c=c_{0}+c^{\prime }$ where $c^{\prime }$ is a $W$-valued
Borel cocycle on $G$.\ Thus, $c^{\prime }$ defines an element of $\mathrm{Ext%
}\left( G,W\right) $ \cite[Proposition 4.31]{casarosa_homological_2026}, the
derived functor of $\mathrm{Hom}$ in the category of pro-Lie Polish abelian
groups. Then by injectivity of vector groups in $\mathbf{ProLiePAb}$, $%
\mathrm{Ext}\left( G,W\right) =0$. Thus, there exists a Borel function $%
s:G\rightarrow W$ such that $c^{\prime }(x,y)=s(x)+s(y)-s\left( x+y\right) $%
. After replacing $\varphi $ with $\varphi -s$, we can assume that $%
c^{\prime }=0$. Thus, $\varphi $ defines a Borel-definable homomorphism $%
G/N\rightarrow H/M_{0}$ that lifts $f$.
\end{proof}

\begin{proposition}
\label{Proposition:morphisms}Let $G/N$ and $H/M$ be groups with a locally
compact Polish cover, where $M$ is an abelian Lie group with no nontrivial
compact connected subgroups. Let $f:G/N\rightarrow H/M$ be a Borel-definable
homomorphism. Then there exist a zero neighborhood $U$ of $G$ and a Borel
function $\varphi :G\rightarrow H$ such that:

\begin{enumerate}
\item $f\left( x+N\right) =\varphi (x)+M$ for $x\in G$,

\item $\varphi |_{U}$ is continuous,

\item $\varphi \left( x+y\right) =\varphi (x)+\varphi (y)$ for $x,y\in U$.
\end{enumerate}
\end{proposition}

\begin{proof}
By Lemma \ref{Lemma:factor-out} and the structure theorem for abelian Lie
groups, we can assume that $M=M_{0}$ is \emph{countable}. Define 
\begin{equation*}
Z:=\left\{ (x,y)\in G\oplus H:f\left( x+N\right) =y+M\right\} \text{.}
\end{equation*}%
As in the proof of \cite[Proposition 6.16]{lupini_looking_2024}, we have a
locally compact Polish group topology on $Z$ for which the sequence%
\begin{equation*}
0\rightarrow \left\{ 0\right\} \oplus M\rightarrow Z\overset{\pi }{%
\rightarrow }G\rightarrow 0\text{.}
\end{equation*}%
is topologically exact. Since $M$ is a countable discrete Lie group, Lemma %
\ref{Lemma:local-section} yields a zero neighborhood $W$ of $G$ and a
continuous function $\sigma :W\rightarrow Z$ such that $\pi \sigma (x)=x$
for every $x\in W$. Shrinking $W$ if necessary, choose a zero neighborhood $%
U_{0}\subseteq W$. Let%
\begin{equation*}
\psi :Z\rightarrow H\text{, }(x,y)\mapsto y\text{.}
\end{equation*}%
Then $\left. \left( \psi \circ \sigma \right) \right\vert
_{U_{0}}:U_{0}\rightarrow H$ is a continuous function such that%
\begin{equation*}
\left( \psi \circ \sigma \right) (x)+M=f\left( x+N\right)
\end{equation*}%
for every $x\in U_{0}$. By \cite[Lemma 3.12]{lupini_looking_2024}, this
shows that $f$ is \emph{locally continuously definable }as in \cite[%
Definition 3.8]{lupini_looking_2024}. By \cite[Proposition 5.3]%
{lupini_looking_2024}, since $M$ is a countable discrete group, there exists
a Borel function $\varphi :G\rightarrow H$ and a compact zero neighborhood $%
U\subseteq U_{0}$ satisfying the desired conclusions.
\end{proof}

\begin{proposition}
\label{Proposition:morphisms-finite-dimensional}Let $G/N$ and $H/M$ be
groups with a locally compact Polish cover, where $M$ is an abelian Lie
group. Let $f:G/N\rightarrow H/M$ be a Borel-definable homomorphism. Then
there exist a zero neighborhood $U$ of $G$ and a Borel function $\varphi
:G\rightarrow H$ such that:

\begin{enumerate}
\item $f\left( x+N\right) =\varphi (x)+M$ for $x\in G$,

\item $\varphi |_{U}$ is continuous.
\end{enumerate}
\end{proposition}

\begin{proof}
Let $C\subseteq M$ be the maximal compact connected subgroup. Then by
Proposition \ref{Proposition:morphisms} applied to $\left( H/C\right)
/\left( M/C\right) $, we obtain a Borel function $\varphi :G\rightarrow H/C$
and a zero neighborhood $U$ of $G$, such that $\varphi \left( x\right)
+M/C=f\left( x\right) $ and $\varphi |_{U}:U\rightarrow H/C$ is continuous,
and $\varphi \left( x+y\right) =\varphi \left( x\right) +\varphi \left(
y\right) $ for $x,y\in U$. Since $C$ is a torus, and in particular a Lie
group, by Lemma \ref{Lemma:local-section} there is a Borel right inverse $%
\psi :H/C\rightarrow H$ for the quotient map $H\rightarrow H/C$ and a zero
neighborhood $V$ of $H/C$ such that $\psi |_{V}$ is continuous. After
replacing $U$ with a smaller zero neighborhood of $G$ we can assume that $%
\varphi \left( U\right) \subseteq V$. Thus, $\psi \circ \varphi $ is the
required function.

Alternatively, one can repeat directly the argument of Proposition \ref%
{Proposition:morphisms}.
\end{proof}

Let us say that a tower $\left( G_{n}\right) $ of Polish groups is \emph{%
epimorphic }if all the bonding maps $G_{n+1}\rightarrow G_{n}$ for $n\in
\omega $ are surjective. Recall that if $G$ is a locally compact Polish
abelian group, then there exists an epimorphic tower $\left( G_{n}\right) $
of abelian Lie groups such that $G$ is isomorphic to the limit $\mathrm{lim}%
_{n}\ G_{n}$.

\begin{lemma}
\label{Lemma:factor}Let $G$ be a locally compact Polish abelian group. Let
also $\left( G_{n}\right) $ be an epimorphic tower of locally compact Polish
groups with $\mathrm{lim}_{n}\ G_{n}\cong G$. Let also $H/M$ be a group with
a locally compact Polish cover, where $M$ is an abelian Lie group with no
nontrivial compact connected subgroups. Suppose that $f:G\rightarrow H/M$ is
a Borel-definable homomorphism. Assume also that, for every $n\in \mathbb{N}$
and every continuous homomorphism%
\begin{equation*}
u:\mathrm{Ker}\left( \pi _{n}:G\rightarrow G_{n}\right) \rightarrow H,
\end{equation*}%
there exists $k\geq n$ such that $u$ vanishes on $\mathrm{Ker}\left( \pi
_{k}\right) $. Then $f$ factors through $\pi _{\ell }:G\rightarrow G_{\ell }$
for some $\ell \in \mathbb{N}$.
\end{lemma}

\begin{proof}
Apply Proposition \ref{Proposition:morphisms} (in the particular case when $%
N=0$) to the Borel-definable homomorphism $f:G\rightarrow H/M$. We obtain a
compact zero neighborhood $U$ of $G$ and a Borel function $\varphi
:G\rightarrow H$ such that:

\begin{enumerate}
\item $f(x)=\varphi (x)+M$ for every $x\in G$;

\item $\varphi |_{U}$ is continuous;

\item $\varphi \left( x+y\right) =\varphi (x)+\varphi (y)$ for every $x,y\in
U$.
\end{enumerate}

Then there exists $n\in \mathbb{N}$ such that $\mathrm{Ker}\left( \pi
_{n}\right) \subseteq U$. Local additivity gives $\varphi \left( 0\right) =0$%
. Since $\mathrm{Ker}\left( \pi _{n}\right) $ is a subgroup contained in $U$%
, it also shows that%
\begin{equation*}
\varphi |_{\mathrm{Ker}\left( \pi _{n}\right) }:\mathrm{Ker}\left( \pi
_{n}\right) \rightarrow H
\end{equation*}%
is a continuous group homomorphism. The hypothesis gives $k\geq n$ such that%
\begin{equation*}
\varphi |_{\mathrm{Ker}\left( \pi _{k}\right) }=0.
\end{equation*}

Suppose now that $x,y\in G$ satisfy $\pi _{k}(x)=\pi _{k}(y)$. Then $x-y\in 
\mathrm{Ker}\left( \pi _{k}\right) \subseteq U$. Since $f$ is a
homomorphism, we have%
\begin{equation*}
f(x)-f(y)=f\left( x-y\right) =\varphi \left( x-y\right) +M=M\text{.}
\end{equation*}%
Thus, $f$ factors through a homomorphism $\bar{f}:G_{k}\rightarrow H/M$.
Since $\pi _{k}$ is a continuous surjective homomorphism between Polish
groups, it admits a Borel right inverse $s:G_{k}\rightarrow G$. The Borel
function $\varphi \circ s:G_{k}\rightarrow H$ is a lift of $\bar{f}$, so $%
\bar{f}$ is Borel-definable. Therefore $f=\bar{f}\circ \pi _{k}$.
\end{proof}

\begin{lemma}
\label{Lemma:Hom(-,V/H)-0}Let $T=\mathbb{T}$, $V=\mathbb{R}$, and let $%
H\subseteq V$ be a countable subgroup. Then $\mathrm{Hom}\left( T,V/H\right) 
$ is countable.
\end{lemma}

\begin{proof}
Let $G$ be the closure of $H$ in $V$. If $G$ is a proper subgroup of $V$,
then $G$ is discrete and hence $H=G$. If $G$ is trivial, then $\mathrm{Hom}%
\left( T,V/H\right) =0$. Otherwise, $H=G$ is a lattice in $V$, so $V/H\cong 
\mathbb{T}$ and $\mathrm{Hom}\left( T,V/H\right) \cong \mathbb{Z}$. Thus, we
can assume that $G=V$ and $H$ is dense in $V$.

Let $f:T\rightarrow V/H$ be a Borel-definable homomorphism. By Proposition %
\ref{Proposition:morphisms-finite-dimensional}, we obtain a Borel lift $%
\varphi :T\rightarrow V$ of $f$ and a zero neighborhood $U$ of $T$ such that 
$\varphi |_{U}$ is continuous and additive. Let $q:\mathbb{R}\rightarrow 
\mathbb{T}=\mathbb{R}/\mathbb{Z}$ be the quotient map. Continuous local
additivity implies that there exist $\varepsilon >0$ and $a\in \mathbb{R}$
such that%
\begin{equation*}
\varphi \left( q(t)\right) =at
\end{equation*}%
whenever $\left\vert t\right\vert <\varepsilon $. For arbitrary $t\in 
\mathbb{R}$, choose a positive integer $n$ such that $\left\vert
t/n\right\vert <\varepsilon $. Since $f$ is a homomorphism, we have%
\begin{equation*}
f\left( q(t)\right) =nf\left( q(t/n)\right) =at+H.
\end{equation*}%
Taking $t=1$ shows that $a\in H$. Thus, $f$ is determined by the element $a$
of the countable group $H$. Therefore $\mathrm{Hom}\left( T,V/H\right) $ is
countable.
\end{proof}

\begin{proposition}
Let $V$ be a finite-dimensional vector group, $T$ a torus group, and let $%
H\subseteq V$ be a locally compact Polish subgroup. Then $\mathrm{Hom}\left(
T,V/H\right) $ is countable.
\end{proposition}

\begin{proof}
We can write $H=W\oplus H_{0}$ where $W$ is a finite-dimensional vector
group and $H_{0}$ has no nontrivial closed vector subgroups. Then $W$ is a
closed subgroup of $V$. After replacing $V$ with $V/W$ and $H$ with $%
H/W\cong H_{0}$, we can assume that $H$ has no nontrivial closed vector
subgroups. Then $H$ has a compact open subgroup $U$. Since $V$ has no
nontrivial compact subgroups, we must have that $U$ is trivial, and $H$ is
countable.

By Lemma \ref{Lemma:factor}, we can assume that $T$ is finite-dimensional.
Thus $T\cong \mathbb{T}^{k}$ for some $k\in \mathbb{N}$, and hence%
\begin{equation*}
\mathrm{Hom}\left( T,V/H\right) \cong \mathrm{Hom}\left( \mathbb{T}%
,V/H\right) ^{k}\text{.}
\end{equation*}%
Thus, we can assume without loss of generality that $T$ is $1$-dimensional.
We prove that the conclusion holds by induction on the dimension $d$ of $V$.
When $d=1$, this follows from Lemma \ref{Lemma:Hom(-,V/H)-0}. Suppose that $%
d>1$. Let $V=W\oplus V_{0}$ where $V_{0}$ is a $1$-dimensional closed
subgroup. Consider the short exact sequence%
\begin{equation*}
0\rightarrow V_{0}/\left( V_{0}\cap H\right) \rightarrow V/H\rightarrow
W/\pi _{W}\left( H\right) \rightarrow 0\text{.}
\end{equation*}%
This induces an exact sequence%
\begin{equation*}
0\rightarrow \mathrm{Hom}\left( T,V_{0}/\left( V_{0}\cap H\right) \right)
\rightarrow \mathrm{Hom}\left( T,V/H\right) \rightarrow \mathrm{Hom}(T,W/\pi
_{W}\left( H\right) )
\end{equation*}%
Then 
\begin{equation*}
\mathrm{Hom}\left( T,V_{0}/\left( V_{0}\cap H\right) \right)
\end{equation*}%
is countable by\ Lemma \ref{Lemma:Hom(-,V/H)-0}, and 
\begin{equation*}
\mathrm{Hom}(T,W/\pi _{W}\left( H\right) )
\end{equation*}%
is countable by the inductive hypothesis. The conclusion follows.
\end{proof}

\section{Extensions of Polish modules as groups with a Polish cover}

\label{Section:module-cocycles}

We now give a module version of the cocycle description from Section~\ref%
{Section:derived}. Throughout this section, $R$ denotes a countable discrete
unital ring, and all modules are left $R$-modules. The intended arithmetic
example is $R=\mathcal{O}$, the ring of integers of a number field. We write 
$\mathbf{LCPMod}_{R}$ for the category of locally compact Polish $R$-modules
and $\mathbf{proLiePMod}_{R}$ for the category of Polish $R$-modules whose
underlying additive groups are pro-Lie. Morphisms are continuous $R$-linear
maps. The strict exact structure is inherited from Polish abelian groups:
the first map in a strict short exact sequence is a closed embedding and the
second is an open surjection.

For $C\in\mathbf{LCPMod}_{R}$ and $A\in\mathbf{proLiePMod}_{R}$, let 
\begin{equation*}
\mathrm{Ext}_{R}^{1}(C,A)
\end{equation*}%
denote the Yoneda--Baer group of strict extensions in $\mathbf{proLiePMod}%
_{R}$. This definition is legitimate because locally compact abelian groups
are pro-Lie and pro-Lie Polish abelian groups form a thick subcategory of
the Polish abelian groups by \cite[Theorems 3.17 and 3.28]%
{casarosa_homological_2026}. Thus, the middle term of an extension of $C$ by 
$A$ is again pro-Lie. Whenever a derived $\mathrm{Ext}^{1}$ has been
constructed, the usual exact-category comparison identifies it with this
Yoneda group.

\subsection{Term-saturated spaces of measurable cochains}

There is a small measure-theoretic point that is absent for $R=\mathbb{Z}$.
If $r\in R$, the continuous endomorphism $x\mapsto rx$ of $C$ need not
preserve the Haar measure class. Consequently, $[f]\mapsto [f(r\,\cdot)]$
need not be defined on the ordinary space of Haar-a.e.\ classes. We use the
following canonical enlargement of the Haar-null ideal.

For $m,n\geq 1$, let $\mathscr{T}_{m,n}(R)$ be the countable collection of $%
R $-term maps $T:C^{m}\rightarrow C^{n}$. Thus, each coordinate of $T$ has
the form 
\begin{equation*}
(x_{1},\ldots,x_{m})\longmapsto \sum_{j=1}^{m}r_{j}x_{j},\qquad r_{j}\in R.
\end{equation*}%
Choose, for every $m\geq1$, a Borel probability measure $\lambda _{m}$ on $%
C^{m}$ equivalent to Haar measure. For every $n\geq1$, choose positive
weights with sum one and let $\mu _{n}$ be a convex combination of all the
probability measures 
\begin{equation}  \label{Equation:module-control-measures}
T_{\ast}\lambda _{m},\qquad m\geq1,\quad T\in\mathscr{T}_{m,n}(R).
\end{equation}
The measure class of $\mu _{n}$ is independent of the choices. Indeed,
equivalent choices of $\lambda _m$ have equivalent pushforwards, and every
measure in \eqref{Equation:module-control-measures} occurs with positive
weight. Notice that the identity term shows that $\mu _n$ dominates the Haar
measure class, while the zero term gives positive mass to the origin.

\begin{lemma}
\label{Lemma:module-control-measures} If $S\in\mathscr{T}_{n,k}(R)$, then 
\begin{equation*}
S_{\ast}\mu _{n}\ll\mu _{k}.
\end{equation*}%
Consequently, pullback along $S$ defines a continuous homomorphism 
\begin{equation*}
S^{\ast}:L^{0}(C^{k},A;\mu _{k})\longrightarrow L^{0}(C^{n},A;\mu _{n}).
\end{equation*}
\end{lemma}

\begin{proof}
Every component of $\mu _n$ has the form $T_{\ast}\lambda _m$. Its
pushforward under $S$ is $(S\circ T)_{\ast}\lambda _m$, which occurs among
the components defining $\mu _k$. This proves absolute continuity.
Convergence in $\mu _k$-measure implies convergence in $S_{\ast}\mu _n$%
-measure by absolute continuity, proving continuity of pullback.
\end{proof}

Set 
\begin{equation*}
\mathscr{L}_{n}^{0}(C,A):=L^{0}(C^{n},A;\mu _{n}).
\end{equation*}%
These are Polish abelian groups. All substitutions obtained from addition,
the zero element, diagonals, permutations of variables, and scalar
multiplication are continuous between the appropriate $\mathscr{L}%
_{n}^{0}(C,A)$. In particular, the value at zero is well defined, and 
\begin{equation*}
\mathscr{L}_{1,0}^{0}(C,A):= \{[f]\in\mathscr{L}_{1}^{0}(C,A):f(0)=0\}
\end{equation*}%
is a closed Polish subgroup.

\subsection{Module cocycles and coboundaries}

Suppose that 
\begin{equation*}
0\longrightarrow A\overset{\alpha}{\longrightarrow}E \overset{\beta}{%
\longrightarrow}C\longrightarrow0
\end{equation*}%
is a strict extension of Polish $R$-modules. Choose a Borel section $%
t:C\rightarrow E$ of $\beta$ with $t(0)=0$. Its additive defect and its
scalar defects are the Borel functions 
\begin{align}
\alpha(c(x,y))&=t(x)+t(y)-t(x+y),  \label{Equation:additive-module-cocycle}
\\
\alpha(\rho _r(x))&=r t(x)-t(rx),\qquad r\in R.
\label{Equation:scalar-module-cocycle}
\end{align}
The pair $(c,\rho)$, where $\rho=(\rho _r)_{r\in R}$, records the entire $R$%
-module extension.

Motivated by these formulas, define the group of module cocycles 
\begin{equation*}
\mathsf{Z}_{R}(C,A)\subseteq \mathscr{L}_{2}^{0}(C,A)\oplus \prod_{r\in R}%
\mathscr{L}_{1}^{0}(C,A)
\end{equation*}%
to consist of the pairs $(c,\rho)$ satisfying the following identities in
the corresponding term-saturated $L^{0}$-spaces: 
\begin{align}
c(x,0)&=0,  \label{Equation:module-cocycle-normalized} \\
c(x,y)&=c(y,x),  \label{Equation:module-cocycle-symmetric} \\
c(x,y)+c(x+y,z)&=c(x,y+z)+c(y,z),
\label{Equation:module-cocycle-associative} \\
\rho _r(x)+\rho _r(y)-\rho _r(x+y) &=r c(x,y)-c(rx,ry),
\label{Equation:module-cocycle-additive-action} \\
\rho _{r+s}(x)&=\rho _r(x)+\rho _s(x)+c(rx,sx),
\label{Equation:module-cocycle-ring-addition} \\
\rho _{rs}(x)&=r\rho _s(x)+\rho _r(sx),
\label{Equation:module-cocycle-ring-multiplication} \\
\rho _0(x)&=0,\qquad \rho _1(x)=0.  \label{Equation:module-cocycle-zero-unit}
\end{align}
Here $r,s$ range over $R$. The first three equations say that $c$ is a
normalized symmetric additive $2$-cocycle. Equation %
\eqref{Equation:module-cocycle-additive-action} says that the proposed
action of $r$ is additive on the extension, while %
\eqref{Equation:module-cocycle-ring-addition}--%
\eqref{Equation:module-cocycle-zero-unit} express the ring and unit axioms.

For $f\in\mathscr{L}_{1,0}^{0}(C,A)$, define 
\begin{align*}
(\delta f)(x,y)&:=f(x)+f(y)-f(x+y), \\
(\partial _r f)(x)&:=r f(x)-f(rx).
\end{align*}%
Lemma~\ref{Lemma:module-control-measures} makes these formulas well defined
and continuous on almost-everywhere classes. A direct calculation shows that 
\begin{equation}  \label{Equation:module-coboundary-map}
D:\mathscr{L}_{1,0}^{0}(C,A)\longrightarrow\mathsf{Z}_{R}(C,A), \qquad
f\longmapsto(\delta f,(\partial _r f)_{r\in R})
\end{equation}
is a continuous homomorphism. Define the group of module coboundaries by 
\begin{equation*}
\mathsf{B}_{R}(C,A):=\mathrm{im}(D).
\end{equation*}

\begin{lemma}[Module strictification]
\label{Lemma:module-cocycle-strictification}Every element of $\mathsf{Z}%
_{R}(C,A)$ has Borel representatives $c$ and $(\rho _r)_{r\in R}$ for which %
\eqref{Equation:module-cocycle-normalized}--%
\eqref{Equation:module-cocycle-zero-unit} hold pointwise. Two such strict
representatives of the same term-saturated almost-everywhere class determine
equivalent extensions.
\end{lemma}

\begin{proof}
This is the countable-signature module form of the Weil--Moore
strictification argument used in Lemmas~\ref{Lemma:Weil-discrete} and~\ref%
{Lemma:Borel-cocycle}; compare \cite[Theorem 5]{moore_group_1976}. We record
the measure-class point, since it is essential here. Choose Borel
representatives of $c$ and the $\rho _r$. For each arity $m$, discard the
union of the exceptional sets for all the identities %
\eqref{Equation:module-cocycle-normalized}--%
\eqref{Equation:module-cocycle-zero-unit} and for all their substitutions by 
$R$-term maps with domain $C^m$. This is a countable union. Moreover, it is
Haar null in $C^m$: an exceptional set has zero pullback measure because the
corresponding term-pushforward of $\lambda_m$ occurs with positive weight in
the appropriate $\mu_n$.

One now runs the proof of the Weil presentation lemma simultaneously for
this countable family of conull sets. In addition to the generators $(x)$
used in Lemma~\ref{Lemma:Weil-discrete}, impose the relations associated
with the unary terms $x\mapsto rx$, $r\in R$. Closure of the family of term
maps under composition guarantees that every finite compatibility diagram
needed in the presentation is valid on a conull set. The additive part of
the presentation is identified with $C$ by Lemma~\ref{Lemma:Weil-discrete};
under this identification the unary relations give precisely the prescribed $%
R$-action. The standard degree-two part of Moore's argument then makes the
additive defect strict, while the unary relations make the scalar defects
strict. Since the construction changes a cochain only on the exceptional
sets just described, and those sets have zero measure for every
term-pushforward occurring in \eqref{Equation:module-control-measures}, the
new cochains represent the original elements of the $\mathscr{L}_{n}^{0}$%
-spaces. This proves the first assertion.

Applying the same argument to the difference of two strict representatives
gives a Borel $f:C\rightarrow A$, with $f(0)=0$, such that the two pairs
differ by 
\begin{equation*}
(\delta f,(r f-f(r\,\cdot))_{r\in R}).
\end{equation*}%
This is precisely the change of cocycle produced by changing a Borel
section, and hence the corresponding extensions are equivalent.
\end{proof}

\begin{proposition}
\label{Proposition:module-cocycle-cover} Let $C\in\mathbf{LCPMod}_{R}$ and $%
A\in\mathbf{proLiePMod}_{R}$. Then:

\begin{enumerate}
\item $\mathsf{Z}_{R}(C,A)$ is a closed subgroup of 
\begin{equation*}
\mathscr{L}_{2}^{0}(C,A)\oplus \prod_{r\in R}\mathscr{L}_{1}^{0}(C,A),
\end{equation*}%
and is therefore a Polish group;

\item $\mathsf{B}_{R}(C,A)$ is a Polishable subgroup of $\mathsf{Z}_{R}(C,A)$%
; explicitly, its Polish topology is the one transported from 
\begin{equation*}
\mathscr{L}_{1,0}^{0}(C,A)/\ker(D);
\end{equation*}

\item the quotient 
\begin{equation*}
\mathsf{Ext}_{R,\mathrm{coc}}^{1}(C,A):= \mathsf{Z}_{R}(C,A)/\mathsf{B}%
_{R}(C,A)
\end{equation*}%
is a group with a Polish cover.
\end{enumerate}
\end{proposition}

\begin{proof}
Every expression occurring in \eqref{Equation:module-cocycle-normalized}--%
\eqref{Equation:module-cocycle-zero-unit} defines a continuous homomorphism
between countable products of the Polish groups $\mathscr{L}_{n}^{0}(C,A)$.
Thus, $\mathsf{Z}_{R}(C,A)$ is the intersection of countably many kernels of
continuous homomorphisms and is closed. The map $D$ in %
\eqref{Equation:module-coboundary-map} is continuous, and its kernel is
closed. Hence its image becomes a Polish group when equipped with the
quotient topology from $\mathscr{L}_{1,0}^{0}(C,A)/\ker(D)$. The induced
injective homomorphism into $\mathsf{Z}_{R}(C,A)$ is Borel, so the image is
a Polishable subgroup. The final assertion is the definition of a group with
a Polish cover.
\end{proof}

\begin{theorem}[Cocycle description of module extensions]
\label{Theorem:module-Ext-cocycles} Let $R$ be a countable discrete unital
ring, $C$ a locally compact Polish $R$-module, and $A$ a pro-Lie Polish $R$%
-module. There is a natural Borel-definable isomorphism of groups with a
Polish cover 
\begin{equation}  \label{Equation:module-Ext-cocycle-isomorphism}
\mathsf{Ext}_{R,\mathrm{coc}}^{1}(C,A) \cong \mathrm{Ext}_{R}^{1}(C,A).
\end{equation}
Under this isomorphism, forgetting the cochains $(\rho _r)_{r\in R}$ is the
map induced by forgetting the $R$-module structures on extensions.
\end{theorem}

\begin{proof}
Given a strict extension of Polish $R$-modules, choose a Borel section $t$
with $t(0)=0$ and form $c$ and $\rho$ as in %
\eqref{Equation:additive-module-cocycle} and %
\eqref{Equation:scalar-module-cocycle}. Associativity and commutativity in
the middle term give \eqref{Equation:module-cocycle-normalized}--%
\eqref{Equation:module-cocycle-associative}. Additivity of scalar
multiplication gives \eqref{Equation:module-cocycle-additive-action}; the
identities $(r+s)e=re+se$, $(rs)e=r(se)$, $0e=0$, and $1e=e$ give %
\eqref{Equation:module-cocycle-ring-addition}--%
\eqref{Equation:module-cocycle-zero-unit}. Replacing $t$ by $t+\alpha\circ f$
changes the pair by 
\begin{equation*}
(c,\rho)\longmapsto (c,\rho)+(\delta f,(r f-f(r\,\cdot))_{r\in R}).
\end{equation*}%
Thus, the extension determines a class in $\mathsf{Ext}_{R,\mathrm{coc}%
}^{1}(C,A)$.

Conversely, use Lemma~\ref{Lemma:module-cocycle-strictification} to choose a
strict Borel representative $(c,\rho)$. On the standard Borel space $%
E:=A\times C$, define 
\begin{align}
(a,x)+(b,y)&:=\bigl(a+b+c(x,y),x+y\bigr),
\label{Equation:module-extension-addition} \\
r(a,x)&:=\bigl(ra+\rho _r(x),rx\bigr).
\label{Equation:module-extension-action}
\end{align}
Equations \eqref{Equation:module-cocycle-normalized}--%
\eqref{Equation:module-cocycle-associative} make %
\eqref{Equation:module-extension-addition} an abelian group operation, and
the remaining equations say exactly that %
\eqref{Equation:module-extension-action} is a unital $R$-module action by
group endomorphisms. The operations are Borel. The realization theorem for
Borel cocycles with pro-Lie coefficients \cite[Proposition 4.31]%
{casarosa_homological_2026}, applied to the additive cocycle $c$, supplies a
compatible pro-Lie Polish group topology on $E$ for which $A\rightarrow
E\rightarrow C$ is a strict extension. Each scalar map in %
\eqref{Equation:module-extension-action} is then a Borel group homomorphism,
and hence is continuous. Therefore 
\begin{equation*}
0\longrightarrow A\longrightarrow E\longrightarrow C\longrightarrow0, \qquad
a\longmapsto(a,0),\quad(a,x)\longmapsto x,
\end{equation*}%
is a strict extension of Polish $R$-modules.

Adding a coboundary gives the equivalent extension through 
\begin{equation*}
E_{c,\rho}\longrightarrow E_{c+\delta f,\,\rho+\partial f},
\qquad(a,x)\longmapsto(a-f(x),x).
\end{equation*}%
The two constructions are inverse and respect Baer sums. Their graphs are
analytic by the same Borel section argument as in the proof of Proposition~%
\ref{Proposition:Yoneda-Ext}; hence the resulting isomorphism is
Borel-definable. Naturality in $C$ and $A$ follows by pullback and pushout
of extensions.
\end{proof}

\begin{remark}[The case $R=\mathbb{Z}$]
\label{Remark:module-cocycles-Z} For an abelian-group extension, the $%
\mathbb{Z}$-module structure on the middle term is forced by its additive
group structure. Accordingly, the scalar defects are determined by $c$. For
example, for $n\geq1$ one has 
\begin{equation*}
\rho _n(x)=\sum_{j=0}^{n-1}c(jx,x),
\end{equation*}%
with the negative scalars determined by additive inverses. Therefore,
projection $(c,\rho)\mapsto c$ identifies %
\eqref{Equation:module-Ext-cocycle-isomorphism}, for $R=\mathbb{Z}$, with
the cocycle presentation of Proposition~\ref{Proposition:Yoneda-Ext}. The
cochains $\rho_r$ are precisely the additional data needed for a genuine
coefficient-ring action.
\end{remark}

\section{Injective and projective objects in categories of locally compact
abelian groups}

\label{Section:injective}In this section, we apply the machinery developed
so far to study injective and projective objects in the left hearts of the
thick subcategories of $\mathbf{LCPAb}$ introduced in Section \ref%
{Subsection:classes}, as well as relative injectivity and projectivity in
the underlying quasi-abelian categories.

\subsection{Injective and projective groups}

Here is a summary of the main results that will be obtained in this section;
see Theorems \ref{Theorem:no-injectives-LCPAb}, \ref%
{Theorem:no-injectives-TDLCPAb}, \ref{Theorem:injectives-FLCPAb(p)}, \ref%
{Theorem:injectives-TorFLCPAb}, and \ref{Theorem:injectives-LH(TorFLCPAb)}.

\begin{theorem}
\label{Theorem:no-injectives}The trivial group is the only injective object
in the following abelian categories:

\begin{enumerate}
\item $\mathrm{LH}\left( \mathbf{LCPAb}\right) $;

\item $\mathrm{LH}\left( \mathbf{LCPAb}_{\mathrm{cg}}\right) $\textrm{;}

\item $\mathrm{LH}\left( \mathbf{FLCPAb}\right) $;

\item $\mathrm{LH}\left( \mathbf{LieAb}\right) $;

\item $\mathrm{LH}\left( \mathbf{TDLCPAb}\right) $;

\item $\mathrm{LH}\left( \mathbf{TorLCPAb}\right) $;

\item $\mathrm{LH}\left( \mathbf{LCPAb}\left( p\right) \right) $ for any
prime number $p$.
\end{enumerate}

In the categories $\mathbf{TDLCPAb}$ and $\mathbf{TorLCPAb}$, an object is
injective if and only if it is countable and divisible. In the category $%
\mathbf{TorFLCPAb}$ an object is injective if and only if it is of the form $%
D\oplus Z$ for some $D$ and $Z$ in $\mathbf{TorFLCPAb}$ where $D$ is
countable and divisible, and $Z$ is a finite direct sum of copies of $%
\mathbb{Q}_{p}$ with $p$ ranging over the primes. An object of $\mathbf{%
TorFLCPAb}$ is injective in $\mathrm{LH}\left( \mathbf{TorFLCPAb}\right) $
if and only if it is a finite direct sum of copies of $\mathbb{Z}\left(
p^{\infty }\right) $ and $\mathbb{Q}_{p}$, with $p$ ranging over the primes.
\end{theorem}

A characterization of the projective objects of each of these categories
will be given in Section \ref{Subsection:projectives}.

\begin{corollary}
The abelian categories \emph{(1)--(7) }from Theorem \ref%
{Theorem:no-injectives} do not have enough injectives.
\end{corollary}

If $\mathcal{A}$ is a class of locally compact Polish abelian groups, then
we say that a locally compact Polish abelian group $G$ is \emph{injective
for $\mathcal{A}$} if $\mathrm{Ext}\left( X,G\right) =0 $ for every group $X$
in $\mathcal{A}$.

\begin{proposition}
Let $G$ be a locally compact Polish abelian group. Then there exists a
locally compact Polish abelian group $I$ that is injective for $\mathbf{LCPAb%
}$ and an injective continuous homomorphism $G\rightarrow I$.
\end{proposition}

\begin{proof}
Suppose first that $G$ is a compactly generated group. Then there exists a
compactly generated group $I$ that is injective for $\mathbf{LCPAb}$ and an
injective continuous homomorphism $G\rightarrow I$.\ Indeed, $G\cong C\oplus 
\mathbb{R}^{n}\oplus \mathbb{Z}^{m}$ for some $m$ and $n$ and compact $C$,
and it will suffice to consider the cases of each of the summands
separately. If $G$ is a vector group there is nothing to prove.\ If $G=%
\mathbb{Z}^{m}$ then $G$ embeds in a vector group as a closed subgroup.
Suppose that $C$ is compact, and hence that $C^{\vee }$ is countable
discrete. Then there is a surjective map $\mathbb{Z}^{\left( \omega \right)
}\rightarrow C^{\vee }$. This induces an injective map $C=C^{\vee \vee
}\rightarrow \left( \mathbb{Z}^{\left( \omega \right) }\right) ^{\vee }=%
\mathbb{T}^{\omega }$.

Let us consider now the case when $G$ is countable. We claim that there
exists an injective homomorphism $G\rightarrow \mathbb{T}^{\omega }$. After
replacing $G$ with its divisible hull $D\left( G\right) $, we can assume
without loss of generality that $G$ is divisible. Thus, it suffices to
consider the cases when $G=\mathbb{Q}$ or $G=\mathbb{Z}\left( p^{\infty
}\right) $. By definition, $\mathbb{Z}\left( p^{\infty }\right) $ is a
subgroup of $\mathbb{T}$. Furthermore, if $\alpha \in \mathbb{R}\setminus 
\mathbb{Q}$, then the homomorphism $\mathbb{Q}\rightarrow \mathbb{T}$, $%
x\mapsto \alpha x$ mod $\mathbb{Z}$, is injective.

Suppose now that $G$ is an arbitrary locally compact Polish abelian group
(i.e., one no longer assumed to be compactly generated). Then $G$ has a
compactly generated open subgroup $U$, and $G/U=D$ is countable discrete. By
the above, we have that $U$ admits an injective continuous homomorphism $%
U\rightarrow I$ for some compactly generated locally compact Polish abelian
group $I$ that is injective for $\mathbf{LCPAb}$. Consider the induced
extension (via pushout)

\begin{tikzcd}
U \arrow[r] \arrow[d] &   G \arrow[r] \arrow[d]  &  D \arrow[d, equal] \\ 
I \arrow[r]  &  Z \arrow[r] &  D   
\end{tikzcd}

where the map $G\rightarrow Z$ is injective and the map $Z\rightarrow D$ is
induced by the universal property of the pushout from the identity map of $D$
and the zero map $I\rightarrow D$. Since $I$ is injective for $\mathbf{LCPAb}
$, we have that $Z\cong I\oplus D$. Furthermore, since $D$ is countable, we
have that $D\subseteq \mathbb{T}^{\omega }$ and hence $G\subseteq I\oplus
D\subseteq I\oplus \mathbb{T}^{\omega }$.
\end{proof}

\subsection{Groups with countable $\mathrm{Hom}$ and $\mathrm{Ext}$}

In this subsection we record several lemmas concerning groups that have
countable $\mathrm{Hom}$ and $\mathrm{Ext}$ groups.\newline

\begin{lemma}
\label{Lemma:Ext(S1,Z)}If $A,B$ are countable finite-rank torsion-free
groups, then $\mathrm{Ext}(A^{\vee },B)$ is countable.
\end{lemma}

\begin{proof}
We begin by considering the case $A=\mathbb{Q}$. In this case, we have a
short exact sequence $\mathbb{Q}\rightarrow \mathbb{A}\rightarrow \mathbb{Q}%
^{\vee }$ and $\mathrm{Ext}(\mathbb{A},\mathbb{Q})=0$ by \cite[Proposition
4.15(vi)]{hoffmann_homological_2007}. Using Remark \ref{Remark:ext-dual}, we
also have that $\mathrm{Hom}(\mathbb{A},\mathbb{Q})\cong \mathrm{Hom}(%
\mathbb{Q}^{\vee },\mathbb{A})=0$ since $\mathbb{A}_{\mathbb{S}^{1}}=0$. It
is now immediate from the relevant part of the $\mathrm{Hom}(-,\mathbb{Q})$
exact sequence 
\begin{equation*}
0=\mathrm{Hom}(\mathbb{A},\mathbb{Q})\rightarrow \mathbb{Q}\cong \mathrm{Hom}%
(\mathbb{Q},\mathbb{Q})\rightarrow \mathrm{Ext}(\mathbb{Q}^{\vee },\mathbb{Q}%
)\rightarrow \mathrm{Ext}(\mathbb{A},\mathbb{Q})=0\text{.}
\end{equation*}%
that $\mathrm{Ext}(\mathbb{Q}^{\vee },\mathbb{Q})\cong \mathbb{Q}$ and $%
\mathrm{Ext}(\mathbb{Q}^{\vee },\mathbb{Q}^{n})\cong \mathbb{Q}^{n}$.

If $A$ is finite rank torsion-free, then there is a short exact sequence $%
A\rightarrow \mathbb{Q}^{n}\rightarrow T$ where $T$ is torsion. We have%
\begin{equation*}
\mathrm{Hom}(\mathbb{Q}^{\vee },T)=0
\end{equation*}%
since $\mathbb{Q}^{\vee }$ is connected. The exact sequence%
\begin{equation*}
0=\mathrm{Hom}(\mathbb{Q}^{\vee },T)\rightarrow \mathrm{Ext}(\mathbb{Q}%
^{\vee },A)\rightarrow \mathrm{Ext}(\mathbb{Q}^{\vee },\mathbb{Q}^{n})=%
\mathbb{Q}^{n}
\end{equation*}%
then implies that $\mathrm{Ext}(\mathbb{Q}^{\vee },A)\cong \mathrm{Ext}%
(A^{\vee },\mathbb{Q})$ is countable.

Suppose now that $B$ is also a finite-rank torsion-free group. Then there is
a short exact sequence $B\rightarrow \mathbb{Q}^{m}\rightarrow T$ with $T$
torsion. This induces a short exact sequence%
\begin{equation*}
\mathrm{Hom}\left( A^{\vee },T\right) \rightarrow \mathrm{Ext}\left( A^{\vee
},B\right) \rightarrow \mathrm{Ext}(A^{\vee },\mathbb{Q}^{m})\cong \mathrm{%
Ext}(A^{\vee },\mathbb{Q})^{m}\text{.}
\end{equation*}%
$\mathrm{Ext}(A^{\vee },\mathbb{Q})^{m}$ is countable by the remarks above,
and $\mathrm{Hom}\left( A^{\vee },T\right) $ is countable because $A^{\vee }$
is compact and $T$ is countable.
\end{proof}

\begin{lemma}
\label{Lemma:Ext(torsion,countable)}If $H$ is a countable abelian group and $%
K$ is a compact topological torsion abelian Polish group, then $\mathrm{Ext}%
\left( K,H\right) $ is countable.
\end{lemma}

\begin{proof}
Recall that a locally compact Polish abelian group is codivisible if its
Pontryagin dual is divisible. By (the dual of) \cite[Proposition 3.8(i)]%
{hoffmann_homological_2007}, we can consider a \emph{codivisible resolution}
for $K$, namely a short exact sequence 
\begin{equation*}
0\rightarrow C^{\prime }\rightarrow C\rightarrow K\rightarrow 0
\end{equation*}%
where $C,C^{\prime }$ are codivisible compact topological torsion groups. We
may also consider a \emph{divisible resolution} for $H$, namely a short
exact sequence%
\begin{equation*}
0\rightarrow H\rightarrow D\rightarrow D^{\prime }\rightarrow 0
\end{equation*}%
where $D,D^{\prime }$ are divisible countable groups and $D^{\prime }$ is
torsion. We have then an exact sequence%
\begin{equation*}
\mathrm{Hom}\left( C,D^{\prime }\right) \rightarrow \mathrm{Ext}\left(
C,H\right) \rightarrow \mathrm{Ext}\left( C,D\right)
\end{equation*}%
Since $C$ is compact and $D^{\prime }$ is countable, $\mathrm{Hom}\left(
C,D^{\prime }\right) $ is countable. This fact and $C_{\mathbb{S}^{1}}=0$
imply that $\mathrm{Ext}\left( C,D\right) =0$ by \cite[Proposition 4.15(iii)]%
{hoffmann_homological_2007}. Thus, $\mathrm{Ext}\left( C,H\right) $ is
countable. We also have an exact sequence%
\begin{equation*}
\mathrm{Hom}\left( C^{\prime },H\right) \rightarrow \mathrm{Ext}\left(
K,H\right) \rightarrow \mathrm{Ext}\left( C,H\right)
\end{equation*}%
As before, $\mathrm{Hom}\left( C^{\prime },H\right) $ is countable, hence $%
\mathrm{Ext}\left( K,H\right) $ is countable.
\end{proof}

\begin{lemma}
\label{Lemma:Ext(compact,countable)}If $H$ is a countable abelian group with
finite $\mathbb{Z}$-rank and $K$ is a compact Polish abelian group with
finite $\mathbb{S}^{1}$-rank, then $\mathrm{Ext}\left( K,H\right) $ is
countable.
\end{lemma}

\begin{proof}
We have an exact sequence%
\begin{equation*}
\mathrm{Ext}\left( K,H_{\mathrm{t}}\right) \rightarrow \mathrm{Ext}\left(
K,H\right) \rightarrow \mathrm{Ext}\left( K,H_{\mathbb{Z}}\right)
\end{equation*}%
as well as exact sequences 
\begin{equation*}
\mathrm{Ext}\left( K_{\mathrm{t}},H_{\mathbb{Z}}\right) \rightarrow \mathrm{%
Ext}\left( K,H_{\mathbb{Z}}\right) \rightarrow \mathrm{Ext}\left( K_{\mathbb{%
S}^{1}},H_{\mathbb{Z}}\right)
\end{equation*}%
and%
\begin{equation*}
\mathrm{Ext}\left( K_{\mathrm{t}},H_{\mathrm{t}}\right) \rightarrow \mathrm{%
Ext}\left( K,H_{\mathrm{t}}\right) \rightarrow \mathrm{Ext}\left( K_{\mathbb{%
S}^{1}},H_{\mathrm{t}}\right) .
\end{equation*}%
By Lemma \ref{Lemma:Ext(torsion,countable)}, $\mathrm{Ext}\left( K_{\mathrm{t%
}},H_{\mathbb{Z}}\right) $ and $\mathrm{Ext}\left( K_{\mathrm{t}},H_{\mathrm{%
t}}\right) $ are countable. By Pontryagin duality and the same lemma, 
\begin{equation*}
\mathrm{Ext}\left(K_{\mathbb{S}^{1}},H_{\mathrm{t}}\right) \cong\mathrm{Ext}%
\left(H_{\mathrm{t}}^{\vee}, K_{\mathbb{S}^{1}}^{\vee}\right)
\end{equation*}%
is countable. Also, $\mathrm{Ext}\left( K_{\mathbb{S}^{1}},H_{\mathbb{Z}%
}\right) $ is countable by Lemma \ref{Lemma:Ext(S1,Z)}. Therefore, $\mathrm{%
Ext}\left( K,H_{\mathbb{Z}}\right) $ and $\mathrm{Ext}\left( K,H_{\mathrm{t}%
}\right) $ are countable, and $\mathrm{Ext}\left( K,H\right) $ is countable,
as claimed.
\end{proof}

\begin{lemma}
\label{Lemma:Ext-Prufer}The following chain of isomorphisms holds: 
\begin{equation*}
\mathbb{Z}_{p}\cong \mathrm{Ext}\left( \mathbb{Z}\left( p^{\infty }\right) ,%
\mathbb{Z}_{p}\right) \cong \mathrm{Hom}(\mathbb{Z}_{p},\mathbb{Z}_{p})\cong 
\mathrm{Hom}\left( \mathbb{Z}\left( p^{\infty }\right) ,\mathbb{Z}\left(
p^{\infty }\right) \right) \cong \mathrm{Ext}\left( \mathbb{Z}\left(
p^{\infty }\right) ,\mathbb{Z}\right)
\end{equation*}
\end{lemma}

\begin{proof}
By \cite[Proposition 4.15(iii)]{hoffmann_homological_2007}, $\mathrm{Ext}(%
\mathbb{Q}_{p},\mathbb{Q}_{p})=0$ and $\mathrm{Ext}(\mathbb{Q}_{p},\mathbb{Z}%
\left( p^{\infty }\right) )=0$. We have an exact sequence%
\begin{equation*}
\mathrm{Ext}(\mathbb{Q}_{p},\mathbb{Z}\left( p^{\infty }\right) )\rightarrow 
\mathrm{Ext}(\mathbb{Q}_{p},\mathbb{Z}_{p})\rightarrow \mathrm{Ext}(\mathbb{Q%
}_{p},\mathbb{Q}_{p})
\end{equation*}%
Therefore $\mathrm{Ext}(\mathbb{Q}_{p},\mathbb{Z}_{p})=0$. By duality, this
implies by Remark \ref{Remark:ext-dual} that $\mathrm{Ext}(\mathbb{Z}\left(
p^{\infty }\right) ,\mathbb{Q}_{p})=0$.

For $z\in \mathbb{Z}_{p}$ define $\varphi _{z}\in \mathrm{Hom}\left( \mathbb{%
Z}_{p},\mathbb{Z}_{p}\right) $ by setting $\varphi _{z}(x)=zx$. If $\varphi
\in \mathrm{Hom}\left( \mathbb{Z}_{p},\mathbb{Z}_{p}\right) $, then setting $%
z:=\varphi \left( 1\right) $ we have that $\varphi \left( n\right)
=nz=\varphi _{z}\left( n\right) $ for every $n\in \mathbb{Z}$ and hence by
density of $\mathbb{Z}$ in $\mathbb{Z}_{p}$ we have that $\varphi =\varphi
_{z}$. Thus, the assignment $z\mapsto \varphi _{z}$ defines a topological
isomorphism $\mathbb{Z}_{p}\rightarrow \mathrm{Hom}\left( \mathbb{Z}_{p},%
\mathbb{Z}_{p}\right) $.

We also have exact sequences%
\begin{equation*}
0=\mathrm{Hom}(\mathbb{Q}_{p},\mathbb{Z}_{p})\rightarrow \mathrm{Hom}\left( 
\mathbb{Z}_{p},\mathbb{Z}_{p}\right) \rightarrow \mathrm{Ext}\left( \mathbb{Z%
}\left( p^{\infty }\right) ,\mathbb{Z}_{p}\right) \rightarrow \mathrm{Ext}(%
\mathbb{Q}_{p},\mathbb{Z}_{p})=0
\end{equation*}%
\begin{equation*}
\begin{aligned} 0=\mathrm{Hom}(\mathbb{Z}\left( p^{\infty }\right)
,\mathbb{Q}_{p}) &\rightarrow \mathrm{Hom}\left( \mathbb{Z}\left( p^{\infty
}\right) , \mathbb{Z}\left( p^{\infty }\right) \right) \\ &\rightarrow
\mathrm{Ext}\left( \mathbb{Z}\left( p^{\infty }\right) ,
\mathbb{Z}_{p}\right) \rightarrow \mathrm{Ext}(\mathbb{Z}\left( p^{\infty
}\right) ,\mathbb{Q}_{p})=0. \end{aligned}
\end{equation*}%
This concludes the proof.
\end{proof}

\begin{lemma}
\label{Lemma:robertson}Suppose that $A$ is a countable $p$-group. Then $%
\mathrm{Ext}\left( A,\mathbb{Z}\left( p\right) \right) =0$ if and only if $%
A=0$.
\end{lemma}

\begin{proof}
Suppose that $\mathrm{Ext}\left( A,\mathbb{Z}\left( p\right) \right) =0$.
Define $C=A^{\vee }$. Using Remark \ref{Remark:ext-dual}, we have that $%
\mathrm{Ext}\left( A,\mathbb{Z}\left( p\right) \right) \cong \mathrm{Ext}%
\left( \mathbb{Z}\left( p\right) ,C\right) \cong C/pC$. Thus $C$ is a
divisible pro-$p$ group. Hence $A=C^{\vee }$ is torsion-free, and therefore
equals $0$, by Robertson's Theorem \cite[Theorem 4.15]%
{armacost_structure_1981} asserting that a locally compact Polish abelian
group has a dense divisible subgroup if and only if its dual is torsion-free.
\end{proof}

\begin{lemma}
\label{Lemma:uncountable-Ext}Let $A$ be a countable reduced $p$-group. If $A$
does not have finite $p$-rank, then $\mathrm{Ext}\left( A,\mathbb{Z}\left(
p\right) \right) $ is uncountable.
\end{lemma}

\begin{proof}
By the basic subgroup theorem \cite[Sections 32--33]{fuchs_infinite_1970}, $%
A $ has a pure subgroup $B$ such that $A/B$ is divisible and%
\begin{equation*}
B\cong \bigoplus_{i\in I}B_{i},
\end{equation*}%
where every $B_{i}$ is a nonzero finite cyclic $p$-group. The set $I$ must
be infinite. Indeed, if $I$ were finite, then $B$ would be finite. Since a
finite pure subgroup is a direct summand, we would have $A\cong B\oplus A/B$%
; the reducedness of $A$ would then force $A/B=0$, contrary to the
assumption that $A$ does not have finite $p$-rank.

The short exact sequence%
\begin{equation*}
0\rightarrow B\rightarrow A\rightarrow A/B\rightarrow 0
\end{equation*}%
induces an exact sequence%
\begin{equation*}
\mathrm{Ext}\left( A,\mathbb{Z}\left( p\right) \right) \rightarrow \mathrm{%
Ext}\left( B,\mathbb{Z}\left( p\right) \right) \rightarrow \mathrm{Ext}%
^{2}\left( A/B,\mathbb{Z}\left( p\right) \right) =0.
\end{equation*}%
The last term vanishes by the higher-$\mathrm{Ext}$ vanishing established in
Subsection \ref{Subsection:derived-lc}. Consequently, the first map is
surjective. Moreover,%
\begin{equation*}
\mathrm{Ext}\left( B,\mathbb{Z}\left( p\right) \right) \cong \prod_{i\in I}%
\mathrm{Ext}\left( B_{i},\mathbb{Z}\left( p\right) \right) \text{.}
\end{equation*}%
For each $i\in I$, we have $B_{i}\cong \mathbb{Z}/p^{n_{i}}\mathbb{Z}$ for
some $n_{i}\geq 1$, and hence%
\begin{equation*}
\mathrm{Ext}\left( B_{i},\mathbb{Z}\left( p\right) \right) \cong \mathbb{Z}%
\left( p\right) \text{.}
\end{equation*}%
Since $I$ is infinite, the product is uncountable. Hence its surjective
preimage $\mathrm{Ext}\left( A,\mathbb{Z}\left( p\right) \right) $ is
uncountable.
\end{proof}

\begin{lemma}
\label{Lemma:countable-Ext(p)}Suppose that $A$ is a countable $p$-group.
Then we have that $\mathrm{Ext}\left( A,\mathbb{Z}\left( p\right) \right) $
is countable if and only if $A$ has finite $p$-rank.
\end{lemma}

\begin{proof}
By \cite[Lemma 2.8(iii)]{hoffmann_homological_2007}, a finite $p$-rank
topological $p$-group is a finite direct sum of finite cyclic $p$-groups and
copies of $\mathbb{Z}_{p}$, $\mathbb{Q}_{p}$, and $\mathbb{Z}\left(
p^{\infty }\right) $. Since $A$ is countable and discrete, the summands $%
\mathbb{Z}_{p}$ and $\mathbb{Q}_{p}$ cannot occur. Using Remark \ref%
{Remark:ext-dual}, we have that%
\begin{equation*}
\mathrm{Ext}\left( \mathbb{Z}(p^{k}),\mathbb{Z}\left( p\right) \right) =%
\mathbb{Z}\left( p\right)
\end{equation*}%
for $k\geq 1$, and%
\begin{equation*}
\mathrm{Ext}\left( \mathbb{Z}\left( p^{\infty }\right) ,\mathbb{Z}\left(
p\right) \right) \cong \mathrm{Ext}\left( \mathbb{Z}\left( p\right) ,\mathbb{%
Z}_{p}\right) \cong \mathbb{Z}_{p}/p\mathbb{Z}_{p}\cong \mathbb{Z}\left(
p\right) \text{.}
\end{equation*}%
This concludes the proof that, if $A$ has finite $p$-rank, then $\mathrm{Ext}%
\left( A,\mathbb{Z}\left( p\right) \right) $ is countable.

Conversely, suppose that $\mathrm{Ext}\left( A,\mathbb{Z}\left( p\right)
\right) $ is countable. We can write $A=D\left( A\right) \oplus B$ where $B$
is reduced, and $D\left( A\right) =\mathbb{Z}\left( p^{\infty }\right)
^{\left( k\right) }$ for some $k\leq \omega $. Then we have that 
\begin{eqnarray*}
\mathrm{Ext}\left( A,\mathbb{Z}\left( p\right) \right) &\cong &\mathrm{Ext}(%
\mathbb{Z}\left( p^{\infty }\right) ^{\left( k\right) },\mathbb{Z}\left(
p\right) )\oplus \mathrm{Ext}\left( B,\mathbb{Z}\left( p\right) \right) \\
&\cong &\mathbb{Z}\left( p\right) ^{k}\oplus \mathrm{Ext}\left( B,\mathbb{Z}%
\left( p\right) \right) \text{.}
\end{eqnarray*}%
This forces $k$ to be finite.\ Furthermore, $\mathrm{Ext}\left( B,\mathbb{Z}%
\left( p\right) \right) $ must also be countable, and hence $B$ has finite $%
p $-rank by Lemma \ref{Lemma:uncountable-Ext}.
\end{proof}

\begin{lemma}
\label{Lemma:torsion-countable-Ext}Suppose that $A$ is a countable torsion
group.

\begin{enumerate}
\item If $\mathrm{Ext}\left( A,\mathbb{Z}\right) $ is countable, then $A$ is
finite.

\item If $\mathrm{Ext}\left( A,T\right) $ is countable for every countable
torsion group $T$ with cyclic primary components, then the $p$-component $%
A_{p}$ of $A$ equals $0$ for all but finitely many primes $p$.

\item If either $\mathrm{Ext}\left( A,\mathbb{Z}_{p}\right) $ or $\mathrm{Ext%
}(A,T)$ is countable, where $T=\bigoplus_{n}\left( \mathbb{Z}/p^{n}\right) $%
, then the $p$-component%
\begin{equation*}
A_{p}:=\left\{ x\in A:\exists n\in \omega \text{, }p^{n}x=0\right\}
\end{equation*}%
of $A$ is reduced.
\end{enumerate}
\end{lemma}

\begin{proof}
(1) We have that%
\begin{equation*}
\mathrm{Ext}\left( A,\mathbb{Z}\right) \cong \mathrm{Ext}\left( \mathbb{T}%
,A^{\vee }\right) \cong A^{\vee }\text{;}
\end{equation*}%
see Remark \ref{Remark:ext-dual}. Thus, if $\mathrm{Ext}\left( A,\mathbb{Z}%
\right) $ is countable then $A^{\vee }$ is countable discrete, hence $A$ is
compact, and therefore finite.

(2) Suppose that there is an infinite set $I$ of primes such that $A_{p}$ is
nonzero for $p\in I$. Then setting $T:=\bigoplus_{p\in I}\mathbb{Z}\left(
p\right) $ we have that%
\begin{equation*}
\mathrm{Ext}\left( A,T\right) \cong \prod_{p\in I}\mathrm{Ext}\left( A_{p},%
\mathbb{Z}\left( p\right) \right) ,
\end{equation*}%
which is uncountable by Lemma \ref{Lemma:robertson}.

(3) It suffices to show that $A_{p}$ does not have $\mathbb{Z}\left(
p^{\infty }\right) $ as a direct summand. This follows from the fact that $%
\mathrm{Ext}\left( \mathbb{Z}\left( p^{\infty }\right) ,\mathbb{Z}%
_{p}\right) \cong \mathbb{Z}_{p}$ by Lemma \ref{Lemma:Ext-Prufer}. We also
know that $\mathrm{Ext}(\mathbb{Z}\left( p^{\infty }\right) ,T)$ is
uncountable where $T\ $is an unbounded reduced $p$-group by \cite[Theorem 4.4%
]{lupini_projective_2025}.
\end{proof}

\begin{lemma}
\label{Lemma:finite-ranks-countable-Ext}Suppose that $A$ is a countable
abelian group. If $\mathrm{Ext}\left( A,\mathbb{Z}\right) $ is countable
then $A\cong F\oplus K$ where $K$ is finite and $F$ is free abelian.
\end{lemma}

\begin{proof}
Conclude from the exact sequence%
\begin{equation*}
0=\mathrm{Hom}\left( A_{\mathrm{t}},\mathbb{Z}\right) \rightarrow \mathrm{Ext%
}\left( A_{\mathbb{Z}},\mathbb{Z}\right) \rightarrow \mathrm{Ext}\left( A,%
\mathbb{Z}\right)
\end{equation*}%
that $\mathrm{Ext}\left( A_{\mathbb{Z}},\mathbb{Z}\right) $ is countable.
After replacing $A$ with $A_{\mathbb{Z}}$, we can assume that $A$ is
torsion-free. In this case $\mathrm{Ext}\left( A_{\mathbb{Z}},\mathbb{Z}%
\right) =\mathrm{PExt}\left( A_{\mathbb{Z}},\mathbb{Z}\right) $ and $\left\{
0\right\} $ is dense in $\mathrm{PExt}\left( A_{\mathbb{Z}},\mathbb{Z}%
\right) $. Thus, if $\mathrm{Ext}\left( A_{\mathbb{Z}},\mathbb{Z}\right) $
is countable, then it is trivial. Hence, it follows from \cite[Theorem 3.2]%
{friedenberg_extensions_2013} that $A_{\mathbb{Z}}$ is a countable free
abelian group.

Thus, we have that $A\cong A_{\mathrm{t}}\oplus A_{\mathbb{Z}}$ and $\mathrm{%
Ext}\left( A_{\mathrm{t}},\mathbb{Z}\right) $ is countable. By Lemma \ref%
{Lemma:torsion-countable-Ext}, $A_{\mathrm{t}}$ is finite; this concludes
the proof.
\end{proof}

\begin{corollary}
\label{Corollary:finite-ranks-countable-Ext}Suppose that $L$ is a compact
abelian group. Suppose that $\mathrm{Ext}\left( \mathbb{T},L\right) $ is
countable. Then $L\cong C\oplus K$ where $K$ is finite and $C$ is a torus.
\end{corollary}

Suppose that $T$ is a countable torsion group. Its first Ulm subgroup is
defined to be%
\begin{equation*}
u_{1}\left( T\right) :=\left\{ x\in T:\forall n\in \mathbb{N}\text{, }x\text{
is }n\text{-divisible}\right\} \text{.}
\end{equation*}%
One defines recursively $u_{\alpha }\left( T\right) $ for every $\alpha
<\omega _{1}$ by setting $u_{\sigma +1}\left( T\right) :=u_{1}\left(
u_{\sigma }\left( T\right) \right) $ for every countable ordinal $\sigma $,
and $u_{\lambda }\left( T\right) :=\bigcap_{\beta <\lambda }u_{\beta }\left(
T\right) $ for every countable limit ordinal $\lambda $. The Ulm rank of $T$
is the least $\alpha <\omega _{1}$ such that $u_{\alpha }\left( T\right)
=u_{\alpha +1}\left( T\right) $, in which case $u_{\alpha }\left( T\right) $
is the largest divisible subgroup of $T$. A countable reduced torsion group
has Ulm rank $1$ if and only if it is a countable direct sum of cyclic
groups \cite[Theorem 17.3]{fuchs_infinite_1970}.

\begin{lemma}
\label{Lemma:torsion-countable-Ext-2}Suppose that $A$ is a countable torsion
group. Suppose that $\mathrm{Ext}\left( A,C\right) $ is countable for every
countable torsion group $C$ of Ulm length at most $1$. Then $A$ is finite.
\end{lemma}

\begin{proof}
By Lemma \ref{Lemma:torsion-countable-Ext}, we have that $A_{p}=0$ for all
but finitely many primes $p$. Furthermore, for every prime $p$, $A_{p}$ has
finite rank and it is reduced. Hence, $A$ is finite.
\end{proof}

\begin{lemma}
\label{Lemma:torsion-free-countable-Ext-2}Suppose that $A$ is a countable
torsion-free group. Suppose that $\mathrm{Ext}\left( A,C\right) $ is
countable for every countable torsion group of Ulm length at most $1$. Then $%
A$ is free.
\end{lemma}

\begin{proof}
By Pontryagin's Theorem, $A$ is free if and only if every finite-rank
subgroup of $A$ is free \cite[Theorem 19.1]{fuchs_infinite_1970}. Since
every finite-rank subgroup of $A$ is torsion-free, we may therefore without
loss of generality assume that $A$ has finite rank $r$.

Let $E\subseteq A$ be an essential free subgroup. For every countable
torsion group $C$ of Ulm rank at most $1$ we have an exact sequence%
\begin{equation*}
C^{r}\cong \mathrm{Hom}\left( E,C\right) \rightarrow \mathrm{Ext}\left(
A/E,C\right) \rightarrow \mathrm{Ext}\left( A,C\right)
\end{equation*}%
Since $C^{r}$ is countable and $\mathrm{Ext}\left( A,C\right) $ is, by
assumption, countable, $\mathrm{Ext}\left( A/E,C\right) $ is countable. By
Lemma \ref{Lemma:torsion-countable-Ext-2} this implies that $A/E$ is finite.
Thus, $A$ is torsion-free and finitely generated. This implies that $A$ is
free.
\end{proof}

\begin{lemma}
\label{Lemma:mixed-countable-Ext2}Suppose that $A$ is a countable group.
Suppose that $\mathrm{Ext}\left( A,C\right) $ is countable for every
countable torsion group of Ulm length at most $1$. Then $A\cong F\oplus K$
where $K$ is finite and $F$ is free abelian.
\end{lemma}

\begin{proof}
Consider the short exact sequence $A_{\mathrm{t}}\rightarrow A\rightarrow A_{%
\mathbb{Z}}$. For every countable torsion group $C$ of Ulm length at most $1$%
, we have an exact sequence%
\begin{equation*}
\mathrm{Ext}\left( A,C\right) \rightarrow \mathrm{Ext}\left( A_{\mathrm{t}%
},C\right) \rightarrow 0
\end{equation*}%
Thus $A_{\mathrm{t}}$ satisfies the assumptions of Lemma \ref%
{Lemma:torsion-countable-Ext-2}. Therefore, $A_{\mathrm{t}}$ is finite. We
also have an exact sequence%
\begin{equation*}
\mathrm{Hom}\left( A_{\mathrm{t}},C\right) \rightarrow \mathrm{Ext}\left( A_{%
\mathbb{Z}},C\right) \rightarrow \mathrm{Ext}\left( A,C\right)
\end{equation*}%
where $\mathrm{Ext}\left( A,C\right) $ is countable by hypothesis, and $%
\mathrm{Hom}\left( A_{\mathrm{t}},C\right) $ is countable since $A_{\mathrm{t%
}}$ is finite. Hence $A_{\mathbb{Z}}$ satisfies the assumptions of Lemma \ref%
{Lemma:torsion-free-countable-Ext-2} and is therefore free abelian. It
follows that $A\cong A_{\mathrm{t}}\oplus A_{\mathbb{Z}}$, concluding the
proof.
\end{proof}

\begin{corollary}
\label{Corollary:product-cyclic-countable-Ext}Suppose that $L$ is a compact
Polish abelian group. Suppose that $\mathrm{Ext}\left( B,L\right) $ is
countable whenever $B$ is a countable product of finite cyclic groups. Then $%
L\cong C\oplus K$ where $K$ is finite and $C$ is a torus.
\end{corollary}

\subsection{Injectives in the left heart of locally compact groups}

We begin with a \emph{structure theorem} providing a \emph{normal form }for
objects in the left heart of the category of compactly generated locally
compact Polish abelian groups.

\begin{proposition}
\label{Proposition:characterize-LH(LCPABcg)}An object of $\mathrm{LH}\left( 
\mathbf{LCPAb}_{\mathrm{cg}}\right) $ is of the form%
\begin{equation*}
\frac{C}{E\oplus W}\oplus F\oplus V
\end{equation*}%
where $C$ is a compact Polish abelian group, $V,W$ are finite-dimensional
vector groups, and $E,F$ are finite-rank free abelian groups.
\end{proposition}

\begin{proof}
Consider $G/H\in \mathrm{LH}\left( \mathbf{LCPAb}_{\mathrm{cg}}\right) $.
After replacing $H$ with $H/U$ and $G$ with $G/U$ for some compact open
subgroup $U$ of $H$, we can assume that $H\cong E\oplus W$ and $G\cong
C\oplus F\oplus V$ for some compact abelian Polish group $C$ and finite-rank
torsion-free abelian groups $E,F$ and vector groups $V,W$.

Let $\left( e_{1},\ldots ,e_{a}\right) $ be a $\mathbb{Z}$-basis of $E$ and $%
\left( f_{1},\ldots ,f_{b}\right) $ an $\mathbb{R}$-vector space basis for $%
W $. For $1\leq i\leq a$, let $\left\langle e_{i}\right\rangle $ be the
subgroup of $E$ generated by $e_{i}$. If the homomorphism%
\begin{equation*}
\left\langle e_{i}\right\rangle \rightarrow H\rightarrow G\rightarrow
F\oplus V
\end{equation*}

\noindent is nonzero, then it is injective, and the homomorphism%
\begin{equation*}
\left\langle e_{i}\right\rangle \rightarrow H\rightarrow G
\end{equation*}%
\noindent has closed image. Thus, after replacing $H$ with $H/\left\langle
e_{i}\right\rangle $ and $G$ with $G/\left\langle e_{i}\right\rangle $, we
can assume that the homomorphism%
\begin{equation*}
\left\langle e_{i}\right\rangle \rightarrow H\rightarrow G\rightarrow
F\oplus V
\end{equation*}

is zero. Repeating this procedure for $1\leq i\leq a$ we conclude that we
can assume that $E\subseteq C$. The same argument applied to the vector
subspace $\left\langle f_{j}\right\rangle _{\mathbb{R}}$ spanned by $f_{j}$
for $1\leq j\leq b$ shows that we can assume without loss of generality that 
$W\subseteq C$. Thus, we have that $H\subseteq C$ and 
\begin{equation*}
G/H\cong \left( C/H\right) \oplus F\oplus V
\end{equation*}%
This concludes the proof.
\end{proof}

\begin{proposition}
\label{Lemma:FLCPAbcg-injective-1}Suppose that $C$ is a compact Polish
abelian group, and consider a group with a Polish cover $X=C/\left( F\oplus
V\right) $ where $F$ is a finite-rank free abelian group and $V$ is a vector
group.\ If $\mathrm{Ext}\left( \mathbb{T},X\right) =0$, then $X\cong
T/(H\oplus W)$ for some torus $T$, free finite-rank subgroup $H$ of $T$, and
vector group $W$.
\end{proposition}

\begin{proof}
Our assumptions determine an exact sequence%
\begin{equation*}
F\cong \mathrm{Ext}\left( \mathbb{T},F\oplus V\right) \rightarrow \mathrm{Ext%
}\left( \mathbb{T},C\right) \rightarrow \mathrm{Ext}\left( \mathbb{T}%
,C/\left( F\oplus V\right) \right) =0.
\end{equation*}%
Hence $\mathrm{Ext}\left( \mathbb{T},C\right) $ is countable. Therefore, by
Corollary \ref{Corollary:finite-ranks-countable-Ext} we have that%
\begin{equation*}
C\cong K\oplus T
\end{equation*}%
for some torus group $T$ and finite group $K$. Proceeding as in the proof of
Proposition \ref{Proposition:characterize-LH(LCPABcg)}, we may assume that
the map $F\rightarrow C\rightarrow K$ is zero. Furthermore, the image of $%
V\rightarrow K\oplus T$ is contained in $T$. Thus 
\begin{equation*}
X=\left( T/(H\oplus W)\right) \oplus K
\end{equation*}%
where $H$ is the image of $F$ in $T$ and $W$ is the image of $V$ in $T$.
Since $\mathrm{Ext}\left( \mathbb{T},X\right) =0$ and $K$ is a direct
summand of $X$, we have that $\mathrm{Ext}\left( \mathbb{T},K\right) =0$.
This then implies that $K=0$ by Proposition \ref{Proposition:injectives},
concluding the proof.
\end{proof}

\begin{proposition}
\label{Proposition:relative-injectives-LH(FLCPAbgc)}Let $X$ be an object of $%
\mathrm{LH}(\mathbf{LCPAb}_{\mathrm{cg}})$. The following assertions are
equivalent:

\begin{enumerate}
\item $X$ is injective for $\mathbf{LCPAb}$;

\item $\mathrm{Ext}\left( \mathbb{T},X\right) =0$;

\item $X$ is isomorphic to $V\oplus T/(E\oplus W)$ for some vector groups $%
V,W$, torus $T$, and finite-rank free subgroup $E$ of $T$.
\end{enumerate}
\end{proposition}

\begin{proof}
The implication (1)$\Rightarrow $(2) is trivial.

(3)$\Rightarrow $(1) Recall that vector groups and torus groups are
injective for $\mathbf{LCPAb}$ by Proposition \ref{Proposition:injectives}.
If $X=V\oplus T/(E\oplus W)$ as in (3), then for every locally compact
Polish abelian group $C$ we have an exact sequence%
\begin{equation*}
0=\mathrm{Ext}\left( C,T\right) \rightarrow \mathrm{Ext}\left( C,T/\left(
E\oplus W\right) \right) \rightarrow \mathrm{Ext}^{2}\left( C,E\oplus
W\right) =0
\end{equation*}%
and hence $\mathrm{Ext}\left( C,T/\left( E\oplus W\right) \right) =0$.
Furthermore, $\mathrm{Ext}\left( C,V\right) =0$.\ Therefore, we have%
\begin{equation*}
\mathrm{Ext}\left( C,X\right) =\mathrm{Ext}\left( C,T/\left( E\oplus
W\right) \right) \oplus \mathrm{Ext}\left( C,V\right) =0\text{.}
\end{equation*}%
As this holds for every locally compact Polish abelian group $C$, this shows
that (1) holds.

(2)$\Rightarrow $(3) Let $X$ be an object of the left heart of $\mathbf{LCPAb%
}_{\mathrm{cg}}$ such that $\mathrm{Ext}\left( \mathbb{T},X\right) =0$. By
Proposition \ref{Proposition:characterize-LH(LCPABcg)}, we can write%
\begin{equation*}
X=\frac{C}{F\oplus V^{\prime }}\oplus F^{\prime }\oplus V
\end{equation*}%
for vector groups $V,V^{\prime }$ and finite-rank torsion-free groups $%
F^{\prime },F$. Since $\mathrm{Ext}\left( \mathbb{T},X\right) =0$ and $%
F^{\prime }$ and $C/\left( F\oplus V^{\prime }\right) $ are direct summands
of $X$, we must have $\mathrm{Ext}\left( \mathbb{T},F^{\prime }\right) =0$
and $\mathrm{Ext}(\mathbb{T},C/\left( F\oplus V^{\prime }\right) )=0$. By
Proposition \ref{Proposition:injectives} this implies that $F^{\prime }=0$.
Furthermore, by Lemma \ref{Lemma:FLCPAbcg-injective-1} we have that $%
C/\left( F\oplus V^{\prime }\right) $ is isomorphic to $T/\left( E\oplus
W\right) $ for some torus group $T$, vector group $W$, and finite-rank free
subgroup $E$ of $T$. Thus, we have%
\begin{equation*}
X=\frac{C}{F\oplus V^{\prime }}\oplus V\cong \frac{T}{E\oplus W}\oplus V%
\text{,}
\end{equation*}%
concluding the proof.
\end{proof}

\begin{lemma}
\label{Lemma:FLCPAbcg-injective-2}Suppose that $X\in \mathrm{LH}\left( 
\mathbf{LCPAb}\right) $. If $\mathrm{Ext}\left( \mathbb{T},X\right) =0$,
then $X\cong \left( V\oplus T\oplus F\right) /(\Lambda \oplus W)$ where $T$
is a torus, $\Lambda $ is countable, $F$ is finite, and $V,W$ are vector
groups.
\end{lemma}

\begin{proof}
Represent $X$ as $G/H$ where $G$ is a locally compact Polish abelian group
and $H\subseteq G$ is a locally compact Polish subgroup. By the structure
theorem for locally compact abelian groups, $H=W\oplus \Lambda $ where $W$
is a vector group and $\Lambda $ has a compact open subgroup. Without loss
of generality, we may assume that $\Lambda $ is countable.

Similarly, $G\cong V\oplus G_{0}$ where $V$ is a vector group and $G_{0}$
has a compact open subgroup $C$. Define $G^{\prime }:=V\oplus C$ and $%
H^{\prime }:=H\cap G^{\prime }$. Consider the short exact sequence%
\begin{equation*}
0\rightarrow G^{\prime }/H^{\prime }\rightarrow G/H\rightarrow \frac{%
G/G^{\prime }}{H/H^{\prime }}\rightarrow 0\text{.}
\end{equation*}%
By the induced long exact sequence%
\begin{equation*}
\Gamma :=\frac{G/G^{\prime }}{H/H^{\prime }}
\end{equation*}%
is a countable group for which $\mathrm{Ext}\left( \mathbb{T},\Gamma \right)
=0$. This implies that $\Gamma =0$ by Proposition \ref%
{Proposition:injectives}. Furthermore, we have that $H^{\prime }$ is a
closed subgroup of $H$. Hence it is also of the form $W^{\prime }\oplus
\Lambda ^{\prime }$ where $W^{\prime }$ is a vector group and $\Lambda
^{\prime }$ is countable. Thus, after replacing $G$ with $G^{\prime }$ and $%
H $ with $H^{\prime }$ we may assume without loss of generality that $%
G=V\oplus C$ and $H=W\oplus \Lambda $ where $V,W$ are vector groups, $C$ is
compact, and $\Lambda $ is countable.

The short exact sequence $\mathbb{Z}\rightarrow \mathbb{R}\rightarrow 
\mathbb{T}$ induces an exact sequence%
\begin{equation*}
0=\mathrm{Hom}\left( \mathbb{R},\Lambda \right) \rightarrow \mathrm{Hom}%
\left( \mathbb{Z},\Lambda \right) \rightarrow \mathrm{Ext}\left( \mathbb{T}%
,\Lambda \right) \rightarrow \mathrm{Ext}\left( \mathbb{R},\Lambda \right) =0%
\text{.}
\end{equation*}%
Thus, we have that $\mathrm{Ext}\left( \mathbb{T},\Lambda \right) \cong 
\mathrm{Hom}\left( \mathbb{Z},\Lambda \right) \cong \Lambda $ and%
\begin{equation*}
\mathrm{Ext}\left( \mathbb{T},H\right) \cong \mathrm{Ext}\left( \mathbb{T}%
,W\right) \oplus \mathrm{Ext}\left( \mathbb{T},\Lambda \right) \cong \mathrm{%
Ext}\left( \mathbb{T},\Lambda \right) \cong \Lambda \text{.}
\end{equation*}

The inclusion $H\rightarrow G$ now induces an exact sequence%
\begin{equation*}
\Lambda \cong \mathrm{Ext}\left( \mathbb{T},H\right) \rightarrow \mathrm{Ext}%
\left( \mathbb{T},V\oplus C\right) \cong \mathrm{Ext}\left( \mathbb{T}%
,C\right) \rightarrow \mathrm{Ext}\left( \mathbb{T},X\right) =0
\end{equation*}%
Thus $\mathrm{Ext}\left( \mathbb{T},C\right) $ is countable. Hence, by
Corollary \ref{Corollary:finite-ranks-countable-Ext}, $C\cong T\oplus F$
where $T$ is a torus and $F$ is finite. The conclusion follows.
\end{proof}

\begin{proposition}
\label{Proposition:structure-injective}Suppose that $X\in \mathrm{LH}\left( 
\mathbf{LCPAb}\right) $. The following assertions are equivalent:

\begin{enumerate}
\item $X$ is injective for $\mathbf{LCPAb}$;

\item $\mathrm{Ext}\left( \mathbb{T},X\right) =0$;

\item $X\cong G/N$ where $G$ is an injective object of $\mathbf{LCPAb}$, and 
$N$ is the topological direct sum of a countable group and a vector group.
\end{enumerate}
\end{proposition}

\begin{proof}
(1)$\Rightarrow $(2) is obvious.

(3)$\Rightarrow $(1) Let $C\in \mathbf{LCPAb}$. The exact sequence $%
0\rightarrow N\rightarrow G\rightarrow G/N\rightarrow 0$ induces an exact
sequence%
\begin{equation*}
0=\mathrm{Ext}\left( C,G\right) \rightarrow \mathrm{Ext}\left( C,G/N\right)
\rightarrow \mathrm{Ext}^{2}\left( C,N\right) =0.
\end{equation*}%
Here the first term vanishes because $G$ is injective for $\mathbf{LCPAb}$,
and the last term vanishes by the vanishing of higher $\mathrm{Ext}$ for
locally compact Polish groups established in Subsection \ref%
{Subsection:derived-lc}. Hence $\mathrm{Ext}\left( C,G/N\right) =0$ for
every $C\in \mathbf{LCPAb}$, proving (1).

(2)$\Rightarrow $(3): By Lemma \ref{Lemma:FLCPAbcg-injective-2} we have that 
$X\cong \left( G\oplus F\right) /M$ where $G$ is an injective object of $%
\mathbf{LCPAb}$, $F$ is finite, and $M$ is the topological direct sum of a
vector group and a countable group. We have an exact sequence%
\begin{equation*}
0\rightarrow G/N\rightarrow X\rightarrow \frac{F}{\pi _{F}\left( M\right) }%
\rightarrow 0,
\end{equation*}%
where $N:=\mathrm{Ker}\left( \pi _{F}|_{M}\right) $ and $\pi _{F}:G\oplus
F\rightarrow F$ is the canonical projection. The long exact sequence
associated with $0\rightarrow N\rightarrow G\rightarrow G/N\rightarrow 0$
places $\mathrm{Ext}^{2}\left( \mathbb{T},G/N\right) $ between $\mathrm{Ext}%
^{2}\left( \mathbb{T},G\right) =0$ and $\mathrm{Ext}^{3}\left( \mathbb{T}%
,N\right) =0$, by the higher-$\mathrm{Ext}$ vanishing established in
Subsection \ref{Subsection:derived-lc}. Therefore, the induced exact sequence%
\begin{equation*}
\mathrm{Ext}\left( \mathbb{T},X\right) \rightarrow \mathrm{Ext}\left( 
\mathbb{T},F/\pi _{F}\left( M\right) \right) \rightarrow 0
\end{equation*}%
shows that%
\begin{equation*}
\mathrm{Ext}\left( \mathbb{T},F/\pi _{F}\left( M\right) \right) \cong 0\text{%
.}
\end{equation*}%
Since 
\begin{equation*}
\mathrm{Ext}\left( \mathbb{T},F/\pi _{F}\left( M\right) \right) \cong F/\pi
_{F}\left( M\right)
\end{equation*}%
we infer that $F/\pi _{F}\left( M\right) =0$. Thus, 
\begin{equation*}
X\cong G/N\text{.}
\end{equation*}%
Since $M$ is the topological sum of a vector group and a countable group,
the same holds for $N$, concluding the proof.
\end{proof}

\begin{lemma}
\label{Lemma:2-zero-Ext}If $X\in \mathbf{LCPAb}$ satisfies $\mathrm{Ext}%
\left( \mathbb{T},X\right) =0$ and $\mathrm{Ext}(\mathbb{Q}^{\vee }/\mathbb{Z%
},X)=0$, then $X=0$.
\end{lemma}

\begin{proof}
{By} Proposition \ref{Proposition:injectives} we must have that $X\cong 
\mathbb{R}^{n}\oplus \mathbb{T}^{k}$ for some $k,n\in \omega $. We have
short exact sequences%
\begin{equation*}
\mathrm{Hom}(\mathbb{Q}^{\vee },\mathbb{T})\cong \mathrm{Hom}(\mathbb{Z},%
\mathbb{Q})\cong \mathbb{Q}\rightarrow \mathrm{Hom}\left( \mathbb{Z},\mathbb{%
T}\right) \cong \mathbb{T}\rightarrow \mathrm{Ext}(\mathbb{Q}^{\vee }/%
\mathbb{Z},\mathbb{T})\rightarrow 0
\end{equation*}%
and%
\begin{equation*}
0=\mathrm{Hom}(\mathbb{Q}^{\vee },\mathbb{R})\rightarrow \mathrm{Hom}\left( 
\mathbb{Z},\mathbb{R}\right) \cong \mathbb{R}\rightarrow \mathrm{Ext}(%
\mathbb{Q}^{\vee }/\mathbb{Z},\mathbb{R})\rightarrow 0
\end{equation*}%
These show that $\mathrm{Ext}(\mathbb{Q}^{\vee }/\mathbb{Z},\mathbb{R})\neq
0 $ and $\mathrm{Ext}(\mathbb{Q}^{\vee }/\mathbb{Z},\mathbb{T})\neq 0$.
Therefore, we must have that $n=k=0$ and $X=0$.
\end{proof}

\begin{lemma}
\label{Lemma:vanishing-Ext2}Suppose that $X/Y$ and $G/H$ are groups with a
locally compact Polish abelian cover. If $\mathrm{Ext}\left( Y,G\right) =0$,
then $\mathrm{Ext}^{2}\left( X/Y,G/H\right) =0$.
\end{lemma}

\begin{proof}
We have an exact sequence%
\begin{equation*}
\mathrm{Ext}\left( Y,G/H\right) \rightarrow \mathrm{Ext}^{2}\left(
X/Y,G/H\right) \rightarrow \mathrm{Ext}^{2}\left( X,G/H\right)
\end{equation*}%
as well as exact sequences%
\begin{equation*}
0=\mathrm{Ext}^{2}\left( X,G\right) \rightarrow \mathrm{Ext}^{2}\left(
X,G/H\right) \rightarrow 0
\end{equation*}%
(since $X$ and $G$ are in $\mathbf{LCPAb}$) and, by assumption 
\begin{equation*}
0=\mathrm{Ext}\left( Y,G\right) \rightarrow \mathrm{Ext}\left( Y,G/H\right)
\rightarrow \mathrm{Ext}^{2}\left( Y,H\right) =0
\end{equation*}%
These imply the vanishing of the outer terms of the first sequence,
concluding the proof.
\end{proof}

Let $\alpha $ be an irrational real number. Consider the subgroup $H_{\alpha
}=\mathbb{Z}+\alpha \mathbb{Z}$ of $\mathbb{R}$. Then $H_{\alpha }\cong 
\mathbb{Z}^{2}$ is a locally compact Polish group with the discrete
topology. Set%
\begin{equation*}
A_{\alpha }:=\mathbb{R}/H_{\alpha }\text{.}
\end{equation*}%
Notice that both $\mathbb{Z}^{2}$ and $\mathbb{R}$ are projective in $%
\mathbf{LCPAb}$ and hence in $\mathrm{LH}\left( \mathbf{LCPAb}\right) $.
Thus,%
\begin{equation*}
0\rightarrow H_{\alpha }\rightarrow \mathbb{R}\rightarrow A_{\alpha
}\rightarrow 0
\end{equation*}%
is a projective resolution for $A_{\alpha }$.\ 

\begin{lemma}
\label{Lemma:zero}Suppose that $G=V\oplus T$ where $V\cong \mathbb{R}^{k}$
and $T\cong \mathbb{T}^{\sigma }$ for some $k\in \mathbb{N}$ and $\sigma
\leq \omega $. Let also $N\subseteq G$ be a locally compact Polish subgroup
such that $N=\Lambda \oplus W$ where $\Lambda $ is a countable group and $W$
is a vector group. Fix also an irrational real number $\alpha $. If $\mathrm{%
Ext}\left( A_{\alpha },G/N\right) =0$, then $G/N=0$.
\end{lemma}

\begin{proof}
We have%
\begin{equation*}
\mathrm{Ext}\left( A_{\alpha },G/N\right) \cong \frac{\mathrm{Hom}\left(
H_{\alpha },G/N\right) }{\mathrm{Hom}\left( \mathbb{R}|H_{\alpha
},G/N\right) }
\end{equation*}%
where $\mathrm{Hom}\left( \mathbb{R}|H_{\alpha },G/N\right) $ denotes the
group of Borel-definable homomorphisms $H_{\alpha }\rightarrow G/N$ that
extend to Borel-definable homomorphisms $\mathbb{R}\rightarrow G/N$. Since $%
H_{\alpha }\cong \mathbb{Z}^{2}$, we have that%
\begin{equation*}
\mathrm{Hom}\left( H_{\alpha },G/N\right) \cong G/N\oplus G/N\text{.}
\end{equation*}%
Via this isomorphism, $\mathrm{Hom}\left( \mathbb{R}|H_{\alpha },G/N\right) $
corresponds to the subgroup%
\begin{equation*}
\left\{ \left( \varphi \left( 1\right) ,\varphi \left( \alpha \right)
\right) :\varphi \in \mathrm{Hom}\left( \mathbb{R},G/N\right) \right\} \text{%
,}
\end{equation*}%
where $\mathrm{Hom}\left( \mathbb{R},G/N\right) $ is the hom-set in the
category $\mathrm{LH}\left( \mathbf{LCPAb}\right) $.

By projectivity of $\mathbb{R}$ in $\mathrm{LH}\left( \mathbf{LCPAb}\right) $%
, every Borel-definable homomorphism $\mathbb{R}\rightarrow G/N$ lifts to a
continuous group homomorphism $\mathbb{R}\rightarrow G$. Thus, denoting by $%
\mathrm{Hom}\left( \mathbb{R},G\right) $ the hom-set in $\mathbf{LCPAb}$,
i.e., the group of continuous group homomorphisms $\mathbb{R}\rightarrow G$,
one obtains that%
\begin{eqnarray*}
\mathrm{Ext}\left( A_{\alpha },G/N\right) &\cong &\frac{G/N\oplus G/N}{%
\left\{ \left( \varphi \left( 1\right) ,\varphi \left( \alpha \right)
\right) :\varphi \in \mathrm{Hom}\left( \mathbb{R},G/N\right) \right\} } \\
&\cong &\frac{G\oplus G}{N\oplus N+\left\{ \left( \gamma \left( 1\right)
,\gamma \left( \alpha \right) \right) :\gamma \in \mathrm{Hom}\left( \mathbb{%
R},G\right) \right\} }\text{.}
\end{eqnarray*}%
The assumption that $\mathrm{Ext}\left( A_{\alpha },G/N\right) =0$ yields%
\begin{equation*}
G\oplus G=N\oplus N+L
\end{equation*}%
where%
\begin{equation*}
L:=\left\{ \left( \gamma \left( 1\right) ,\gamma \left( \alpha \right)
\right) :\gamma \in \mathrm{Hom}\left( \mathbb{R},G\right) \right\} \text{.}
\end{equation*}%
Fix now $d\leq \sigma $ such that $d<\omega $, and set 
\begin{equation*}
G_{d}:=\mathbb{R}^{k}\oplus \mathbb{T}^{d}\subseteq \mathbb{R}^{k}\oplus 
\mathbb{T}^{\sigma }=G\text{.}
\end{equation*}%
Let also $\hat{G}_{d}\cong \mathbb{R}^{k}\oplus \mathbb{R}^{d}\cong \mathbb{R%
}^{k+d}$ be a vector group, and $p_{d}:\hat{G}_{d}\rightarrow G_{d}$ be the
canonical universal covering map. Notice that every continuous homomorphism $%
\gamma :\mathbb{R}\rightarrow G_{d}$ lifts to a continuous homomorphism $%
\hat{\gamma}:\mathbb{R}\rightarrow \hat{G}_{d}$ by projectivity of $\mathbb{R%
}$. Let also $\pi _{d}:G\rightarrow G_{d}$ be the canonical projection.
Define:

\begin{itemize}
\item $K_{d}:=\mathrm{Ker}\left( \pi _{d}|_{N}\right) $;

\item $N_{d}:=N/K_{d}$;

\item $q_{d}:N\rightarrow N_{d}$, the canonical quotient map.
\end{itemize}

We can identify $N_{d}$ with a locally compact Polish subgroup of $G_{d}$
via its image under the continuous homomorphism induced by $\pi _{d}|_{N}$.

Let $\varpi _{V}:G\rightarrow V$ be the canonical projection. Then $\mathrm{%
Ker}\left( \varpi _{V}\right) =T$. Thus, $\mathrm{Ker}\left( \varpi
_{V}|_{W}\right) =T\cap W=0$ since $W$ has no nontrivial compact subgroups.
It follows that $K_{d}\cap W=\left\{ 0\right\} $ and $q_{d}|_{W}$ is
injective. Thus, $q_{d}|_{W}$ establishes a topological isomorphism between $%
W$ and its image $W_{d}\subseteq N_{d}$. This shows that $\mathrm{\dim }%
\left( W_{d}\right) =\mathrm{\dim }\left( W\right) $, and $W_{d}$ is an open
subgroup of $N_{d}$ by openness of the quotient map $q_{d}:N\rightarrow
N_{d} $. By injectivity of vector groups, we have a topological direct sum $%
N_{d}=W_{d}\oplus \Lambda _{d}$ for a countable discrete group $\Lambda _{d}$%
.

Finally, let%
\begin{equation*}
L_{d}:=\left\{ \left( \gamma \left( 1\right) ,\gamma \left( \alpha \right)
\right) :\gamma \in \mathrm{Hom}\left( \mathbb{R},G_{d}\right) \right\}
\end{equation*}%
and%
\begin{equation*}
\hat{L}_{d}:=\left\{ \left( \gamma \left( 1\right) ,\gamma \left( \alpha
\right) \right) :\gamma \in \mathrm{Hom}(\mathbb{R},\hat{G}_{d})\right\} 
\text{.}
\end{equation*}%
Notice that%
\begin{equation*}
L_{d}=\left( \pi _{d}\oplus \pi _{d}\right) \left( L\right) =\left(
p_{d}\oplus p_{d}\right) (\hat{L}_{d})\text{.}
\end{equation*}%
Furthermore,%
\begin{equation*}
\hat{L}_{d}=\left\{ \left( x,\alpha x\right) :x\in \hat{G}_{d}\right\}
\end{equation*}%
and hence%
\begin{equation*}
\dim (\hat{L}_{d})=\dim (\hat{G}_{d})=k+d\text{.}
\end{equation*}%
Then we have%
\begin{equation*}
G_{d}\oplus G_{d}=N_{d}\oplus N_{d}+L_{d}\text{.}
\end{equation*}

The inclusion 
\begin{equation*}
W_{d}\rightarrow N_{d}\rightarrow G_{d}
\end{equation*}%
lifts by projectivity of $W_{d}$ to a continuous homomorphism $%
s_{d}:W_{d}\rightarrow \hat{G}_{d}$. Set%
\begin{equation*}
S_{d}:=s_{d}\left( W_{d}\right) \text{.}
\end{equation*}%
Define%
\begin{equation*}
r_{d}:=\mathrm{\dim }\left( S_{d}\right) \text{.}
\end{equation*}%
Notice that%
\begin{equation*}
\dim \left( S_{d}\right) \leq \dim \left( W_{d}\right) =\dim \left( W\right) 
\text{.}
\end{equation*}%
Since%
\begin{equation*}
\left( S_{d}\oplus S_{d}\right) \cap \hat{L}_{d}=\left\{ \left( x,\alpha
x\right) :x\in S_{d}\right\}
\end{equation*}%
we have%
\begin{equation*}
\dim ((S_{d}\oplus S_{d})\cap \hat{L}_{d})=r_{d}\text{.}
\end{equation*}%
Thus, we have that%
\begin{equation*}
\hat{G}_{d}\oplus \hat{G}_{d}=S_{d}\oplus S_{d}+\hat{L}_{d}+\left(
p_{d}\oplus p_{d}\right) ^{-1}\left( \Lambda _{d}\oplus \Lambda _{d}\right) 
\text{.}
\end{equation*}%
Since $\Lambda _{d}$ and $\mathrm{Ker}\left( p_{d}\right) $ are countable, 
\begin{equation*}
\left( p_{d}\oplus p_{d}\right) ^{-1}\left( \Lambda _{d}\oplus \Lambda
_{d}\right)
\end{equation*}%
is also countable, and $S_{d}\oplus S_{d}+\hat{L}_{d}$ is an \emph{open}
subgroup $\hat{G}_{d}\oplus \hat{G}_{d}$. Since $\hat{G}_{d}\oplus \hat{G}%
_{d}$ is also connected, this yields%
\begin{equation*}
\hat{G}_{d}\oplus \hat{G}_{d}=S_{d}\oplus S_{d}+\hat{L}_{d}\text{.}
\end{equation*}%
Therefore, we have%
\begin{eqnarray*}
2\left( k+d\right) &=&\dim (\hat{G}_{d}\oplus \hat{G}_{d}) \\
&=&\mathrm{\dim }(S_{d}\oplus S_{d}+\hat{L}_{d}) \\
&=&\dim \left( S_{d}\oplus S_{d}\right) +\dim (\hat{L}_{d})-\dim
((S_{d}\oplus S_{d})\cap \hat{L}_{d}) \\
&=&2r_{d}+\left( k+d\right) -r_{d}\text{.}
\end{eqnarray*}%
Hence,%
\begin{equation*}
k\leq k+d=r_{d}\leq \dim \left( W\right) \text{.}
\end{equation*}%
Since this holds for every $d<\sigma $, and $W$ is finite-dimensional, we
must have that $\sigma <\omega $.

Set now $d:=\sigma <\omega $. Notice that in this case, we have that $%
N_{d}=N $ and $W_{d}=W$. We also set:

\begin{itemize}
\item $\hat{G}:=\hat{G}_{d}$;

\item $\hat{L}:=\hat{L}_{d}$;

\item $p=p_{d}:\hat{G}\rightarrow G$;

\item Let $s:W\rightarrow \hat{G}$ be a lift of the inclusion $W\rightarrow
G $;

\item $S:=s\left( W\right) \subseteq \hat{G}$;

\item $r:=\dim \left( S\right) \leq \dim \left( W\right) $.
\end{itemize}

Then by the above we have%
\begin{equation*}
\hat{G}\oplus \hat{G}=S\oplus S+\hat{L}
\end{equation*}%
and $k+d=r$. Since $S$ is a vector subspace of $\hat{G}$ and $\dim \left(
S\right) =\dim (\hat{G})$, this implies $S=\hat{G}$.

If $\hat{x}\in \mathbb{Z}^{d}$ is nonzero and $\left( 0,\hat{x}\right) \in 
\mathbb{R}^{k}\oplus \mathbb{R}^{d}=\hat{G}$ then we have that there exists $%
x\in W$ such that $s(x)=\left( 0,\hat{x}\right) $ and hence $x=p\left( 0,%
\hat{x}\right) =0$ and $\left( 0,\hat{x}\right) =s(x)=0$ and $\hat{x}=0$.
This forces $d=0$. Thus, we have that $\hat{G}=G=V=\mathbb{R}^{k}$. Again
the surjectivity of $s:W\rightarrow \hat{G}$ yields that $s$ is an
isomorphism, and hence $W=G$ and $G/N=G/W$ is trivial.
\end{proof}

\begin{theorem}
\label{Theorem:no-injectives-LCPAb}Suppose that $X\in \mathrm{LH}\left( 
\mathbf{LCPAb}\right) $. Let also $\alpha $ be an irrational real number. If 
$\mathrm{Ext}\left( \mathbb{T},X\right) =0$ and $\mathrm{Ext}\left(
A_{\alpha },X\right) =0$, then $X=0$.
\end{theorem}

\begin{proof}
By Proposition \ref{Proposition:structure-injective} we have $X\cong G/N$,
where 
\begin{equation*}
G=V\oplus T\text{\quad and\quad }N=\Lambda \oplus W\text{.}
\end{equation*}%
Here, $V,W$ are finite-dimensional vector groups, $T$ is a torus group, and $%
\Lambda $ is a countable group. Thus, by Lemma \ref{Lemma:zero}, $X=0$.
\end{proof}

\begin{corollary}
If $X\in \mathrm{LH}\left( \mathbf{LCPAb}\right) $ is injective, then $X=0$.
\end{corollary}

\subsection{Injectives in the heart of topological torsion groups}

Recall the following notations for subcategories of $\mathbf{LCPAb}$: $%
\mathbf{LCPAb}\left( p\right) $ denotes the category of topological $p$%
-groups, $\mathbf{FLCPAb}\left( p\right) $ denotes the category of finite $p$%
-rank topological $p$-groups, $\mathbf{TorFLCPAb}$ denotes the category of
topological torsion groups of finite ranks, and $\mathbf{TDLCPAb}$ denotes
the category of totally disconnected locally compact Polish abelian groups.

\begin{lemma}
\label{Lemma:almost-TDLCPAb-injective}Let $G$ be a locally compact Polish
abelian group. Then the following assertions are equivalent:

\begin{enumerate}
\item $\mathrm{Ext}\left( X,G\right) $ is countable for every countable
product $X$ of finite cyclic groups;

\item $G\cong V\oplus T\oplus A$ where $A$ is a countable group, $V$ is a
vector group, and $T$ is a torus group.
\end{enumerate}
\end{lemma}

\begin{proof}
(2)$\Rightarrow $(1) follows from Lemma \ref{Lemma:Ext(compact,countable)}
(as well as from Lemma \ref{Lemma:Ext(torsion,countable)}), together with
injectivity of $V$ and $T$ in $\mathbf{LCPAb}$.

(1)$\Rightarrow $(2) We can write $G=V\oplus H$ where $V$ is a vector group,
and the connected component $C$ of $H$ is compact. Set $A:=H/C$. Let $X$ be
a countable product of finite cyclic groups. The short exact sequence%
\begin{equation*}
0\rightarrow C\rightarrow H\rightarrow A\rightarrow 0
\end{equation*}%
induces an exact sequence 
\begin{equation*}
\mathrm{Ext}\left( X,H\right) \rightarrow \mathrm{Ext}\left( X,A\right)
\rightarrow \mathrm{Ext}^{2}\left( X,C\right) =0\text{.}
\end{equation*}%
Since $G=V\oplus H$ and $V$ is injective, the hypothesis implies that $%
\mathrm{Ext}\left(X,H\right)$ is countable. The displayed exact sequence
therefore implies that $\mathrm{Ext}\left(X,A\right)$ is countable. Let now $%
K\subseteq A$ be a compact open subgroup of $A$. Then the short exact
sequence%
\begin{equation*}
0\rightarrow K\rightarrow A\rightarrow A/K\rightarrow 0
\end{equation*}%
induces an exact sequence%
\begin{equation*}
\mathrm{Hom}\left( X,A/K\right) \rightarrow \mathrm{Ext}\left( X,K\right)
\rightarrow \mathrm{Ext}\left( X,A\right) \text{.}
\end{equation*}%
Since $X$ is compact and $A/K$ is countable discrete, $\mathrm{Hom}\left(
X,A/K\right) $ is countable. Since $\mathrm{Ext}\left( X,A\right) $ is also
countable, we infer that $\mathrm{Ext}\left( X,K\right) $ is also countable.
Since $X$ is an arbitrary product of finite cyclic groups, by Corollary \ref%
{Corollary:product-cyclic-countable-Ext},%
\begin{equation*}
K\cong T_{K}\oplus F
\end{equation*}%
where $T_{K}$ is a torus group and $F$ is a finite group. Since $K$ is
totally disconnected, we must have that $T_{K}=0$ and hence $K=F$ is finite.
Since $A/K$ is countable, this implies that $A$ is countable. Let again $X$
be an arbitrary product of finite cyclic groups. The exact sequence%
\begin{equation*}
\mathrm{Hom}\left( X,A\right) \rightarrow \mathrm{Ext}\left( X,C\right)
\rightarrow \mathrm{Ext}\left( X,H\right)
\end{equation*}%
implies that $\mathrm{Ext}\left( X,C\right) $ is countable. Applying again
Corollary \ref{Corollary:product-cyclic-countable-Ext}, since $C$ is
connected, we infer that $C$ is a torus group. By injectivity of torus
groups, this implies that $H=C\oplus A$ and hence $G=V\oplus C\oplus A$,
concluding the proof.
\end{proof}

\begin{proposition}
\label{Proposition:injectives-TDLC}Let $G$ be a locally compact Polish
abelian group. The following assertions are equivalent:

\begin{enumerate}
\item $G$ is injective for $\mathbf{TDLCPAb}$;

\item $G\cong V\oplus T\oplus D$ where $V$ is a vector group, $T$ is a
torus, and $D$ is countable and divisible.
\end{enumerate}
\end{proposition}

\begin{proof}
(2)$\Rightarrow $(1) If $X$ is totally disconnected, then $X_{\mathbb{S}%
^{1}}=0$. Thus, $D$ is injective for $\mathbf{TDLCPAb}$ by \cite[Proposition
4.15(vi)]{hoffmann_homological_2007}. We also have that $V$ and $T$ are
injective for $\mathbf{LCPAb}$.

(1)$\Rightarrow $(2) By Lemma \ref{Lemma:almost-TDLCPAb-injective}, we have
that $G\cong V\oplus T\oplus A$ where $V$ is a vector group, $T$ is a torus
group, and $A$ is a countable group. Since $G$ is injective for $\mathbf{%
TDLCPAb}$, the same holds for $A$. In particular, $A$ is injective in the
category of countable abelian groups. Hence, $A$ is divisible.
\end{proof}

\begin{proposition}
\label{Proposition:almost-LCPAb-injective-LH(LCPAb)}Suppose that $G/H\in 
\mathrm{LH}\left( \mathbf{LCPAb}\right) $. Then the following assertions are
equivalent:

\begin{enumerate}
\item $\mathrm{Ext}(X,G/H)$ is countable for every countable product $X$ of
finite cyclic groups;

\item $G/H\cong L/\left( N\oplus W\right) $ where $L\cong V\oplus T\oplus A$%
, for some vector groups $V,W$, torus group $T$, and countable discrete
groups $A$ and $N$.
\end{enumerate}
\end{proposition}

\begin{proof}
After quotienting both $G$ and $H$ by a compact subgroup of $H$, we can
assume that $H=N\oplus W$ where $W$ is a vector group and $N$ is countable.

(1)$\Rightarrow $(2) Let $X$ be a countable product of finite cyclic groups.
Consider the exact sequence%
\begin{equation*}
\mathrm{Ext}\left( X,H\right) \rightarrow \mathrm{Ext}\left( X,G\right)
\rightarrow \mathrm{Ext}\left( X,G/H\right) \text{.}
\end{equation*}%
Then $\mathrm{Ext}\left( X,G/H\right) $ is countable by hypothesis, and $%
\mathrm{Ext}\left( X,H\right) $ is countable by Lemma \ref%
{Lemma:almost-TDLCPAb-injective}. Thus, $\mathrm{Ext}\left( X,G\right) $ is
countable. Since this holds for every countable product $X$ of finite cyclic
groups, by Lemma \ref{Lemma:almost-TDLCPAb-injective} we have that $G\cong
V\oplus T\oplus A$, for some vector groups $V,W$, torus group $T$, and
countable discrete groups $A$ and $N$.

(2)$\Rightarrow $(1) Suppose that $G=V\oplus T\oplus A$ and $H=N\oplus W$.
Let $X$ be a countable product of finite cyclic groups. Consider the exact
sequence%
\begin{equation*}
\mathrm{Ext}\left( X,G\right) \rightarrow \mathrm{Ext}\left( X,G/H\right)
\rightarrow \mathrm{Ext}^{2}\left( X,H\right) =0\text{.}
\end{equation*}%
Then by Lemma \ref{Lemma:almost-TDLCPAb-injective}, $\mathrm{Ext}\left(
X,G\right) $ is countable. This implies that $\mathrm{Ext}\left(
X,G/H\right) $ is also countable. We have an exact sequence%
\begin{equation*}
\mathrm{Ext}\left( X,H\right) \cong \mathrm{Ext}\left( X,N\right)
\rightarrow \mathrm{Ext}\left( X,G\right) \rightarrow \mathrm{Ext}\left(
X,G/H\right)
\end{equation*}%
where $\mathrm{Ext}\left( X,N\right) $ is countable by either Lemma \ref%
{Lemma:Ext(torsion,countable)} or Lemma \ref{Lemma:Ext(compact,countable)}.
Thus $\mathrm{Ext}\left( X,G\right) $ is countable.
\end{proof}

\begin{theorem}
\label{Theorem:LCPAb-injective-LH(LCPAb)}Suppose that $G/H\in \mathrm{LH}%
\left( \mathbf{LCPAb}\right) $. The following assertions are equivalent:

\begin{enumerate}
\item $G/H$ is injective for $\mathbf{TDLCPAb}$;

\item $G/H\cong L/\left( N\oplus W\right) $ where $L$ is a Polish abelian
group that is injective for $\mathbf{TDLCPAb}$, $W$ is a vector group, and $%
N $ is countable;

\item $G/H\cong D\oplus J/\left( N^{\prime }\oplus W^{\prime }\right) $
where $D$ is countable and divisible, $J$ is injective for $\mathbf{LCPAb}$,
and $N^{\prime }$ is countable and $W^{\prime }$ is a vector group.
\end{enumerate}
\end{theorem}

\begin{proof}
We can write $H=W\oplus N$ where $W$ is a vector group and $N$ has a compact
open subgroup. By quotienting by that subgroup, we may without loss of
generality assume that $N$ is countable. Furthermore, since $N$ embeds into
a countable divisible group, by considering a suitable pushout we can assume
that $N$ is \emph{divisible}; see Lemma \ref{Lemma:presentation}.

(1)$\Rightarrow $(2): Since $N$ and $W$ are injective for $\mathbf{TDLCPAb}$%
, the long exact sequence relating $\mathrm{Hom}$ and $\mathrm{Ext}$ implies
that $L$ is also injective for $\mathbf{TDLCPAb}$.

(2)$\Rightarrow $(3): By Propositions \ref{Proposition:injectives-TDLC} and %
\ref{Proposition:injectives}, $L=D\oplus J$ where $D$ is countable and
divisible, and $J$ is injective for $\mathbf{LCPAb}$. Let $\pi _{D}:D\oplus
J\rightarrow D$ be the canonical projection, and observe that $\pi
_{D}\left( W\right) =0$. Thus, we have a short exact sequence%
\begin{equation*}
J/\left( J\cap (N\oplus W)\right) \rightarrow L/\left( N\oplus W\right)
\rightarrow D/\pi _{D}\left( N\right)
\end{equation*}%
Since $J$ is injective for $\mathbf{LCPAb}$, the same holds for $J/\left(
J\cap (N\oplus W)\right) $, and this exact sequence splits. Thus%
\begin{equation*}
L/\left( N\oplus W\right) \cong D/\pi _{D}\left( N\right) \oplus J/(J\cap
\left( N\oplus W)\right)
\end{equation*}%
where $D/\pi _{D}\left( N\right) $ is countable and divisible; clearly $%
J\cap \left( N\oplus W\right) \cong N^{\prime }\oplus W^{\prime }$ for some
countable group $N^{\prime }$ and vector group $W^{\prime }$ as well.

(3)$\Rightarrow $(1): We need to prove that $D\oplus J/\left( N^{\prime
}\oplus W^{\prime }\right) $ is injective for $\mathbf{TDLCPAb}$. By
Proposition \ref{Proposition:injectives-TDLC} we have that $D$ is injective
for $\mathbf{TDLCPAb}$. Furthermore, if $X$ is a totally disconnected
locally compact Polish abelian group, then the exact sequence%
\begin{equation*}
0=\mathrm{Ext}\left( X,J\right) \rightarrow \mathrm{Ext}\left( X,J/\left(
N^{\prime }\oplus W^{\prime }\right) \right) \rightarrow \mathrm{Ext}%
^{2}\left( X,N^{\prime }\oplus W^{\prime }\right) =0
\end{equation*}%
shows that $\mathrm{Ext}\left( X,J/\left( N^{\prime }\oplus W^{\prime
}\right) \right) =0$; as $X$ is arbitrary, this implies (1). 
\end{proof}

\begin{lemma}
\label{Lemma:LH(TDLCPAb)}Consider $X=G/N\in \mathrm{LH}\left( \mathbf{LCPAb}%
\right) $. Suppose that 
\begin{equation*}
G=V\oplus T\oplus D\text{\quad and\quad }N=\Lambda \oplus W
\end{equation*}%
where $V,W$ are vector groups, $T\ $is a torus, and $D,\Lambda $ are
countable. If $G/N$ is isomorphic to an object of $\mathrm{LH}\left( \mathbf{%
TDLCPAb}\right) $, considered as a full subcategory of $\mathrm{LH}\left( 
\mathbf{LCPAb}\right) $, then $V\oplus T\subseteq N$, and $X$ is isomorphic
in $\mathrm{LH}\left( \mathbf{LCPAb}\right) $ to a countable discrete group.
\end{lemma}

\begin{proof}
By assumption there exist totally disconnected locally compact Polish groups 
$C$ and $N^{\prime }$ such that $C/N^{\prime }\in \mathrm{LH}\left( \mathbf{%
LCPAb}\right) $ is isomorphic in $\mathrm{LH}\left( \mathbf{LCPAb}\right) $
to $G/N$. Let $U$ be a finite-dimensional vector group. By projectivity of $%
U $ in $\mathrm{LH}\left( \mathbf{LCPAb}\right) $, the homomorphism%
\begin{equation*}
\mathrm{Hom}\left( U,C\right) \rightarrow \mathrm{Hom}\left( U,X\right)
\end{equation*}%
is surjective. Since $C$ is totally disconnected, $\mathrm{Hom}\left(
U,C\right) =0$. If $g\in V\oplus T$ then there exists a homomorphism $%
\varphi :\mathbb{R}\rightarrow V\oplus T$ that contains $g$ in its image.
Considering the induced homomorphism%
\begin{equation*}
\mathbb{R}\rightarrow V\oplus T\rightarrow G\rightarrow G/N\cong C/N^{\prime
}
\end{equation*}%
we conclude that it must be zero. This forces $g\in N\subseteq G$. Since
this holds for every $g\in V\oplus T$, it implies that $V\oplus T\subseteq N$%
. Thus, the projection $\pi _{D}:G\rightarrow D$ induces an isomorphism%
\begin{equation*}
X=G/N\cong \frac{D}{\pi _{D}\left( N\right) }
\end{equation*}%
between $X$ and a quotient of $D$.
\end{proof}

\begin{theorem}
\label{Theorem:LCPAb-injective-LH(TDLCPAb)}Suppose that $G/H\in \mathrm{LH}%
\left( \mathbf{TDLCPAb}\right) $. The following assertions are equivalent:

\begin{enumerate}
\item $G/H$ is injective for $\mathbf{TDLCPAb}$;

\item $G/H$ is injective for $\mathbf{TorLCPAb}_{\mathrm{c}}$;

\item $G/H\cong D$ where $D$ is a countable divisible group.
\end{enumerate}
\end{theorem}

\begin{proof}
(2)$\Rightarrow $(1) Put $E:=G/H$. Fix $n\in \mathbb{N}$. Consider $\mathbb{Z%
}/n\mathbb{Z}\in \mathbf{TorLCPAb}_{\mathrm{c}}$.\ Then it has a projective
resolution%
\begin{equation*}
0\rightarrow n\mathbb{Z}\rightarrow \mathbb{Z}\rightarrow \mathbb{Z}/n%
\mathbb{Z}\rightarrow 0\text{.}
\end{equation*}%
This produces%
\begin{equation*}
0=\mathrm{Ext}\left( \mathbb{Z}/n,E\right) \cong E/nE\text{.}
\end{equation*}%
This shows that $E$ is divisible. Let $D$ be a countable discrete group. Let%
\begin{equation*}
0\rightarrow R\rightarrow F\rightarrow D\rightarrow 0
\end{equation*}%
be a projective resolution in the category of countable discrete abelian
groups. Thus, $F$ is a free abelian group, and $R\subseteq F$ is a free
abelian subgroup. Since $E$ is divisible, every homomorphism $R\rightarrow E$
extends to a homomorphism $F\rightarrow E$. Thus, we obtain that%
\begin{equation*}
\mathrm{Hom}\left( F,E\right) \rightarrow \mathrm{Hom}\left( R,E\right)
\end{equation*}%
is onto. Since $F$ is a countable free abelian group, we also have 
\begin{equation*}
\mathrm{Ext}\left( F,E\right) =0\text{.}
\end{equation*}

Considering the exact sequence%
\begin{equation*}
\mathrm{Hom}\left( F,E\right) \rightarrow \mathrm{Hom}\left( R,E\right)
\rightarrow \mathrm{Ext}\left( D,E\right) \rightarrow \mathrm{Ext}\left(
F,E\right)
\end{equation*}%
we infer that 
\begin{equation*}
\mathrm{Ext}\left( D,E\right) =0\text{.}
\end{equation*}%
Suppose now that $X\in \mathbf{TDLCPAb}$.\ Let $K\subseteq X$ be a compact
open subgroup. Then $K\in \mathbf{TDLCPAb}_{\mathrm{c}}$, and $X/K$ is
countable. The exact sequence%
\begin{equation*}
0=\mathrm{Ext}\left( X/K,E\right) \rightarrow \mathrm{Ext}\left( X,E\right)
\rightarrow \mathrm{Ext}\left( K,E\right) =0
\end{equation*}%
shows that $\mathrm{Ext}\left( X,E\right) =0$.

(1)$\Rightarrow $(3) By Theorem \ref{Theorem:LCPAb-injective-LH(LCPAb)} we
have that $G/H\cong D\oplus \left( V\oplus T\right) /\left( \Lambda \oplus
W\right) $ where $D$ is countable and divisible, $V$ and $W$ are vector
groups, $T$ is a torus, and $\Lambda $ is countable. Since by hypothesis $%
G/H\in \mathrm{LH}\left( \mathbf{TDLCPAb}\right) $, we have by Lemma \ref%
{Lemma:LH(TDLCPAb)} that $\left( V\oplus T\right) /\left( \Lambda \oplus
W\right) =0$ and $G/H\cong D$.

(3)$\Rightarrow $(2) follows from Theorem \ref%
{Theorem:LCPAb-injective-LH(LCPAb)}.
\end{proof}

For a prime number $p$, set%
\begin{equation*}
\Xi _{p}:=\left( \mathbb{Z}/p\mathbb{Z}\right) ^{\omega }/\left( \mathbb{Z}/p%
\mathbb{Z}\right) ^{\left( \omega \right) }\text{.}
\end{equation*}

\begin{lemma}
\label{Lemma:Xi-p}We have that $\mathrm{Ext}(\Xi _{p},\mathbb{T})\cong 
\mathrm{Ext}(\Xi _{p},\mathbb{Z}\left( p^{\infty }\right) )\cong \Xi _{p}$.
\end{lemma}

\begin{proof}
Since $\mathbb{T}$ is injective in $\mathbf{LCPAb}$, we have that $\mathrm{%
Ext}\left( \left( \mathbb{Z}/p\mathbb{Z}\right) ^{\omega },\mathbb{T}\right)
=0$. Likewise, since countable divisible groups are injective in $\mathbf{%
TDLCPAb}$, we have that $\mathrm{Ext}\left( \left( \mathbb{Z}/p\mathbb{Z}%
\right) ^{\omega },\mathbb{Z}\left( p^{\infty }\right) \right) =0$.
Therefore,%
\begin{eqnarray*}
\mathrm{Ext}(\Xi _{p},\mathbb{T}) &\cong &\mathrm{coker}\left( \mathrm{Hom}%
\left( \left( \mathbb{Z}/p\mathbb{Z}\right) ^{\omega },\mathbb{T}\right)
\rightarrow \mathrm{Hom}(\left( \mathbb{Z}/p\mathbb{Z}\right) ^{\left(
\omega \right) },\mathbb{T})\right) \\
&\cong &\mathrm{coker}\left( \mathrm{Hom}\left( \left( \mathbb{Z}/p\mathbb{Z}%
\right) ^{\omega },\mathbb{Z}\left( p^{\infty }\right) \right) \rightarrow 
\mathrm{Hom}(\left( \mathbb{Z}/p\mathbb{Z}\right) ^{\left( \omega \right) },%
\mathbb{Z}\left( p^{\infty }\right) )\right) \\
&\cong &\mathrm{Ext}(\Xi _{p},\mathbb{Z}\left( p^{\infty }\right) )\text{.}
\end{eqnarray*}%
Furthermore,%
\begin{equation*}
\mathrm{coker}\left( \mathrm{Hom}\left( \left( \mathbb{Z}/p\mathbb{Z}\right)
^{\omega },\mathbb{Z}\left( p^{\infty }\right) \right) \rightarrow \mathrm{%
Hom}(\left( \mathbb{Z}/p\mathbb{Z}\right) ^{\left( \omega \right) },\mathbb{Z%
}\left( p^{\infty }\right) )\right) \cong \Xi _{p}
\end{equation*}
\end{proof}

\begin{lemma}
\label{Lemma:must-be-tf}Suppose that $D$ is a countable divisible group such
that $\mathrm{Ext}\left( \Xi _{p},D\right) $ is countable for every prime $p$%
. Then $D$ is torsion-free.
\end{lemma}

\begin{proof}
This is an immediate consequence of Lemma \ref{Lemma:Xi-p}, considering the
structure theorem for countable divisible groups \cite[Theorem 23.1]%
{fuchs_infinite_1970}.
\end{proof}

\begin{theorem}
\label{Theorem:no-injectives-TDLCPAb}Suppose that $G/H\in \mathrm{LH}\left( 
\mathbf{TDLCPAb}\right) $. The following assertions are equivalent:

\begin{enumerate}
\item $G/H$ is injective for $\mathrm{LH}(\mathbf{TDLCPAb})$;

\item $G/H=0$.
\end{enumerate}
\end{theorem}

\begin{proof}
(1)$\Rightarrow $(2) By Theorem \ref{Theorem:LCPAb-injective-LH(TDLCPAb)}, $%
G/H=D$ where $D$ is a countable divisible group. For every prime $p$, $%
\mathrm{Ext}\left( \Xi _{p},D\right) =0$. Thus, by Lemma \ref%
{Lemma:must-be-tf} $D$ is torsion-free, and hence $D$ is a direct sum of
copies of $\mathbb{Q}$.

Fix now a prime $p\in \mathbb{N}$ and consider the object%
\begin{equation*}
Y:=\mathbb{Q}_{p}/\mathbb{Q}\in \mathrm{LH}(\mathbf{TDLCPAb})\text{.}
\end{equation*}%
Considering the exact sequence%
\begin{equation*}
0\rightarrow \mathbb{Q}\rightarrow \mathbb{Q}_{p}\rightarrow Y\rightarrow 0
\end{equation*}%
gives an exact sequence%
\begin{equation*}
\mathrm{Hom}(\mathbb{Q}_{p},D)\rightarrow \mathrm{Hom}(\mathbb{Q}%
,D)\rightarrow \mathrm{Ext}\left( Y,D\right) \rightarrow \mathrm{Ext}(%
\mathbb{Q}_{p},D)\text{.}
\end{equation*}%
Then $\mathrm{Ext}(\mathbb{Q}_{p},D)=0$ by injectivity of $D$ for $\mathbf{%
TDLCPAb}$, and $\mathrm{Hom}(\mathbb{Q}_{p},D)=0$ since $\mathbb{Q}_{p}$ is
a topological $p$-group and $D$ is countable discrete and torsion-free. Thus,%
\begin{equation*}
0=\mathrm{Ext}(Y,D)\cong \mathrm{Hom}(\mathbb{Q},D)\text{.}
\end{equation*}%
Since $D$ is a direct sum of copies of $\mathbb{Q}$, this implies $D=0$.
\end{proof}

A similar proof of Theorem \ref{Theorem:no-injectives-TDLCPAb} gives the
following results.

\begin{lemma}
\label{Lemma:almost-injectives-LCPAb(p)}Suppose that $G/H\in \mathrm{LH}%
\left( \mathbf{LCPAb}\left( p\right) \right) $. The following assertions are
equivalent:

\begin{enumerate}
\item $G/H$ is injective for $\mathbf{LCPAb}\left( p\right) $;

\item $G/H\cong D$ where $D$ is a countable divisible $p$-group.
\end{enumerate}
\end{lemma}

\begin{proof}
(2)$\Rightarrow $(1) Let $X$ be a locally compact Polish topological $p$%
-group. Suppose that $D$ is a countable divisible $p$-group. We need to
prove that $\mathrm{Ext}\left( X,D\right) =0$. To this purpose, let%
\begin{equation*}
0\rightarrow D\rightarrow E\overset{\pi }{\rightarrow }X\rightarrow 0
\end{equation*}%
be an extension of $X$ by $D$. We identify $D$ with the closed subgroup $%
\mathrm{Ker}\left( \pi \right) $ of $E$. Since $D$ is a countable divisible $%
p$-group, we have that $D$ is injective in the category of discrete groups
with group homomorphisms as morphisms. Thus, there exists a (not necessarily
continuous) group homomorphism $t:X\rightarrow E$ such that $\pi t(x)=x$ for
every $x\in X$.\ Since $D$ is a closed subgroup of $E$, and $D$ is a
countable discrete group, and $E$ is a topological $p$-group, there exists
an open subgroup $U$ of $E$ such that $U\cap D=\left\{ 0\right\} $. Thus, we
have that $W:=\pi \left( U\right) \subseteq X$ is an open subgroup, and $\pi
|_{U}:U\rightarrow W$ is a topological isomorphism. Define $\sigma :=\left(
\pi |_{U}\right) ^{-1}:W\rightarrow E$. Notice that both $\sigma $ and $%
t|_{W}$ are right inverses for $\pi |_{U}$, and hence $a:=\sigma -t|_{W} $
is a homomorphism $W\rightarrow D$. By injectivity of $D$ in the category of
groups, such a homomorphism extends to a homomorphism $\tilde{a}%
:X\rightarrow D$. Define now $s:=t+\tilde{a}$. Then $s:X\rightarrow E$ is a
homomorphism that is a right inverse for $\pi $. Furthermore, $s|_{W}=\sigma
|_{W}$ is continuous. Since $W$ is an open subgroup of $X$, this implies
that $s$ is continuous. This shows that $\pi $ is a split epimorphism in the
category $\mathbf{LCPAb}\left( p\right) $, concluding the proof that $%
\mathrm{Ext}\left( X,D\right) =0$.

(1)$\Rightarrow $(2) As in the proof of Theorem \ref%
{Theorem:no-injectives-TDLCPAb}, we can assume without loss of generality
that $H$ is a countable divisible $p$-group. Since $\mathrm{Ext}\left( 
\mathbb{Z}/p\mathbb{Z},G/H\right) =0$, we infer that $G/H$ is divisible.
Furthermore, as in the proof of Theorem \ref{Theorem:no-injectives-TDLCPAb},
if $K$ is a compact open subgroup of $G$, one obtains that $\mathrm{Ext}%
\left( B,K\right) $ is countable for every compact topological $p$-group $B$%
. Thus, for every countable $p$-group $P$,%
\begin{equation*}
\mathrm{Ext}\left( K^{\vee },P\right) \cong \mathrm{Ext}\left( P^{\vee
},K\right)
\end{equation*}%
is countable, since $P^{\vee }$ is a compact topological $p$-group. Taking $%
P=\mathbb{Z}\left( p\right) $ one obtains that $K^{\vee }$ has finite $p$%
-rank by Lemma \ref{Lemma:countable-Ext(p)}. Taking 
\begin{equation*}
P:=\bigoplus_{n\in \mathbb{N}}\mathbb{Z}/p^{n}\mathbb{Z}
\end{equation*}%
one obtains that $K^{\vee }$ is reduced by Lemma \ref%
{Lemma:torsion-countable-Ext}(3). This implies that $K$ is finite, and hence 
$G$ is countable. Thus, $G/H$ is a countable divisible $p$-group.
\end{proof}

\begin{theorem}
\label{Theorem:no-injectives-LCPAb(p)}Suppose that $X\in \mathrm{LH}\left( 
\mathbf{LCPAb}\left( p\right) \right) $. The following assertions are
equivalent:

\begin{enumerate}
\item $X$ is injective in $\mathrm{LH}(\mathbf{LCPAb}\left( p\right) )$;

\item $X=0$.
\end{enumerate}
\end{theorem}

\begin{proof}
(1)$\Rightarrow $(2) Let $X$ be injective in $\mathrm{LH}(\mathbf{LCPAb}%
\left( p\right) )$. Then by\ Lemma \ref{Lemma:almost-injectives-LCPAb(p)} we
have that $X$ is a countable divisible $p$-group. Suppose by contradiction
that $X$ is nonzero. Since $X$ is a countable divisible $p$-group, it must
have a direct summand isomorphic to $\mathbb{Z}\left( p^{\infty }\right) $.
Since $X$ is injective in $\mathrm{LH}(\mathbf{LCPAb}\left( p\right) )$, the
same holds for $\mathbb{Z}\left( p^{\infty }\right) $.\ However,%
\begin{equation*}
\mathrm{Ext}\left( \Xi _{p},\mathbb{Z}\left( p^{\infty }\right) \right)
\cong \Xi _{p}\neq 0
\end{equation*}%
by Lemma \ref{Lemma:Xi-p}. This is a contradiction.
\end{proof}

\subsection{Injectives in the heart of topological torsion groups of finite
ranks}

We denote by $\mathbf{FLCPAb}\left( p\right) $ the full subcategory of $%
\mathbf{LCPAb}\left( p\right) $ of \emph{finite }$p$-\emph{rank} topological 
$p$-groups. By \cite[Lemma 2.8(iii)]{hoffmann_homological_2007}, an object $%
A $ of $\mathbf{FLCPAb}\left( p\right) $ is a finite sum of cyclic $p$%
-groups and copies of $\mathbb{Z}\left( p^{\infty }\right) $, $\mathbb{Z}%
_{p} $, and $\mathbb{Q}_{p}$. Let $F$ be a finite set of primes. For a
topological torsion group $G$, we write 
\begin{equation*}
G=G_{F}\oplus G_{F}^{\bot }
\end{equation*}%
where 
\begin{equation*}
G_{F}:=\bigoplus_{p\in F}G_{p}
\end{equation*}%
is the sum of $p$-adic components of $G$ for $p\in F$, and $G_{F}^{\bot }$
has trivial $p$-adic component for $p\in F$.

\begin{lemma}
\label{Lemma:divisible-codivisible}Let $X$ be a finite $p$-rank topological $%
p$-group.

\begin{enumerate}
\item $X$ is divisible iff it is of the form $\mathbb{Z}\left( p^{\infty
}\right) ^{k}\oplus \mathbb{Q}_{p}^{n}$ for some $k,n\in \omega $;

\item $X$ is codivisible iff it is of the form $\mathbb{Z}_{p}^{k}\oplus 
\mathbb{Q}_{p}^{n}$ for some $k,n\in \omega $; equivalently, $X$ is
torsion-free.
\end{enumerate}
\end{lemma}

\begin{proof}
Notice that finite cyclic groups are neither divisible nor codivisible. Note
also that $\mathbb{Z}\left( p^{\infty }\right) $ is divisible but $\mathbb{Z}%
\left( p^{\infty }\right) ^{\vee }\cong \mathbb{Z}_{p}$ is not divisible.
Finally, note that $\mathbb{Q}_{p}$ is divisible and%
\begin{equation*}
\mathbb{Q}_{p}^{\vee }\cong \mathrm{lim}_{n}\left( \mathbb{Z}\left(
p^{\infty }\right) ,\times p\right)
\end{equation*}%
is divisible. The conclusion follows since a finite $p$-rank topological $p$%
-group is a finite direct sum of finite cyclic groups and copies of $\mathbb{%
Z}\left( p^{\infty }\right) $, $\mathbb{Z}_{p}$, and $\mathbb{Q}_{p}$.
\end{proof}

\begin{theorem}
\label{Theorem:injectives-FLCPAb(p)}For an object $X$ of $\mathbf{FLCPAb}%
\left( p\right) $, the following assertions are equivalent:

\begin{enumerate}
\item $X$ is injective for $\mathbf{FLCPAb}\left( p\right) $;

\item $X$ is divisible.
\end{enumerate}
\end{theorem}

\begin{proof}
(2)$\Rightarrow $(1) By \cite[Proposition 4.15(v)]{hoffmann_homological_2007}%
, we have that $\mathrm{Ext}\left( A,X\right) =0$ when $A$ is codivisible,
since $X_{\mathbb{Z}}=0$. In particular, $\mathrm{Ext}\left( \mathbb{Z}%
_{p},X\right) =0$ and $\mathrm{Ext}(\mathbb{Q}_{p},X)=0$. We also have that $%
\mathrm{Ext}\left( \mathbb{Z}/p^{k},X\right) \cong X/p^{k}=0$. We have a
short exact sequence%
\begin{equation*}
\mathrm{Hom}(\mathbb{Q}_{p},X)\rightarrow \mathrm{Hom}\left( \mathbb{Z}%
_{p},X\right) \rightarrow \mathrm{Ext}\left( \mathbb{Z}\left( p^{\infty
}\right) ,X\right) \rightarrow \mathrm{Ext}(\mathbb{Q}_{p},X)=0\text{.}
\end{equation*}%
As $X$ is divisible, 
\begin{equation*}
\mathrm{Hom}(\mathbb{Q}_{p},X)\rightarrow \mathrm{Hom}\left( \mathbb{Z}%
_{p},X\right)
\end{equation*}%
is surjective by \cite[Proposition 3.3]{hoffmann_homological_2007}. Thus, $%
\mathrm{Ext}\left( \mathbb{Z}\left( p^{\infty }\right) ,X\right) =0$. This
concludes the proof by \cite[Lemma 2.8(iii)]{hoffmann_homological_2007}.

(1)$\Rightarrow $(2) If $X$ is injective in $\mathbf{FLCPAb}\left( p\right) $%
, then $\mathrm{Ext}\left( \mathbb{Z}/p\mathbb{Z},X\right) =0$. As $\mathrm{%
Ext}\left( \mathbb{Z}/p\mathbb{Z},X\right) \cong X/pX$, we obtain that $pX=X$
and $X$ is $p$-divisible and hence divisible.
\end{proof}

In a similar fashion, we have the following:

\begin{lemma}
\label{Lemma:injectives-TorFLCPAb}Let $X$ be an object of $\mathbf{TorFLCPAb}
$. Suppose that for every compact group $C$ in $\mathbf{TorFLCPAb}$, $%
\mathrm{Ext}\left( C,X\right) $ is countable. Then there exists a finite set 
$F$ of primes such that $X_{F}^{\bot }$ is countable.
\end{lemma}

\begin{proof}
Observe that for every compact open subgroup $U$ of $X$, and for every
compact group $C$ in $\mathbf{TorFLCPAb}$, $\mathrm{Ext}\left( C,U\right) $
is countable.

Set $Y:=X^{\vee }$. Choose a compact open subgroup $U\subseteq Y$. Then for
every countable torsion group $T$ with cyclic primary components, the
sequence%
\begin{equation*}
\mathrm{Hom}\left( U,T\right) \rightarrow \mathrm{Ext}\left( Y/U,T\right)
\rightarrow \mathrm{Ext}\left( Y,T\right)
\end{equation*}%
is exact. Since $U$ is a compact Polish group and $T$ is a countable
discrete group, $\mathrm{Hom}\left( U,T\right) $ is countable. Likewise, $%
\mathrm{Ext}\left( Y,T\right) \cong \mathrm{Ext}\left( T^{\vee },X\right) $
is countable. Thus, the same holds for $\mathrm{Ext}\left( Y/U,T\right) $.
Thus, by\ Lemma \ref{Lemma:torsion-countable-Ext}, there exists a finite set 
$F$ of primes such that $\left( Y/U\right) _{p}=Y_{p}/U_{p}=0$ for all
primes $p$ \emph{not }in $F$. This implies that for all primes $p$ not in $F$%
, $U_{p}=Y_{p}$, and hence%
\begin{equation*}
Y_{F}^{\bot }=\prod_{p\notin F}\left( Y_{p}:U_{p}\right) =\prod_{p\notin
F}Y_{p}
\end{equation*}%
is compact. This proves that $X_{F}^{\bot }=\left( Y_{F}^{\bot }\right)
^{\vee }$ is countable.
\end{proof}

\begin{theorem}
\label{Theorem:injectives-TorFLCPAb}For any object $X$ of $\mathbf{TorFLCPAb}
$, the following assertions are equivalent:

\begin{enumerate}
\item $X$ is injective in $\mathbf{TorFLCPAb}$;

\item $X\cong D\oplus Z$ for $D,Z$ in $\mathbf{TorFLCPAb}$, where $D$ is
countable divisible, $Z$ is both divisible and codivisible, and all but
finitely many $p$-components of $Z$ are zero.
\end{enumerate}
\end{theorem}

\begin{proof}
(1)$\Rightarrow $(2): If $X$ is an arbitrary object of $\mathbf{TorFLCPAb}$,
then by \cite[Lemma 2.8(iii)]{hoffmann_homological_2007} we have that for
every prime $p$, $X_{p}$ is a finite sum of cyclic $p$-groups and copies of $%
\mathbb{Z}\left( p^{\infty }\right) $, $\mathbb{Z}_{p}$, and $\mathbb{Q}_{p}$%
. Suppose that $X$ is injective. Then $\mathrm{Ext}\left( C,X\right) =0$ for
every $C$ in $\mathbf{TorFLCPAb}$. In particular, for a prime $p$ and $C=%
\mathbb{Z}\left( p\right) $, $X/pX\cong \mathrm{Ext}\left( \mathbb{Z}\left(
p\right) ,X\right) =0$ and $X$ is $p$-divisible. Furthermore, by Lemma \ref%
{Lemma:injectives-TorFLCPAb}, $X_{p}$ is countable for all but finitely many
primes $p$. Choose a finite set $F$ of primes such that $D_{0}:=X_{F}^{\bot
} $ is countable. Since $X$ is divisible, $D_{0}$ is divisible. The group $%
X_{F}$ has only finitely many nonzero $p$-components, so Lemma \ref%
{Lemma:divisible-codivisible}(1) yields a decomposition $X_{F}\cong
D_{1}\oplus Z$, where $D_{1}$ is countable and divisible, and $Z$ is both
divisible and codivisible and has nonzero $p$-components only for $p\in F$.
Thus $X\cong \left( D_{0}\oplus D_{1}\right) \oplus Z$, as required.

(2)$\Rightarrow $(1) From (3)$\Rightarrow $(1) in Theorem \ref%
{Theorem:LCPAb-injective-LH(TDLCPAb)}, it follows that $D$ is injective in $%
\mathbf{TorFLCPAb}$. We need to prove that $Z$ is injective in $\mathbf{%
TorFLCPAb}$. For this purpose, it suffices to prove that $\mathbb{Q}_{p}$ is
injective in $\mathbf{TorFLCPAb}$ for every prime $p$. Let $G$ be an object
of $\mathbf{TorFLCPAb}$, let $N\subseteq G$ be a closed subgroup, and let $%
f:N\rightarrow \mathbb{Q}_{p}$ be a continuous homomorphism. Decompose%
\begin{equation*}
N=N_{p}\oplus N^{\prime }\text{\quad and\quad }G=G_{p}\oplus G^{\prime },
\end{equation*}%
where $N^{\prime }$ and $G^{\prime }$ have no nontrivial $p$-components.
Then $f$ vanishes on $N^{\prime }$ and is determined by its restriction to $%
N_{p}$. The induced map $N_{p}\rightarrow G_{p}$ is a monomorphism in the
abelian category $\mathbf{FLCPAb}\left( p\right) $. By Theorem \ref%
{Theorem:injectives-FLCPAb(p)}, the restriction $N_{p}\rightarrow \mathbb{Q}%
_{p}$ extends to $G_{p}$. Composing with the projection $G\rightarrow G_{p}$
gives an extension of $f$ to $G$. Thus $\mathbb{Q}_{p}$ is injective in $%
\mathbf{TorFLCPAb}$.
\end{proof}

\begin{remark}
The groups%
\begin{equation*}
\prod_{p}\left( \mathbb{Z}\left( p^{\infty }\right) :\mathbb{Z}\left(
p\right) \right)
\end{equation*}%
and%
\begin{equation*}
\prod_{p}(\mathbb{Q}_{p}:\mathbb{Z}_{p})
\end{equation*}%
where the (reduced) product is taken over all primes, are examples of
divisible groups in $\mathbf{TorFLCPAb}$ that are not injective.
\end{remark}

\begin{corollary}
The category $\mathbf{TorFLCPAb}$ does not have enough injectives.
\end{corollary}

\begin{proof}
Consider%
\begin{equation*}
K:=\prod_{p}\mathbb{Z}\left( p\right)
\end{equation*}%
where $p$ ranges over all primes. Consider an injective object of the form $%
D\oplus Z$ where $D$ is countable divisible, and $Z$ is divisible and
codivisible and with only finitely many nonzero $p$-components. By Lemma \ref%
{Lemma:divisible-codivisible}, $Z$ is of the form $\bigoplus_{q\in F}\mathbb{%
Q}_{q}^{n_{q}}$ for a finite set $F$ of primes and, in particular, is
torsion-free. Since $K$ is a torsion group, $\mathrm{Hom}\left( K,Z\right)
=0 $. Thus,%
\begin{equation*}
\mathrm{Hom}\left( K,D\oplus Z\right) \cong \mathrm{Hom}\left( K,D\right) 
\text{.}
\end{equation*}%
Let $\varphi :K\rightarrow D$ be a continuous homomorphism. Since $D$ is
countable discrete, the kernel of $\varphi $ is an open subgroup of $K$.
Since $K$ is not discrete, $\ker \left( \varphi \right) $ is nonzero. Thus, $%
\varphi $ is not injective.
\end{proof}

We first record two lemmas that will be used to show that every injective
object in the left heart is represented by an object of $\mathbf{TorFLCPAb}$.

\begin{lemma}
\label{Lemma:dense-monomorphism-TorFLCPAb}Suppose that $H$ and $K$ are
objects of $\mathbf{TorFLCPAb}$ such that $H$ is a dense Polish subgroup of $%
K$ and $K$ is compact. There are \emph{finite} $p$-groups $A_{p}$ such that,
in $\mathrm{LH}\left( \mathbf{TorFLCPAb}\right) $, 
\begin{equation*}
K/H\cong \frac{\prod_{p}A_{p}}{\bigoplus_{p}A_{p}}.
\end{equation*}
\end{lemma}

\begin{proof}
After replacing $H$ and $K$ with a quotient by a compact open subgroup of $H$%
, we can assume that $H$ is \emph{countable}. Thus, $H=\bigoplus_{p}H_{p}$,
while $K\cong \prod_{p}K_{p}$. For every prime $p$, we have $H=H_{p}\oplus
H_{p}^{\bot }$ and $K=K_{p}\oplus K_{p}^{\bot }$. Since $H$ is dense in $K$, 
$H_{p}$ is dense in $K_{p}$. By the structure of compact finite $p$-rank
topological $p$-groups, $K_{p}\cong \mathbb{Z}_{p}^{n_{p}}\oplus F_{p}$ for
some $n_{p}\in \omega $ and some finite $p$-group $F_{p}$. Since $H_{p}$ is
torsion, its projection onto $\mathbb{Z}_{p}^{n_{p}}$ is trivial. Density
therefore forces $n_{p}=0$, and hence $K_{p}$ is finite and $H_{p}=K_{p}$.
Setting $A_{p}:=K_{p}$ for every prime $p$ concludes the proof.
\end{proof}

The following lemma in particular shows that \emph{reduced products }(with
respect to the Fr\'{e}chet filter) of finite $p$-groups, indexed over all
primes, regarded as objects of $\mathrm{LH}\left( \mathbf{TorFLCPAb}\right) $%
, are not injective.

\begin{lemma}
\label{Lemma:reduced-product-obstruction}Let $S$ be an infinite set of
primes and let $A_{p}$ be a nonzero finite $p$-group for every $p\in S$. Set 
\begin{equation*}
P:=\prod_{p\in S}A_{p},\qquad P_{0}:=\bigoplus_{p\in S}A_{p}\text{.}
\end{equation*}%
Then there exists an object $B$ of $\mathbf{TorFLCPAb}$ and a continuous
homomorphism $j:P\rightarrow B$ with closed image such that $\mathrm{Ext}%
\left( B/j\left( P\right) ,P/P_{0}\right) $ is uncountable.
\end{lemma}

\begin{proof}
Let $\left( S^{\left( n\right) }\right) $ be a partition of $S$ into
infinitely many infinite sets of primes. For every prime $p$ we can write $%
A_{p}$ as the direct sum of cyclic $p$-groups $A_{p}\left( i\right) $ of
order $p^{a_{p}\left( i\right) }$ for $i<r_{p}\in \mathbb{N}$ with 
\begin{equation*}
a_{p}\left( 0\right) \leq a_{p}\left( 1\right) \leq \cdots \leq a_{p}\left(
r_{p}-1\right) :=e_{p}\text{.}
\end{equation*}%
Define now $B_{p}$ to be the direct sum of cyclic $p$-groups $B_{p}\left(
i\right) $ of order $p^{a_{p}\left( i\right) +e_{p}}$ for $i<r_{p}$. Define $%
j_{p}:A_{p}\rightarrow B_{p}$, 
\begin{equation*}
\left( x_{i}\right) _{i<r_{p}}\mapsto \left( p^{e_{p}}x_{i}\right) _{i<r_{p}}%
\text{.}
\end{equation*}%
Thus $j_{p}$ is injective. Furthermore, every homomorphism $%
u:B_{p}\rightarrow A_{p}$ satisfies $u\circ j_{p}=0$, because $%
p^{e_{p}}A_{p}=0$. Set 
\begin{equation*}
B:=\prod_{p\in S}B_{p},\qquad j:=\prod_{p\in S}j_{p}:P\longrightarrow B.
\end{equation*}%
Both $P$ and $B$ are objects of $\mathbf{TorFLCPAb}$, are compact, and $j$
has closed image.

Let $q:P\rightarrow P/P_{0}$ be the quotient map. Let us initially show that
there exists no Borel-definable homomorphism $\Phi :B\rightarrow P/P_{0}$
such that $\Phi \circ j=q$. Suppose that $\Phi $ is such a morphism. Since $%
P_{0}$ is countable and discrete, Proposition \ref{Proposition:morphisms}
provides an open subgroup $V$ in $B$ and a Borel lift $\varphi :B\rightarrow
P$ of $\Phi $ such that $\varphi |_{V}$ is continuous and additive. There is
a finite subset $E$ of $S$ such that 
\begin{equation*}
B_{E}:=\prod_{p\in S\setminus E}B_{p}\subseteq V,
\end{equation*}%
where the omitted coordinates are set equal to zero. The restriction 
\begin{equation*}
h:=\varphi |_{B_{E}}:B_{E}\rightarrow P
\end{equation*}%
is a continuous homomorphism.

Every coordinate homomorphism $h_{p}:B_{E}\rightarrow A_{p}$ factors through
a finite product of the groups $B_{q}$ for $q\in S\setminus E$. If $p\in E$,
then all such primes $q$ are distinct from $p$, and hence $h_{p}=0$. If $%
p\in S\setminus E$, then homomorphisms between a $q$-group and a $p$-group
vanish for $q\neq p$, while the restriction of $h_{p}$ to $B_{p}$ vanishes
on $j_{p}(A_{p})$.

Choose $x=(x_{p})\in P$ such that $x_{p}\neq 0$ for $p\in S\setminus E$ and $%
x_{p}=0$ for $p\in E$. Thus, $j\left( x\right) \in B_{E}$. It follows that 
\begin{equation*}
h(j(x))=0\text{.}
\end{equation*}%
On the other hand, since $\varphi $ lifts $\Phi $ and $\Phi \circ j=q$, we
have 
\begin{equation*}
q(\varphi (j(x)))=\Phi (j(x))=q(x).
\end{equation*}%
Thus $q(x)=0$, contradicting the fact that $x\notin P_{0}$.

Let now $\left( S^{\left( n\right) }\right) $ be a partition of $S$ into
infinitely many infinite sets. Fix $n\in \mathbb{N}$. Let $P^{\left(
n\right) }$, $P_{0}^{\left( n\right) }$, $B^{\left( n\right) }$, $j^{\left(
n\right) }$, and $q^{\left( n\right) }$ be defined as above, replacing $S$
with $S^{\left( n\right) }$. We identify $P$ with the product of $P^{\left(
n\right) }$ for $n\in \mathbb{N}$. Let $\alpha \in 2^{\mathbb{N}}$ where $%
2=\left\{ 0,1\right\} $. Define%
\begin{equation*}
q_{\alpha }:P\rightarrow P/P_{0}\text{, }\left( x^{\left( n\right) }\right)
_{n\in \mathbb{N}}\mapsto \left( \alpha \left( n\right) x^{\left( n\right)
}\right) _{n\in \mathbb{N}}+P_{0}\in P/P_{0}\text{.}
\end{equation*}%
Then, by applying the remarks above after replacing $S$ with $S^{\left(
n\right) } $ for any $n\in \mathbb{N}$, we obtain that, for any two distinct 
$\alpha ,\beta \in 2^{\mathbb{N}}$, $q_{\alpha }-q_{\beta }$ does not extend
alongside $j$.

Consider the exact sequence%
\begin{equation*}
\mathrm{Hom}\left( B,P/P_{0}\right) \rightarrow \mathrm{Hom}\left(
P,P/P_{0}\right) \rightarrow \mathrm{Ext}\left( B/j\left( P\right)
,P/P_{0}\right) \text{.}
\end{equation*}%
Consider also the function%
\begin{equation*}
2^{\mathbb{N}}\rightarrow \mathrm{Hom}\left( P,P/P_{0}\right) \text{, }%
\alpha \mapsto q_{\alpha }\text{.}
\end{equation*}%
The remarks above show that the composed function%
\begin{equation*}
2^{\mathbb{N}}\rightarrow \mathrm{Hom}\left( P,P/P_{0}\right) \rightarrow 
\mathrm{Ext}\left( B/j\left( P\right) ,P/P_{0}\right)
\end{equation*}%
is injective. Since $2^{\mathbb{N}}$ is uncountable, the same holds for $%
\mathrm{Ext}\left( B/j\left( P\right) ,P/P_{0}\right) $.
\end{proof}

\begin{lemma}
\label{Lemma:heart-injective-is-topological}Suppose that $X=G/N$ is an
object of $\mathrm{LH}\left( \mathbf{TorFLCPAb}\right) $. Assume that $%
\mathrm{Ext}\left( C,G/N\right) $ is \emph{countable }for every compact
group $C$ in $\mathbf{TorFLCPAb}$. Then $X$ is isomorphic in $\mathrm{LH}%
\left( \mathbf{TorFLCPAb}\right) $ to an object of $\mathbf{TorFLCPAb}$.
\end{lemma}

\begin{proof}
We need to prove that $N$ is closed in $G$. This is equivalent to the
assertion that $N\cap U$ is closed in $U$ for some (or, equivalently, every)
open subgroup $U$ of $G$. Likewise, since $\mathbf{FLCPAb}\left( p\right) $
is an abelian category, it is equivalent to the assertion that $N_{F}^{\bot
} $ is closed in $G_{F}^{\bot }$ for some (equivalently, every) finite set
of primes $F$.

Suppose that $U\subseteq G$ is a compact open subgroup. Then the short exact
sequence%
\begin{equation*}
0\rightarrow U/\left( U\cap N\right) \rightarrow G/N\rightarrow
G/U\rightarrow 0
\end{equation*}%
induces for every compact object $C$ of $\mathbf{TorFLCPAb}$ an exact
sequence%
\begin{equation*}
\mathrm{Hom}\left( C,G/U\right) \rightarrow \mathrm{Ext}\left( C,U/\left(
U\cap N\right) \right) \rightarrow \mathrm{Ext}\left( C,G/N\right)
\end{equation*}%
Since $G/U$ is countable discrete, every homomorphism $C\rightarrow G/U$
factors through a discrete (hence, finite) quotient $F$ of $C$.\ Since for
every such quotient $F$, $\mathrm{Hom}\left( F,G/U\right) $ is countable,
and there are countably many such finite quotients, we conclude that $%
\mathrm{Hom}\left( C,G/U\right) $ is countable. Hence $\mathrm{Ext}\left(
C,U/\left( U\cap N\right) \right) $ is also countable. Thus, $U/\left( U\cap
N\right) $ satisfies the same hypothesis as $G/N$.

There is a short exact sequence in the left heart 
\begin{equation*}
0\longrightarrow \overline{N}^{\,G}/N\longrightarrow G/N\overset{\pi }{%
\longrightarrow }G/\overline{N}^{\,G}\longrightarrow 0.
\end{equation*}

We observe that $G/\overline{N}^{\,G}$ satisfies the hypotheses of Lemma \ref%
{Lemma:injectives-TorFLCPAb}. Indeed, for a compact group $T\in \mathbf{%
TorFLCPAb}$, we have a corresponding exact sequence 
\begin{equation*}
\mathrm{Ext}(T,G/N)\longrightarrow \mathrm{Ext}(T,G/\overline{N}%
^{\,G})\longrightarrow \mathrm{Ext}^{2}(T,\overline{N}^{\,G})=0\text{.}
\end{equation*}%
The first term is countable by hypothesis. The last term is zero by Lemma %
\ref{Lemma:vanishing-Ext2}. Hence $\mathrm{Ext}(T,G/\overline{N}^{\,G})$ is
countable for every $T$, as claimed. Lemma \ref{Lemma:injectives-TorFLCPAb}
therefore gives a finite set $F$ of primes such that $(G/\overline{N}%
^{\,G})_{F}^{\bot }$ is countable. After replacing $G$ with $G_{F}^{\bot }$,
we can assume without loss of generality that $G_{F}=0$. Thus, $(G/\overline{%
N}^{\,G})_{F}=0$ and $G/\overline{N}^{G}$ is countable. This shows that $%
\overline{N}^{G}$ is an open subgroup of $G$. Thus, after replacing $G$ with 
$\overline{N}^{G}$ we can assume that $N$ is \emph{dense }in $G$.
Furthermore, after replacing $G$ with a compact open subgroup $U$ and $N$
with $N\cap U$, we can assume that $G$ is compact. It then follows that $G/N$
is trivial by Lemma \ref{Lemma:reduced-product-obstruction} and Lemma \ref%
{Lemma:dense-monomorphism-TorFLCPAb}.
\end{proof}

We can now characterize all injective objects in the left heart of $\mathbf{%
TorFLCPAb}$.

\begin{theorem}
\label{Theorem:injectives-LH(TorFLCPAb)}For any object $X$ of $\mathrm{LH}(%
\mathbf{TorFLCPAb})$, the following assertions are equivalent:

\begin{enumerate}
\item $X$ is injective in $\mathrm{LH}\left(\mathbf{TorFLCPAb}\right)$;

\item $X$ is isomorphic to an object of $\mathbf{TorFLCPAb}$ that is
divisible and whose $p$-components are zero for all but finitely many primes 
$p$;

\item $X$ is isomorphic to a finite direct sum of copies of $\mathbb{Z}%
\left( p^{\infty }\right) $ and $\mathbb{Q}_{p}$, for primes $p$.
\end{enumerate}
\end{theorem}

\begin{proof}
The equivalence of (2) and (3) follows from the structure of finite-rank
topological $p$-groups.

(3)$\Rightarrow $(1) It suffices to consider a divisible finite-rank
topological $p$-group $X$ for a fixed prime $p$. Let $i:N\rightarrow G$ be
an arbitrary monomorphism in $\mathbf{TorFLCPAb}$, not necessarily strict,
and let $f:N\rightarrow X$ be a continuous homomorphism. We have topological
direct-sum decompositions%
\begin{equation*}
N=N_{p}\oplus N^{\prime }\text{\quad and\quad }G=G_{p}\oplus G^{\prime },
\end{equation*}%
where $N^{\prime }$ and $G^{\prime }$ have no nontrivial $p$-components. The
homomorphism $i$ restricts to a monomorphism $i_{p}:N_{p}\rightarrow G_{p}$
in the abelian category $\mathbf{FLCPAb}\left( p\right) $. Moreover, $%
\mathrm{Hom}\left( N^{\prime },X\right) =0$, so $f$ is determined by its
restriction $f_{p}:N_{p}\rightarrow X$. By Theorem \ref%
{Theorem:injectives-FLCPAb(p)}, there is a continuous homomorphism $%
g_{p}:G_{p}\rightarrow X$ such that $g_{p}i_{p}=f_{p}$. Composing $g_{p}$
with the projection $G\rightarrow G_{p}$ gives an extension $g:G\rightarrow
X $ of $f$.

We now verify injectivity in the left heart. Let $Y=G/N$ be an arbitrary
object of $\mathrm{LH}\left( \mathbf{TorFLCPAb}\right) $. The exact sequence%
\begin{equation*}
\mathrm{Hom}\left( G,X\right) \rightarrow \mathrm{Hom}\left( N,X\right)
\rightarrow \mathrm{Ext}\left( Y,X\right) \rightarrow \mathrm{Ext}\left(
G,X\right)
\end{equation*}%
has a surjective first map by the preceding paragraph. The same extension
property for strict monomorphisms shows that $X$ is injective in $\mathbf{%
TorFLCPAb}$, and hence $\mathrm{Ext}\left( G,X\right) =0$. Therefore $%
\mathrm{Ext}\left( Y,X\right) =0$. Since $Y$ was arbitrary, $X$ is injective
in $\mathrm{LH}\left( \mathbf{TorFLCPAb}\right) $. The conclusion for the
groups in (3) follows because finite direct sums of injective objects are
injective.

(1)$\Rightarrow $(2) Suppose that $X$ is injective in $\mathrm{LH}(\mathbf{%
TorFLCPAb})$. By Lemma \ref{Lemma:heart-injective-is-topological}, $X$ is
isomorphic to an object of $\mathbf{TorFLCPAb}$. It is injective in $\mathbf{%
TorFLCPAb}$ as well, since every strict monomorphism in $\mathbf{TorFLCPAb}$
is a monomorphism in its left heart. Hence, by Theorem \ref%
{Theorem:injectives-TorFLCPAb}, $X$ is isomorphic to $D\oplus Z$ where 
\begin{equation*}
D=\bigoplus_{p}\mathbb{Z}\left( p^{\infty }\right) ^{d_{p}}\text{\quad
and\quad }Z=\bigoplus_{p\in F}\mathbb{Q}_{p}^{n_{p}}\text{,}
\end{equation*}%
for some finite set $F$ of primes and $d_{p},n_{p}\in \mathbb{N}$. Suppose,
towards a contradiction, that the set $\left\{ p:d_{p}>0\right\} $ is
infinite, and let $S$ be an infinite subset of it. Choosing one copy of $%
\mathbb{Z}\left( p^{\infty }\right) $ for every $p\in S$, we obtain a direct
summand of $X$ of the form%
\begin{equation*}
D_{S}:=\bigoplus_{p\in S}\mathbb{Z}\left( p^{\infty }\right)
\end{equation*}%
Since a direct summand of an injective object is injective, $D_{S}$ would be
injective in $\mathrm{LH}(\mathbf{TorFLCPAb})$. Consider%
\begin{equation*}
T_{S}:=\bigoplus_{p\in S}\mathbb{Z}\left( p\right)
\end{equation*}%
and%
\begin{equation*}
C_{S}:=\prod_{p\in S}\mathbb{Z}\left( p\right) \text{.}
\end{equation*}%
Then $T_{S}$ is a subgroup of $C_{S}$. For $p\in S$, identify $\mathbb{Z}%
\left( p\right) $ as a subgroup of $\mathbb{Z}\left( p^{\infty }\right) $.
This also allows one to regard $T_{S}$ as a subgroup of $D_{S}$. If $D_{S}$
were injective, the inclusion $T_{S}\rightarrow D_{S}$ would extend across
the dense monomorphism $T_{S}\rightarrow C_{S}$ to a continuous homomorphism 
$\varphi :C_{S}\rightarrow D_{S}$. Since $D_{S}$ is discrete, its kernel $U$
is an open subgroup of $C_{S}$. Since $C_{S}$ is not discrete, $U$ is
nonzero; and since $T_{S}$ is dense in $C_{S}$, the intersection $U\cap
T_{S} $ is nonzero. This contradicts the fact that $\varphi |_{T_{S}}$ is
the inclusion $T_{S}\rightarrow D_{S}$. Therefore $d_{p}=0$ for all but
finitely many primes $p$.
\end{proof}

\subsection{Projectives\label{Subsection:projectives}}

We conclude with a characterization of projective objects in the left hearts
of the categories we have considered so far. The following corollary is a
particular case of a general fact about quasi-abelian categories \cite[%
Proposition 1.3.24]{schneiders_quasi-abelian_1999}.

\begin{corollary}
\label{Corollary:projectives}Suppose that $\mathcal{A}$ is a thick
quasi-abelian subcategory of the category of locally compact Polish abelian
groups. An object $P\in \mathrm{LH}\left( \mathcal{A}\right) $ is projective
in $\mathrm{LH}\left( \mathcal{A}\right) $ if and only if $P$ is isomorphic
to an object of $\mathcal{A}$ that is projective in $\mathcal{A}$.
\end{corollary}

\begin{proof}
If $P\in \mathcal{A}$ is projective in $\mathcal{A}$, then it is projective
in $\mathrm{LH}\left( \mathcal{A}\right) $ by \cite[Proposition 1.3.24(a)]%
{schneiders_quasi-abelian_1999}. Conversely, if $P\in \mathrm{LH}\left( 
\mathcal{A}\right) $ is projective in $\mathrm{LH}\left( \mathcal{A}\right) $%
, then the proof of \cite[Proposition 1.3.24(b)]%
{schneiders_quasi-abelian_1999} shows that $P$ is isomorphic to an object of 
$\mathcal{A}$. The assumption that $\mathcal{A}$ has enough projective
objects from the statement of \cite[Proposition 1.3.24(b)]%
{schneiders_quasi-abelian_1999} is not used in the proof of this fact.
\end{proof}

Using Corollary \ref{Corollary:projectives} and the results about injectives
we have presented so far, we obtain a characterization of projective objects.

\begin{proposition}
We have the following:

\begin{enumerate}
\item The categories $\mathrm{LH}\left( \mathbf{LCPAb}\right) $ and $\mathrm{%
LH}\left( \mathbf{LieAb}\right) $ have the same projective objects, which
are the Polish abelian groups of the form $V\oplus F$ where $V$ is a vector
group and $F$ is a countable free abelian group;

\item The projective objects in $\mathrm{LH}\left( \mathbf{TDLCPAb}\right) $
are the countable free abelian groups;

\item The projective objects in $\mathrm{LH}\left( \mathbf{TorLCPAb}\right) $
are the codivisible compact topological torsion groups;

\item The projective objects in $\mathrm{LH}\left( \mathbf{LCPAb}\left(
p\right) \right) $ are the codivisible compact topological $p$-groups;

\item The projective objects in $\mathbf{FLCPAb}\left( p\right) $ are
codivisible finite $p$-rank topological $p$-groups;

\item The projective objects in $\mathrm{LH}\left( \mathbf{TorFLCPAb}\right) 
$ are precisely those of the form $C\oplus W$ for some $C$ and $W$ in $%
\mathbf{TorFLCPAb}$ such that $C$ is compact codivisible and $W$ is
divisible and codivisible and all but finitely many $p$-components of $W$
are zero.
\end{enumerate}
\end{proposition}

\begin{proof}
We will appeal to Corollary \ref{Corollary:projectives} for all the items of
the proof.

(1) This follows from Proposition \ref{Proposition:injectives}.

(2) This follows from \cite[Theorem 4.1]{fulp_extensions_1971}.

(3) This follows from Proposition \ref{Proposition:injectives-TDLC}.

(4) This follows from (3).

(5) This follows from Theorem \ref{Theorem:injectives-FLCPAb(p)}.

(6) This follows from Theorem \ref{Theorem:injectives-TorFLCPAb}.
\end{proof}

\bibliographystyle{amsalpha}
\bibliography{bibliography}

\end{document}